\input amstex
\documentstyle{amsppt}
\newcount\mgnf\newcount\tipi\newcount\tipoformule\newcount\greco
\tipi=2          
\tipoformule=0   

\global\newcount\numsec\global\newcount\numfor
\global\newcount\numapp\global\newcount\numcap
\global\newcount\numfig\global\newcount\numpag
\global\newcount\numnf
\global\newcount\numtheo

\def\SIA #1,#2,#3 {\senondefinito{#1#2}%
\expandafter\xdef\csname #1#2\endcsname{#3}\else
\write16{???? ma #1,#2 e' gia' stato definito !!!!} \fi}

\def \FU(#1)#2{\SIA fu,#1,#2 }

\def\etichetta(#1){(\veroparagrafo.\veraformula)%
\SIA e,#1,(\veroparagrafo.\veraformula) %
\global\advance\numfor by 1%
\write15{\string\FU (#1){\equ(#1)}}%
\write16{ EQ #1 ==> \equ(#1)  }}

\def\etichettat(#1){\veroparagrafo.\veratheorema:%
\SIA e,#1,{\veroparagrafo.\veratheorema} %
\global\advance\numtheo by 1%
\write15{\string\FU (#1){\thu(#1)}}%
\write16{ TH #1 ==> \thu(#1)  }}

\def\etichettaa(#1){(A\veraappendice.\veraformula)
 \SIA e,#1,(A\veraappendice.\veraformula)
 \global\advance\numfor by 1
 \write15{\string\FU (#1){\equ(#1)}}
 \write16{ EQ #1 ==> \equ(#1) }}
\def\getichetta(#1){Fig. \verafigura
 \SIA g,#1,{\verafigura}
 \global\advance\numfig by 1
 \write15{\string\FU (#1){\graf(#1)}}
 \write16{ Fig. #1 ==> \graf(#1) }}
\def\retichetta(#1){\numpag=\pgn\SIA r,#1,{\verapagina}
 \write15{\string\FU (#1){\rif(#1)}}
 \write16{\rif(#1) ha simbolo  #1  }}
\def\etichettan(#1){(n\verocapitolo.\veranformula)
 \SIA e,#1,(n\verocapitolo.\veranformula)
 \global\advance\numnf by 1
\write16{\equ(#1) <= #1  }}

\newdimen\gwidth
\gdef\profonditastruttura{\dp\strutbox}
\def\senondefinito#1{\expandafter\ifx\csname#1\endcsname\relax}
\def\BOZZA{
\def\alato(##1){
 {\vtop to \profonditastruttura{\baselineskip
 \profonditastruttura\vss
 \rlap{\kern-\hsize\kern-1.2truecm{$\scriptstyle##1$}}}}}
\def\galato(##1){ \gwidth=\hsize \divide\gwidth by 2
 {\vtop to \profonditastruttura{\baselineskip
 \profonditastruttura\vss
 \rlap{\kern-\gwidth\kern-1.2truecm{$\scriptstyle##1$}}}}}
\def\verapagina{
{\romannumeral\number\numcap}.\number\numsec.\number\numpag}}

\def\alato(#1){}
\def\galato(#1){}
\def\veroparagrafo{\number\numsec}\def\veraformula{\number\numfor}
\def\veraappendice{\number\numapp}
\def\verapagina{\number\pageno}\def\veranformula{\number\numnf}
\def\verafigura{{\romannumeral\number\numcap}.\number\numfig}
\def\verocapitolo{\number\numcap}\def\veranformula{\number\numnf}
\def\veratheorema{\number\numtheo}
\def\Eqn(#1){\eqno{\etichettan(#1)\alato(#1)}}
\def\eqn(#1){\etichettan(#1)\alato(#1)}
\def\TH(#1){{\etichettat(#1)\alato(#1)}}
\def\thv(#1){\senondefinito{fu#1}$\clubsuit$#1\else\csname fu#1\endcsname\fi}
\def\thu(#1){\senondefinito{e#1}\thv(#1)\else\csname e#1\endcsname\fi}

\def\Eq(#1){\eqno{\etichetta(#1)\alato(#1)}}
\def\eq(#1){\etichetta(#1)\alato(#1)}
\def\Eqa(#1){\eqno{\etichettaa(#1)\alato(#1)}}
\def\eqa(#1){\etichettaa(#1)\alato(#1)}
\def\dgraf(#1){\getichetta(#1)\galato(#1)}
\def\drif(#1){\retichetta(#1)}

\def\eqv(#1){\senondefinito{fu#1}$\clubsuit$#1\else\csname fu#1\endcsname\fi}
\def\equ(#1){\senondefinito{e#1}\eqv(#1)\else\csname e#1\endcsname\fi}
\def\graf(#1){\senondefinito{g#1}\eqv(#1)\else\csname g#1\endcsname\fi}
\def\rif(#1){\senondefinito{r#1}\eqv(#1)\else\csname r#1\endcsname\fi}
\def\bib[#1]{[#1]\numpag=\pgn
\write13{\string[#1],\verapagina}}

\def\include#1{
\openin13=#1.aux \ifeof13 \relax \else
\input #1.aux \closein13 \fi}

\openin14=\jobname.aux \ifeof14 \relax \else
\input \jobname.aux \closein14 \fi
\openout15=\jobname.aux
\openout13=\jobname.bib

\let\EQ=\Eq

\ifnum\tipoformule=1\let\Eq=\eqno\def\eq{}\let\Eqa=\eqno\def\eqa{}
\def\equ{}\fi


{\count255=\time\divide\count255 by 60 \xdef\hourmin{\number\count255}
        \multiply\count255 by-60\advance\count255 by\time
   \xdef\hourmin{\hourmin:\ifnum\count255<10 0\fi\the\count255}}

\def\oramin{\hourmin }

\def\data{\number\day/\ifcase\month\or january \or february \or march \or
april \or may \or june \or july \or august \or september
\or october \or november \or december \fi/\number\year;\ \oramin}

\def\titdate{ \ifcase\month\or January \or February \or March \or
April \or May \or June \or July \or August \or September
\or October \or November \or December \fi \number\day, \number\year;\ \oramin}

\def\titdatebis{ \ifcase\month\or January \or February \or March \or
April \or May \or June \or July \or August \or September
\or October \or November \or December \fi \number\day, \number\year}


\newcount\pgn \pgn=1
\def\foglio{\number\numsec:\number\pgn
\global\advance\pgn by 1}
\def\foglioa{A\number\numsec:\number\pgn
\global\advance\pgn by 1}

\footline={\rlap{\hbox{\copy200}}\hss\tenrm\folio\hss}


\global\newcount\numpunt

\magnification=\magstephalf
\baselineskip=12pt
\parskip=6pt

\hoffset=1.0truepc
\hsize=6.1truein
\vsize=8.4truein 

\def\a{\alpha}
\predefine\barunder{\b}
\redefine\b{\beta}
\def\d{\delta}
\def\e{\epsilon}

\def\f{\phi}
\def\g{\gamma}

\def\l{\lambda}

\def\s{\sigma}
\def\t{\tau}

\def\vp{\varphi}

\def\o{\omega}

\def\L{\Lambda}
\def\G{\Gamma}
\def\O{\Omega}

\def\[{{[\kern-.18em{[}}}
\def\]{{]\kern-.18em{]}}}

\def\1{{1\kern-.25em\roman{I}}}
\def\eu{{1\kern-.25em\roman{I}}}
\def\f1{{1\kern-.25em\roman{I}}}

\def\R{{\Bbb R}}  
\def\N{{\Bbb N}}  
\def\P{{\Bbb P}}  
\def\Z{{\Bbb Z}}  
\def\E{{\Bbb E}}  




\let\cal=\Cal

\def\BB{{\cal B}}
\def\CC{{\cal C}}
\def\DD{{\cal D}}
\def\EE{{\cal E}}
\def\FF{{\cal F}}
\def\GG{{\cal G}}

\def\LL{{\cal L}}

\def\PP{{\cal P}}

\def\SS{{\cal S}}

\def\VV{{\cal V}}

\def\VV{{\cal V}}

\def\chap #1#2{\line{\ch #1\hfill}\numsec=#2\numfor=1\numtheo=1}

\def\wt{\widetilde}
\def\wh{\widehat}


\def\note#1{\footnote{#1}}

\def\frac#1#2{{#1\over #2}}
\def\sfrac#1#2{{\textstyle{#1\over #2}}}

\def\text#1{\quad{\hbox{#1}}\quad}
\def\newpage{\vfill\eject}
\def\proposition #1{\noindent{\thbf Proposition #1}}

\def\theo #1{\noindent{\thbf Theorem {#1} }}

\def\lemma #1{\noindent{\thbf Lemma {#1} }}
\def\definition #1{\noindent{\thbf Definition {#1} }}

\def\corollary #1{\noindent{\thbf Corollary #1 }}
\def\proof{{\noindent\pr Proof: }}
\def\proofof #1{{\noindent\pr Proof of #1: }}
\def\endproof{\hfill$\square$}
\def\remark{\noindent{\bf Remark: }}
\def\thanks{\noindent{\bf Acknowledgements: }}

\font\pr=cmbxsl10

\font\thbf=cmbxsl10 scaled\magstephalf

\font\ch=cmbx12

\font\it=cmti10
\font\bf=cmbx10

\def\PRM{\hbox{\rm PRM}}
\def\asl{\hbox{\rm Asl}}



\overfullrule=0 pt


\font\tit=cmbx12
\font\aut=cmbx12

\def\s{\char'31}

\centerline{\tit AGING IN REVERSIBLE DYNAMICS OF DISORDERED SYSTEMS.}
\smallskip
\centerline{\tit I. emergence of the arcsine law in Bouchaud's asymmetric }
\centerline{\tit trap model on the complete graph}
\vskip.2truecm
\vskip1truecm
\vskip.5cm

\centerline{\aut V\'eronique  Gayrard
\note{
 CMI, LATP, Universit\'e de Provence
 39 rue F. Joliot-Curie, 13453 Marseille cedex 13\hfill\break
e-mail: gayrard\@cmi.univ-mrs.fr}
}

\vskip0.5cm
\centerline{\titdatebis}

\vskip.5cm

\vskip0.5truecm\rm
\def\s{\sigma}
\noindent {\bf Abstract:} In this paper the celebrated arcsine aging scheme of G\. Ben Arous and  J\. \v Cern\'y is taken up.
Using a brand new approach based on point processes and weak convergence techniques,
this scheme is implemented in a wide class of Markov processes
that can best be described as Glauber dynamics of discrete disordered systems.
More specifically, conditions are given for the underlying clock process (a partial sum process that
measures the total time elapsed along paths of a given length) to converge to a subordinator,
and this subordinator  is constructed explicitly.
This approach is illustrated on Bouchaud's asymmetric trap model on the
complete graph for which aging is for the first time proved, and the full, optimal picture,  obtained.

\noindent {\it Keywords:} Aging, trap model, subordinators.

\noindent {\it AMS Subject  Classification:}
82C44, 
60K35, 
82D30, 
60F17. 
\vfill
$ {} $

\newpage



\chap{1. Introduction}1

This paper is made of two parts. In a first, abstract one,
we place ourselves in the general framework of Glauber dynamics of discrete disordered systems and give
sufficient conditions for a so-called {\it arcsine aging regime} to occur.
In a second, applied one, we use these results to study  a
specific model, namely, Bouchaud's asymmetric trap model on the complete graph,
for which aging is for the first time proved, and the full, optimal picture,  obtained.
These also are the first aging results for trap model of mean field type which
is not a time change of a simple random walk.

To motivate
our goals let us first introduce a popular class of dynamics that often comes under the name of {\it ``trap models''}.
These are sequences of Markov jump processes, $X_n$,
that evolve in random landscapes made of traps, and depend on a parameter $0\leq a<1$.
To define them we first choose a graph $G_n(\VV_n, \EE_n)$ with set of vertices $\VV_n$ and set of edges $\EE_n$.
To each vertex $x\in\VV_n$ we attach a positive random variable, $\t_n(x)$, that represents the depth of a trap at $x$.
Then, given $0\leq a<1$,
$X_n$ behaves as follows: when it is at site $x$ it waits there an exponential time of parameter (proportional to)
$$
\l_n(x)=
(\t_n(x))^{-(1-a)}\sum_{y:(x,y)\in\EE_n}\t_n^{a}(y)\,,\quad\forall x\in\VV_n\,,
\Eq(3.1.2)
$$
and when it jumps, it chooses the next site, $y$, with probability
$$
p_n(x,y)=\frac{\t_n^{a}(y)}{\sum_{y:(x,y)\in\EE_n}\t_n^{a}(y)}\,,\text{if} (x,y)\in\EE_n\,,
\Eq(3.1.3)
$$
and $p_n(x,y)=0$ otherwise.  Let $J_n$ be the discrete time Markov chain
with transition kernel \eqv(3.1.3). We see that
when $a=0$, $J_n$ simply
is the homogeneous random walk on $G_n$, whereas when  $a>0$, $J_n$ favors
jumps to the neighboring traps of largest depths.
Models with $a>0$ will be called {\it asymmetric} as opposed to the {\it symmetric} ones where $a=0$.

Trap models have played a special role in the understanding of the aging phenomenon.
In fact almost all models for which the existence of an (arcsine) aging
regime has so far been proved belongs to this class of processes.
Let us review the key results
(exhaustive reviews can be found in the recent works \cite{BBC,BC4}, as well as \cite{BC2}).
Historically, symmetric trap models were introduced by Bouchaud {\it et al\.} \cite{B,BD}
as simple phenomenological models for the aging behavior of mean field spin glasses
(see \cite{BBG3,BBC2} for more on their derivation).
Taking $G_n$ as the complete graph with $n$ vertices, and letting $(\t_n(x), x\in\VV_n)$ be
i.i.d\. heavy tailed r.v\.'s, yields a model for the aging dynamics of the REM.
It is proved in \cite{BD} (see also \cite{BF,BBG2}) that this model exhibits an {\it arcsine aging regime} in the
sense that the probability $\Pi_n(t,s)$ that no jump occurs in the time interval $(t, t+s)$
behaves, for large $n$ and large times $t,s$, as the generalized arcsine distribution function
evaluated at $t/s$.
Going back to the REM itself let $\EE_n$ be the set of edges of the hypercube  $\VV_n=\{-1,1\}^n$,
and let $(\t_n(x), x\in\VV_n)$ be the Boltzman weights of the REM.
With theses choices, \eqv(3.1.2)-\eqv(3.1.3) define a special Glauber dynamics of the REM
which is known, in the symmetric case, as {\it random hopping time dynamics} (hereafter RHT dynamics).
The first connection between the REM dynamics and its trap version
was established in \cite{BBG1,BBG2} where it is proved that a discrete time version of the
RHT dynamics of the REM has the same arcsine aging regime as Bouchaud's symmetric trap model
on the complete graph.
Meanwhile, in another direction of research,
symmetric trap models on $\Z^d$ with i.i.d\. heavy tailed landscapes where studied in depth
\cite{BC2,BC3,BC4,BCM,FIN}. From this it emerged that aging in $\Z^d$, $d\geq 2$,
is the same as in Bouchaud's symmetric trap model on the complete graph.
In a landmark paper \cite{BC4}, G\. Ben Arous and  J\. \v Cern\'y proposed a scheme that
explains this apparent universality by linking the existence of an arcsine aging regime
to the arcsine law for subordinators.

This scheme centers on a certain partial sum process  $S_n$, called the {\it clock process}, that
measures the total time elapsed along paths of a given length. Namely, given two scaling sequences,
$a_n$ and $c_n$, set
$$
S_n(t)=c_n^{-1}\sum_{i=1}^{\lfloor a_n t\rfloor}\lambda_n^{-1}(J_n(i))e_{n,i}\,,\quad t>0,
\Eq(1.clock)
$$
where $(e_{n,i}\,,n\in\N, i\in\N)$ is family of independent mean one exponential
random variables, independent of $J_n$. The idea now is that if the clock process converges
to a subordinator, and if this subordinator satisfies the regular variation conditions
of Dynkin and Lamperti arcsine law, then the probability that the range of $S_n$ intersects
the time interval $(t, t+s)$  converges (in a sense to be made precise) to the generalized
arcsine distribution function evaluated at $t/s$. Now this is the signature of arcsine aging.
To put this scheme to practice one has to face two difficulties: the clock process is a random
process on the probability space of the random landscape and, for fixed realization
of the landscape variables, it is a partial sum process of dependent summands.

In \cite{BC4} the authors solve this problem in the setting of symmetric ($a=0$) trap models.
They give a set of abstract conditions that ensure that the clock process
converges to a stable subordinator. Technically, these  conditions bear, mainly, on the potential
theory of the chain $J_n$ and on the distribution of the random landscape.
By way of illustration, these results are then applied to the RHT dynamics of the REM
for which aging in proved on shorter time scales and higher temperatures than those
considered in \cite{BBG1,BBG2}.

At this point all models for which an arcsine aging regime had been proved shared two main non physical
features:

\item{(1)} the landscape in made of independent and identically distributed traps, and
\item{(2)}the dynamics is symmetric ($a=0$),
implying that the chain $J_n$ is a homogeneous random walk, independent of the trapping landscape;
this is to be contrasted with  the asymmetric case ($a>0$),
where $J_n$ favors jumps to the neighboring traps of largest depths,
as would be the case in a classical Glauber dynamics.

\noindent Moreover all know results were obtained either almost surely or in probability with respect
to law of the random landscape
(it is important to keep in mind that almost sure results do not always hold and that
in probability results sometimes are, as in \cite{BBG2}, the strongest statement possible):
let us momentarily stretch the terminology and call such results {\it quenched}.

In \cite{BBC} a model with correlation was for the first time considered, namely the
$p$-spin SK spin glass model, evolving, as in the REM, under the RHT dynamics.
If the abstract results of \cite{BC4} do in principle allow to treat situations with correlations,
too little is known about the random landscape of the $p$-spin SK model to actually carry them through.
To circumvent this difficulty the authors propose to take the ``view point of the particle''
and, rather than looking for  quenched results, average over the landscape variables
while conditioning on the trajectories of the chain $J_n$.
Then, adapting the arcsine aging scheme to this framework, they prove that,
for appropriate choices of time scales and parameters,
aging is again the same as in Bouchaud's symmetric trap model on the complete graph.

In the present paper we adopt yet another approach, completely different from those of
\cite{BC4} and \cite{BBC}, which will allow us to both
implement the arcsine aging scheme in the general setting that we called earlier
``Glauber dynamics of discrete disordered systems'', and obtain quenched results.
This approach is based on a powerful and illuminating method developed
by S\. Resnick and R\. Durrett \cite{DuRe}
to prove functional limit theorems for dependent variables.
By extending the framework of \cite{DuRe} to our random setting, and specializing it
to processes of the form \eqv(1.clock), we give simple sufficient conditions for $S_n$
to converge to a subordinator. An important aspect of the method is that it yields
an explicit expression of the limiting subordinator in terms
of the two-dimensional Poisson point process
that describes its jumps sizes and jumps times.
This result is the content of Theorem \thv(1.3.theo1) and the core of the paper.

A description of the organization of the paper is now in order.
As we have already announced it is made of two distinct parts: an abstract one
(that consists of Sections 1 and 2) and an applied one (formed of Sections 3 to 8).
In the rest of this first section we introduce our general setting (Subsection 1.1),
the necessary notions and definitions about aging (Subsection 1.2), and state
our main results on convergence of the clock process and its associated time-time correlation function
(Subsection 1.3): we will distinguish the {\it pure} process, whose initial increment is zero
(see Theorem \thv(1.3.theo1)), from  the {\it full} or {\it delayed} process, whose initial increment
depends on the initial distribution  (Theorem \thv(1.3.theo2)).
Section 2 contains the proofs of these results. It also contains the statement and proofs
of their counterparts
for the asymmetric trap model on the complete graph,
a model for which the notion of
convergence to renewal processes, and not only to subordinators, is relevant
(see Theorem  \thv(2.4.theo1) and Theorem  \thv(2.4.theo2)).
Section 3 begins the investigation of Bouchaud's asymmetric trap model on the complete
graph proper. It contains a separate introduction, and the statement of the results. Their
proofs occupy the rest of the paper, up to the appendix.


\bigskip
\line{\bf 1.1. The setting.\hfill}

Let $G_n(\VV_n,\EE_n)$, $n\in\N$, be a sequence of  connected graphs
with set of vertices $\VV_n$ and set of (non oriented) edges $\EE_n$.
A {\it random landscape} on $\VV_n$ (or {\it random environment}) is
a family $(\t_n(x), x\in\VV_n)$ on $\VV_n$ non-negative random variables.
As we shall want to take $n\uparrow\infty$ limits we assume that the sequence of these families
can be defined on a common probability space $(\O^{\t}, \FF^{\t}, \P)$.
Note that we do not assume a priori that the $\t_n(x)$'s are i.i.d\..
Using the random landscape a positive random measure $\t_n$
is defined on $\VV_n$ by,
$$
\t_n=\sum_{x\in\VV_n}\t_n(x)\d_x\,,
\Eq(1.1.1)
$$
where $\d_x$ is the point mass at $x$.
We call $\t_n$ the non-normalized Gibbs measure and, whenever $\t_n$ has finite total mass,
define the Gibbs measure through
$$
\GG_n=\sum_{x\in\VV_n}\frac{\t_n(x)}{\sum_{x\in\VV_n}\t_n(x)}\d_x\,.
\Eq(1.1.2)
$$

On $\VV_n$ we consider a continuous time Markov chain $(X_n(t), t\geq 0)$ that moves along the edges
of $G_n$ and is reversible w.r.t\. the measure $\t_n$. We may describe this
chain using its infinitesimal generator matrix, $\L_n=(\l_n(x,y))_{x,y\in\VV_n}$,
by requiring that all transition rates off $\EE_n$ are zero, that is $\l_n(y,x)=0$ for all $(x,y)\notin\EE_n\,,x\neq y$,
whereas on $\EE_n$, they satisfy the detailed balance condition
$$
\t_n(x)\l_n(x,y)=\t_n(y)\l_n(y,x)\,,\quad\forall\, (x,y)\in\EE_n\,,x\neq y\,.
\Eq(1.1.3)
$$
There are clearly many ways to choose such $\l_n$'s. For $\L_n$ to be an infinitesimal generator matrix
they must obey the constraint
$$
\l_n(x):=\sum_{y\in\VV_n}\l_n(x,y)<\infty\,,\quad\forall x\in\VV_n\,,
\Eq(1.1.4)
$$
and the diagonal elements $\l_n(x,x)$ must be set to $-\l_n(x)$.
(Here we assumed
that the graph $G_n$ contains no loops, i.e\.
$(x,x)\notin\EE_n$. If this is not the case
one must first suppress them using the appropriate time change.)
Finally, we make the extra assumption that $\l_n(x)\neq 0$ for all $x\in\VV_n$.

An alternative to the above construction is to describe the chain $X_n$ in terms of the joint distribution of its
jump chain and holding times. The jump chain  of $X_n$ is a discrete time Markov chain
$(J_n(k)\,,k\in\N)$  with transition
probability matrix $\G_n=(p_n(x,y))_{x,y\in\VV_n}$,
$$
p_n(x,y)=
\cases
\lambda_n(x,y)/\lambda_n(x)
&\hbox{if $(x,y)\in\EE_n\,,x\neq y$} , \,\,\,\cr
0 , &\hbox{otherwise.}\,\,\,\cr
\endcases
%
\Eq(1.1.5)
$$
It describes the sequence of states visited by $X_n$, the length of each visit at a given site, say $x$,
being exponential with parameter $\l_n(x)$.
To make this precise let the {\it clock process} of $X_n$ be defined through
$$
\wt S_n(k)=\sum_{i=0}^{k}\l_n^{-1}(J_n(i))e_{n,i}\,,\quad k\in \N\,,
\Eq(1.1.6)
$$
where $(e_{n,i}\,,n\in\N, i\in\N)$ is a family of independent mean one exponential random variables, independent of $J_n$.
Then, if
$X_n$ has initial distribution $\mu_n$, $J_n$ has initial distribution $\mu_n$ and
$$
X_n(t)= J_n(i) \text{if} \wt S_n(i)\leq t<\wt S_n(i+1) \text{for some} i\,.
\Eq(1.1.7)
$$
Given an initial distribution $\mu_n$ we write $\PP_{\mu_n}$ for the law of $X_n$
and $P_{\mu_n}$ for the law of $J_n$.
In view of taking $n\uparrow\infty$ limits we assume that
the sequences of chains  $X_n$, resp\. $J_n$, can be constructed on a common
probability space $(\O^{X}, \FF^{X}, \PP)$, resp\. $(\O^{J}, \FF^{J}, P)$.
We refer to \cite{FeGa} for an explicit construction.
Expectation with respect to $\P$, $P$, and $\PP$ will be denoted respectively by $\E$, $E$, and $\EE$.

\vfill\eject
\line{\bf 1.2. Aging.\hfill}

To study aging one needs to choose three ingredients:

(1) An initial distribution, which we denote by $\mu_n$.\hfill\break
\null\hskip.45truecm (2) A time scale, $c_n$, on which to observe $X_n$;
$c_n$ can either be a constant (in which case we may take $c_n=1$) or a positive
increasing sequence satisfying $c_n\uparrow\infty$ as $n\uparrow\infty$.\hfill\break
\null\hskip.45truecm
(3) A time-time correlation function, that is, a function $\CC_{n}(t,s)$ that gives
some interesting information on how much $X_n(c_n(t+s))$ depends on $X_n(c_nt)$ for $t, s\geq 0$.
A list of the functions commonly used in the literature can be found in \cite{BC4}.
 In Theorem \thv(1.3.theo1) below we will make the following choice:
$$
\CC_{n}(t,s)=\PP_{\mu_n}\left(\left\{c_n^{-1}\wt S_n(i)\,, i\in\N\right\}
\cap (t, t+s)=\emptyset\right)\,,\quad 0\leq t<t+s\,.
\Eq(1.1.8)
$$
Namely, this is the probability that the range of the re-scaled clock process $c_n^{-1}\wt S_n$
does not intersect the time interval $(t, t+s)$. In the arcsine aging scheme of \cite{BC4}, one aims
at controlling this probability asymptotically, in the limit of large $n$ and/or long times $t,s$,
using the Dynkin-Lamperti arcsine law for subordinators
(see Theorem \thv(A.2.theo2) of Appendix A.2). With this in mind we make the following definitions.

\definition{\TH(1.1.def1)} {We say that a time-time correlation function $\CC_{n}$ exhibits normal aging
on time scale $c_n$ if
one of the following three relations holds true:
$$
\lim_{t\rightarrow 0}\lim_{n\rightarrow\infty}
\CC_{n}(t,\rho t)=\CC_{\infty}(\rho)\,,
\Eq(1.1.14)
$$
$$
\lim_{n\rightarrow\infty}
\CC_{n}(t,\rho t)=\CC_{\infty}(\rho)\,,\,\,\,t>0\,\,\,\hbox{\rm arbitrary,}
\Eq(1.1.13)
$$
$$
\lim_{t\rightarrow\infty}\lim_{n\rightarrow\infty}
\CC_{n}(t,\rho t)=\CC_{\infty}(\rho)\,,
\Eq(1.1.12)
$$
for all $\rho\geq 0$, some non trivial limiting function
\note{
In all generality the r\.h\.s\. of \eqv(1.1.14)-\eqv(1.1.12) are not necessarily the same.
}
$\CC_{\infty}: [0,\infty)\mapsto [0,1]$, and for some convergence mode w.r.t\.
the probability law $\P$ of the random landscape.

We are now equipped to give a formal definition of what we  called earlier an arcsine aging regime.
Let $\asl_{\a}$ denote the distribution function of the generalized arcsine law with parameter $0<\a<1$,
$$
\asl_{\a}(u)=\frac{\sin \a\pi}{\pi}\int_0^u (1-x)^{-\a}x^{\a-1}dx\,,\quad 0\leq u\leq 1\,.
\Eq(1.1.16)
$$

\definition{\TH(1.1.def2)} {
We say that the process $X_n$ has an arcsine aging regime with parameter $\a$ whenever one can find
a time-time correlation function $\CC_{n}$ exhibiting normal aging with
$$
\CC_{\infty}(\rho)=\asl_{\a}(1/1+\rho)\,.
\Eq(1.1.15)
$$
}

\bigskip
\line{\bf 1.3. Convergence of the clock process to a subordinator.\hfill}

As we will see the first increment of the clock process plays a special role. For this reason we define
$$
\s_n=c_n^{-1}\wt S_n(0)\,,\quad
\overline S_n(k)=
\cases
c_n^{-1}\sum_{i=1}^{k}\lambda_n^{-1}(J_n(i))e_{n,i}
&\hbox{if $k\geq 1$} , \,\,\,\cr
0 , &\hbox{otherwise.}\,\,\,\cr
\endcases
\Eq(1.3.2'')
$$
Given a  positive (possibly constant) sequence $a_n$ we then set, for $t\geq 0$,
$$
S_n(t)=\overline S_n(\lfloor a_n t\rfloor)\,,
\Eq(1.3.2)
$$
and
$$
\wh S_n(t)=\s_n+S_n(t)\,.
\Eq(1.3.2')
$$
The re-scaled clock processes $S_n(t)$ and $\wh S_n(t)$ will be called, respectively, {\it pure} and {\it full} or {\it delayed}.
Note that $\bigl\{c_n^{-1}\wt S_n(i)\,, i\in\N\bigr\}=\bigl\{\wh S_n(u)\,,u>0\bigr\}$, that is, the  processes
$c_n^{-1}\wt S_n$ and $\wh S_n$ have identical range.
Also note that \eqv(1.1.8) may be rewritten as
$$
\CC_{n}(t,s)=\PP_{\mu_n}\left(\left\{\wh S_n(u)\,,u>0\right\}
\cap (t, t+s)=\emptyset\right)\,,\quad 0\leq t<t+s\,.
\Eq(1.3.3)
$$

We now state three conditions, (A1)-(A3),
that ensure that the pure process $S_n$ converges to a subordinator.
Because this process is a random variable on the probability space $(\O^{\t}, \FF^{\t}, \P)$
of the landscape (our random environment) we must first decide in which sense to seek convergence on
that space.
The relevant convergence modes (those which will be needed in practice)
are almost sure convergence and convergence in probability. This means
that one of the following statements should be in force:

\noindent{\it Almost sure convergence:}
There exists a subset $\wt\O^{\t}\subset\O^{\t}$ such that $\P(\wt\O^{\t})=1$
and such that, for all $\o\in \wt\O^{\t}$, for all large enough $n$, (A1)-(A3) are verified.

\noindent{\it Convergence in probability:} There exists a sequence $\wt\O^{\t}_n\subset\O^{\t}$ such that
$\lim_{n\rightarrow\infty}\P(\wt\O^{\t}_n)=1$
and such that, for all large enough $n$, (A1)-(A3) are verified  for all $\o\in \wt\O^{\t}_n$.


\noindent We now state our three conditions for fixed $\o$ and make this explicit by adding  the superscript $\o$
to landscape dependent quantities.
Since these conditions depend on
the choice of the initial distribution $\mu_n$, and of the sequences $a_n$ and $c_n$,
their formulation must thus be preceded by the statement:
``Given a sequence  of initial distributions $\mu_n$, there exist positive sequences $a_n$ and $c_n$ such that the following holds.''


\noindent{\bf Condition (A1).}
There exists a $\s$-finite measure $\nu$ on $(0,\infty)$ satisfying
$\int_{(0,\infty)}(1\wedge u)\nu(du)<\infty$ such that, for all $t>0$ and all $u>0$,\note{
The set $\wt\O^{\t}$ (respectively the sequence of sets $\wt\O^{\t}_n$)
for which convergence w\.r\.t\. the environment holds almost surely
(respectively in probability) is (are) the same for all $t>0$ and $u>0$.
}
$$
P^{\o}\left(
\left|
\sum_{j=1}^{\lfloor a_n t\rfloor}
\sum_{x\in\VV_n}p_n^{\o}(J_n^{\o}(j-1),x)
e^{-uc_n\l^{\o}_n(x)}
-t\nu(u,\infty)
\right|
<\e
\right)=1-o(1)\,,\quad\forall\e>0\,.
\Eq(1.A1)
$$

\noindent{\bf Condition (A2).}  For all $u>0$ and all $t>0$,
$$
P^{\o}\left(
\sum_{j=1}^{\lfloor a_n t\rfloor}\left[
\sum_{x\in\VV_n}p_n^{\o}(J_n^{\o}(j-1),x)
e^{-uc_n\l^{\o}_n(x)}
\right]^2
<\e
\right)=1-o(1)\,,\quad\forall\e>0\,.
\Eq(1.A2)
$$

\noindent{\bf Condition (A3).}  There exists a sequence of functions $\varepsilon_n\geq 0$  satisfying
$
\displaystyle\lim_{\d\rightarrow 0}\limsup_{n\rightarrow \infty}\varepsilon_n(\d)=0
$
such that 
for some $0<\d_0\leq 1$, for all $0<\d\leq\d_0$ and all $t>0$,
$$
E^{\o}\left(\int_{0}^{\d}du
\sum_{j=1}^{\lfloor a_n t\rfloor}\sum_{x\in\VV_n}p_n^{\o}(J_n^{\o}(j-1),x)e^{-uc_n\l^{\o}_n(x)}
\right)
\leq t\varepsilon_n(\d)\,.
\Eq(1.A3)
$$

\theo{\TH(1.3.theo1)} {\it For all sequences of initial distributions  $\mu_n$ and all sequences  $a_n$ and $c_n$
for which Conditions (A1), (A2), and (A3) are verified,
either $\P$-almost surely or in $\P$-probability,
the following holds w.r.t\. the same convergence mode:
%
Let $\{(t_k,\xi_k)\}$ be the points of a Poisson random measure of intensity measure
$dt\times d\nu$.  We have,
$$
S_n(\cdot)\Rightarrow S(\cdot)=\sum_{t_k\leq \cdot}\xi_k\,,
\Eq(1.3.theo1.1)
$$
in the sense of weak convergence
in the space $D([0,\infty))$ of c\`adl\`ag functions on
$[0,\infty)$ equipped with the Skorohod $J_1$-topology\note{
see e\.g\. \cite{W} p\. 83 for the definition of convergence in $D([0,\infty))$.
}.
}

\remark Although we do not make this explicit in the notation,
note that the limiting subordinator $S(\cdot)$ may remain a random variable on the
probability space $(\O^{\t}, \FF^{\t}, \P)$ of the random landscape (or some representation
of this space). We will see an example of this in the context of the asymmetric trap model
on the complete graph (see Proposition \thv(4.prop3) of Section 4.3).

\remark A sufficient condition for (A3) is given in Lemma \thv(2.2.A3').

To obtain convergence of the full re-scaled clock process $\wh S_n(\cdot)$ of \eqv(1.3.2'),
we still need to control the initial increment $\s_n$.
For this we introduce a separate condition. With the same notations and conventions as before:

\noindent{\bf Condition (A0).}
There exists a continuous distribution function $F^{\o}$ on $[0,\infty)$ such that,
for all $v\geq 0$,
$$
\left|
\sum_{x\in\VV_n}\mu_n^{\o}(x)e^{-vc_n\l^{\o}_n(x)}
-(1-F^{\o}(v))\right|=o(1)\,.
\Eq(1.A0)
$$

\theo{\TH(1.3.theo2)} {\it For all sequences of initial distributions  $\mu_n$ and all sequences  $a_n$ and $c_n$
for which Conditions (A0), (A1), (A2), and (A3) are verified,
either $\P$-almost surely or in $\P$-probability,
the following holds w.r.t\. the same convergence mode.
For $S(\cdot)$ defined in \eqv(1.3.theo1.1):

\item{(i)} Let $\s$ denote the random variable of (possibly random) distribution function $F$. Then,
$$
\wh S_n(\cdot)\Rightarrow \wh S(\cdot)=\s+S(\cdot)\,,
\Eq(1.3.theo2.1)
$$
(where $\Rightarrow$  has the same meaning as in \eqv(1.3.theo1.1)).

\item{(ii)}
Set
$$
\CC_{\infty}(t,s)=\PP\left(\left\{S(u)\,,u>0\right\}
\cap (t, t+s)=\emptyset\right)\,,\quad 0\leq t<t+s\,.
\Eq(1.3.theo1.3)
$$
If, for each
$\o\in\O^{\t}$, $\s$ and $S(\cdot)$ in \eqv(1.3.theo2.1) are independent
r\.v\.'s
on $(\O^{X}, \FF^{X}, \PP)$, then, for all $0\leq t<t+s$,
$$
\lim_{n\rightarrow\infty}\CC_{n}(t,s)=1-F(t+s)+\int_{0}^t\CC_{\infty}(t-v,s)dF(v)\,.
\Eq(1.3.theo1.4')
$$
In particular, if $\s=0$,
$$
\lim_{n\rightarrow\infty}\CC_{n}(t,s)=\CC_{\infty}(t,s)\,.
\Eq(1.3.theo1.4)
$$
}

In words, Theorem \thv(1.3.theo1) states that the process $S(\cdot)$
is a subordinator of L\'evy measure $\nu$.
Thus, by assertion (i) of Theorem \thv(1.3.theo2),
if $\s$ and $S(\cdot)$ in \eqv(1.3.theo2.1) are independent,
the  process $\wh S(\cdot)$ is a delayed subordinator.
Since the arcsine law for subordinators
(stated in Theorem \thv(A.2.theo2) of the Appendix)
provides us with necessary and sufficient conditions for $\CC_{\infty}(t,s)$ to be, or converge
to, the distribution function of the generalized arcsine law,
assertion (ii) of Theorem \thv(1.3.theo2) yields criteria for the process
$X_n$ to have an arcsine aging regime.

In trying to verify the two conditions (A1) and (A2) one should be guided by the fact that
they are kinds of ergodic theorems in a random environment.
The asymmetric trap model on the complete graph is not a good working ground to illustrate
this idea as ergocitity
is obtained trivially
(see the remark below \eqv(1.AA2) in the proof of Theorem \thv(2.4.theo1)).
A more involved model where this idea is clearly illustrated is
the random hopping time dynamics of the REM studied in \cite{G1} (see Section 1.4).

Let us finally note that the form of the relation \eqv(1.3.theo1.4'), where the role of the initial
distribution $\mu_n$ is made explicit, is new.
For all models where the existence of an arcsine aging regime has been
proved so far, the initial distribution was chosen in such a way that $\s=0$.
In Section 3.4 we will give examples of initial distributions such that $\s\neq 0$,
and for which the arcsine aging regime still prevails.
A full investigation of the impact of the initial distribution
on the aging phenomenon will be carried out in \cite{G3}.

\remark In line with the remark following Theorem \thv(1.3.theo1), let us recall that
$\s$ and/or $S(\cdot)$ may be random variables on $(\O^{\t}, \FF^{\t}, \P)$. Thus
both the limiting functions in \eqv(1.3.theo1.4') and \eqv(1.3.theo1.4)
may be random variables on that space. (We will see an instance of this in
Theorem \thv(4.theo4) of Section 3.2.).
This is why we assume in Theorem \thv(1.3.theo1), (ii), that
$\s$ and $S(\cdot)$ are independent for each $\o\in\O^{\t}$.
We could make weaker assumptions: this one is designed to cover the needs of Section 3.

\remark Clearly, Conditions (A1)-(A3) can be verified only if $a_n$ is an increasing and diverging sequence. In the case of
constant sequence, say $a_n=1$, time remains discrete in the limit $n\uparrow\infty$:
if the clock process converges to a limiting object, the latter has to be a process of partial sums.
We will give an example of this in Theorem \thv(2.4.theo1)
where we will see that the limiting partial-sum process is a renewal process.

\vfill\eject

\bigskip


\chap{2. Convergence of the clock process and related results}2

This section is divided in four parts.
In Subsection 2.1 we state a result by Durrett and Resnick \cite {DuRe}
that is central to the proof of Theorem \thv(1.3.theo1) and Theorem \thv(1.3.theo2).
The proofs of the latter ones are given in Subsection 2.2
(which focuses on convergence of the processes $S_n$ and
$\wh S_n$ to subordinators) and Subsection 2.3
(where convergence of the time-time correlation function $\CC_{n}(t,s)$
is established).
In Subsection 2.4 we specialize Theorem \thv(1.3.theo1) and Theorem \thv(1.3.theo2)
to the asymmetric trap model on the complete graph, and complement these results
with sufficient conditions for convergence of the re-scaled clock process
to a partial-sum process in the case, not covered by Theorem \thv(1.3.theo1) and Theorem \thv(1.3.theo2),
where the auxiliary time scale $a_n$ is a constant
(see Theorem \thv(2.4.theo1) and Theorem \thv(2.4.theo2)).
What we need to know about subordinators and renewal processes is summarized in Section A.2 of the Appendix.


\bigskip
\line{\bf 2.1. A result by Durrett and Resnick \hfill}



In \cite{DuRe} a method is developed for proving convergence of partial sums processes
with dependent increments to L\'evy processes.
This method consists of two steps. In the first step, one shows that a sequence of point processes associated
with the increments converges weakly to a two dimensional Poisson process. Then,
applying appropriate functionals (to `sum up the points') and continuity arguments,
one obtains weak convergence of the sum to a limiting L\'evy process.

In this section we specialize this result, namely Theorem 4.1 of \cite {DuRe},
to the case of processes with non-negative increments.
Our framework is the following.
Let $\{Z_{n,i}, n\geq 1, i\geq 1\}$, $Z_{n,i}\geq 0$, be an array of random variables defined on a probability space
$(\O,\FF,\PP)$ and let $\{\FF_{n,i}, n\geq 1, i\geq 0\}$ be an array of sub-sigma fields of $\FF$ such that
for each $n$ and $i\geq 1$, $Z_{n,i}$ is $\FF_{n,i}$ measurable and $\FF_{n,i-1}\subset\FF_{n,i}$.
Let $k_n(t)$ be a nondecreasing right continuous function with range $\{0,1,2,\dots\}$ and assume that for
each $t>0$ $k_n(t)$ is a stopping time.
Set
$$
\wt S_{n,k}=\sum_{i=1}^k Z_{n,i}\,,
\Eq(2.1.2)
$$
for $k\geq 1$, $\wt S_{n,0}=0$, and  define
$$
S_{n}(t)=\wt S_{n,k_n(t)}\,.
\Eq(2.1.3)
$$
The next theorem gives conditions for $S_n$ to converge to a subordinator.
%
%
To state it we will need the following extra notation:
for $\d\geq 0$ set $Z^{\d}_{n,i}=Z_{n,i}\1_{\{Z_{n,i}\leq\d\}}$; further set
$$
\eqalign{
&\wt S^{\d}_{n,k}=\sum_{i=1}^k Z^{\d}_{n,i}\,,\cr
}
\Eq(2.1.4)
$$
for $k\geq 1$, $\wt S^{\d}_{n,0}=0$, and define
$$
S^{\d}_{n}(t)=\wt S^{\d}_{n,k_n(t)}\,.
\Eq(2.1.5)
$$

\theo{\TH(2.1.theo2)} (Durrett and Resnick). {\it
Let $\nu$ be a $\s$-finite measure on $(0,\infty)$ satisfying
$\int_{(0,\infty)}(1\wedge x)\nu(dx)<\infty$, and let $\{S(t), t\geq 0\}$ be the subordinator
of Laplace exponent
$
\Phi(\theta)=\int_{(0,\infty)}\left(1-e^{-\theta x}\right)\nu(dx)
$,
$\theta\geq 0$.
If, as $n\rightarrow\infty$,
\item{(D1)} For all $t>0$ and for $x>0$ such that $\nu(\{x\})=0$,
$$
\sum_{i=1}^{k_n(t)}
\PP\left(Z_{n,i}>x \mid \FF_{n,i-1}\right)
@>proba>>
t\nu(x,\infty)\,,
\Eq(2.1.theo2.1)
$$
\item{(D2)} For all $t>0$ and and all $\e>0$,
$$
\sum_{i=1}^{k_n(t)}
\left[\PP\left(Z_{n,i}>\e \mid \FF_{n,i-1}\right)\right]^2
@>proba>>0\,,
\Eq(2.1.theo2.2)
$$

\noindent and
\item{(D3)} For all $t>0$ and all $\e>0$,
$$
\lim_{\d\rightarrow 0}\limsup_{n\rightarrow\infty}
\PP\left(S^{\d}_{n}(t)>\e\right)=0\,,
\Eq(2.1.theo2.3)
$$

\noindent then
$
S_n\Rightarrow S
$
in the space $D([0,\infty))$ of c\`adl\`ag functions on $[0,\infty)$ equipped with the Skorohod topology.

}

\remark In \cite{DuRe}, Conditions (D2) and (D3) are stated for $t$ fixed and equal to one.
This does not seem to be correct.

\bigskip
\line{\bf 2.2. Convergence to subordinators.
\hfill}

In this subsection we prove Theorem \thv(1.3.theo1) and the first assertion of Theorem \thv(1.3.theo2),
and give an alternative to Condition (A3).

\proofof{Theorem \thv(1.3.theo1)}
Our aim is to apply Theorem \thv(2.1.theo2) to the sum
$$
S_n(t)=c_n^{-1}\sum_{i=1}^{\lfloor a_n t\rfloor}\l_n^{-1}(J_n(i))e_{n,i}\,.
$$
Let us first do this for a fixed realization $\o\in\O^{\tau}$ of the
environment.
Set $k_n(t)={\lfloor a_n t\rfloor}$,
$
Z_{n,i}=(c_n\l_n(J_n(i)))^{-1}e_{n,i}\,,
$
and (with obvious notation) define
$
\FF_{n,i}
=\BB(J_n(0),\dots,J_n(i), e_{n,0},\dots,e_{n,i})\,.
$
Clearly, for each $n$ and $i\geq 1$, $Z_{n,i}$ is $\FF_{n,i}$ measurable and $\FF_{n,i-1}\subset\FF_{n,i}$.
Next observe that
$$
\eqalign{
\PP_{\mu_n}\left(J_n(i)=x, Z_{n,i}>z\mid\FF_{n,i-1}\right)
&=\PP_{\mu_n}\left(J_n(i)=x, Z_{n,i}>z\mid J_{n}(i-1)\right)
\cr
&=p_n(J_n(i-1),x)\PP_{\mu_n}\left((\l_n(x))^{-1}e_{n,i}>z\right)
\cr
&=p_n(J_n(i-1),x) \exp\{-zc_n\l_n(x)\}\,.
}
\Eq(2.2.1)
$$
From this it follows that
$$
\sum_{i=1}^{k_n(t)}
\PP_{\mu_n}\left(Z_{n,i}>z \mid \FF_{n,i-1}\right)
=
\sum_{i=1}^{\lfloor a_n t\rfloor}
\sum_{x\in\VV_n}p_n(J_n(i-1),x)
\exp\{-zc_n\l_n(x)\}\,,
\Eq(2.2.2)
$$
and
$$
\sum_{i=1}^{k_n(t)}
\left[\PP_{\mu_n}\left(Z_{n,i}>\e \mid \FF_{n,i-1}\right)\right]^2
=
\sum_{i=1}^{\lfloor a_nt\rfloor}\left[
\sum_{x\in\VV_n}p_n(J_n(i-1),x)
\exp\{-uc_n\l_n(x)\}\right]^2
\,,
\Eq(2.2.3)
$$
so Condition (A2) and (A1) of Theorem \thv(1.3.theo1) are, respectively,
Conditions (D2) and condition (D1) of Theorem \thv(2.1.theo2).

We will now show that Condition (A1) together with Condition (A3) imply Condition (D3).
To simplify the notation in Conditions (A1)-(A3) we write $\bar\nu(u)\equiv\nu(u,\infty)$, and set
$$
\bar\nu_n^{J,t}(u)=
\sum_{j=1}^{\lfloor a_n t\rfloor}
\sum_{x\in\VV_n}p_n(J_n(j-1),x)
\exp\{-uc_n\l_n(x)\}\,.
\Eq(2.2.5)
$$
Consider now Condition (D3).
By Tchebychev inequality
$
\PP_{\mu_n}\left(S^{\d}_{n}(t)>\e\right)
\leq
\e^{-1}\EE_{\mu_n} S^{\d}_{n}(t)
$.
Expressed in terms of the truncated variables
$
Z^{\d}_{n,i}=Z_{n,i}\1_{\{Z_{n,i}<\d\}}
$,
$\d\geq 0$, the latter expectation becomes,
$$
\EE_{\mu_n} S^{\d}_{n}(t)
=\EE_{\mu_n}{{\sum_{j=1}^{\lfloor a_n t \rfloor}}}Z^{\d}_{n,i}
=E_{\mu_n}{{\sum_{j=1}^{\lfloor a_n t \rfloor}}}
\EE_{\mu_n}\left(Z^{\d}_{n,j}\,\Big|\, J_{n}(j-1)\right)\,.
\Eq(2.1.theo2.A3.3)
$$
Integrating by parts,
$$
\eqalign{
\EE_{\mu_n}\left(Z^{\d}_{n,i} \,\Big|\, J_{n}(i-1)\right)
=&\int_0^{\infty}\PP_{\mu_n}\left(Z^{\d}_{n,i}(J_n(i))>y \mid J_{n}(i-1)\right)dy
\cr
=&\int_0^{\d}\PP_{\mu_n}\left(Z_{n,i}\geq z \mid J_{n}(i-1)\right)dz
-\d\PP_{\mu_n}\left(Z_{n,i}>\d \mid J_{n}(i-1)\right)\,,
}
\Eq(2.1.theo2.A3.4)
$$
and since
$
\sum_{i=1}^{\lfloor a_n t\rfloor }
\PP\left(Z_{n,i}>z \mid J_{n}(i-1)\right)
=\bar\nu_n^{J,t}(u)
$,
as follows from \eqv(2.2.2) and \eqv(2.2.5), we arrive at
$$
\EE_{\mu_n} S^{\d}_{n}(t)=E_{\mu_n}\left(\int_{0}^{\d}du\bar\nu_n^{J,t}(u)-\d\bar\nu_n^{J,t}(\d)\right)\,.
\Eq(2.1.theo2.A3.5)
$$
Now by Condition (A1),
$
E_{\mu_n}\d\bar\nu_n^{J,t}(\d)\leq t\d\bar\nu(\d)+o(1)
$
and $\lim_{\d\rightarrow 0}\d\bar\nu(\d)=0$,
whereas Condition (A3) states that
$
E_{\mu_n}\left(\int_{0}^{\d}du\bar\nu_n^{J,t}(u)\right)\leq t\varepsilon_n(\d)
$,
where $\lim_{\d\rightarrow 0}\limsup_{n\rightarrow \infty}\varepsilon_n(\d)=0$.
Hence, if both these conditions are satisfied,
$
\lim_{\d\rightarrow 0}\limsup_{n\rightarrow \infty}\EE_{\mu_n} S^{\d}_{n}(t)=0
$,
so that Condition (D3) also is satisfied.

We may now conclude the proof of Theorem \thv(1.3.theo1).
We proved that  (A1)$\Rightarrow$(D1),
(D2)$\Leftrightarrow$(A2),
and (A1)$\&$(A3) $\Rightarrow$(D3).
Therefore, by Theorem \thv(2.1.theo2),
$
S_n\Rightarrow S
$
in $D([0,\infty))$
where $S$ is the subordinator \eqv(1.3.theo1.1).

So far we kept $\o\in\O^{\tau}$ fixed, i.e\.
we worked with a fixed realization of the
environment.
Let us now introduce the subsets $\O^{\tau}_{n,1},\O^{\tau}_{n,2}\subset\O^{\tau}$,
$$
\eqalign{
\O^{\tau}_{n,1}=&
\left\{
\forall t>0,\forall u>0,\forall\e>0,
P\left(
\left|
\sum_{j=1}^{\lfloor a_n t\rfloor}
\sum_{x\in\VV_n}p_n^{\o}(J_n^{\o}(j-1),x)
e^{-uc_n\l^{\o}_n(x)}
-t\nu(u,\infty)
\right|
<\e
\right)=1-o(1)
\right\}
\cr
\O^{\tau}_{n,2}=&
\left\{
\forall t>0,\forall\e>0,
P\left(
\sum_{j=1}^{\lfloor a_n t\rfloor}\left[
\sum_{x\in\VV_n}p_n^{\o}(J_n^{\o}(j-1),x)
e^{-uc_n\l^{\o}_n(x)}
\right]^2
<\e
\right)=1-o(1)
\right\}
\cr
}\,,
\Eq(2.2.14)
$$
and set $\wt\O^{\tau}_{n}=\O^{\tau}_{n,1}\cap\O^{\tau}_{n,2}$.
By definition of weak convergence what we have  just established is
that for each $\o\in\wt\O^{\tau}_{n}$, and large enough $n$,
$$
\left|\EE\left(f(S_n)\right)-\EE\left(f(S)\right)\right|=o(1)\,,
\Eq(2.2.15)
$$
for each continuous  bounded function $f$ on the space $D([0,\infty))$
equipped with Skorohod metric $\rho_\infty$.
If it holds true that
$
\P\left(\bigcup_m\bigcap_{n>m}\wt\O^{\tau}_{n}\right)=1
$,
then
$
S_n\Rightarrow S
$
$\P$-almost surely. If instead we have
$\lim_{n\rightarrow\infty}\P(\wt\O^{\tau}_{n})=1$,
then
$
S_n\Rightarrow S
$
in $\P$-probability. Theorem \thv(1.3.theo1) is thus proved. \endproof

\proofof{assertion (i) of Theorem \thv(1.3.theo2)}
 As in the proof of Theorem \thv(1.3.theo1) we first establish
\eqv(1.3.theo2.1) for a fixed realization $\o\in\O^{\tau}$ of the environment.
Note that the additional Condition (A0)
is designed to guarantee that $\s_n$ converges in distribution to $\s$. Indeed, since
$
\s_n=c_n^{-1}\wt S_n(0)=c_n^{-1}\lambda_n^{-1}(J_n(0))e_{n,0}
$,
we have
$
1-\PP_{\mu_n}(\s_n<v)=\sum_{x\in\VV_n}\mu_n(x)e^{-vc_n\l_n(x)}
$,
so that \eqv(1.A0)
becomes
$
\left|\PP_{\mu_n}(\s_n<v)-F(v)\right|=o(1)
$.
Thus, supplementing Conditions (A1) and (A2) with Condition (A0), it follows from
Theorem \thv(1.3.theo1)
that,
viewing $\s_n$ as a constant function in $D([0,\infty))$,
the pairs $(\s_n, S_n(\cdot))$ jointly converge, weakly, to the pair $(\s, S(\cdot))$,
in $D^2([0,\infty))$.
It next follows from the continuous mapping theorem, upon adding $\s_n$ and $ S_n(\cdot)$, that
$
\s_n+ S_n(\cdot)\Rightarrow \wh S(\cdot)=\s+ S(\cdot)
$
in $D([0,\infty))$
(see \cite{W}, p\. 84, last paragraph of Section 3.3, for the continuity of the
addition of an arbitrary element of $D([0,\infty))$ and the constant function).
Eq\. \eqv(1.3.theo2.1) being established for a fixed realization $\o\in\O^{\tau}$,
we conclude the proof proceeding exactly as in the proof of Theorem \thv(1.3.theo1)
\note{
see the paragraph beginning above \eqv(2.2.14).},
introducing the extra subsets
$
\O^{\tau}_{n,3}=
\left\{
\left|\sum_{x\in\VV_n}\mu_n(x)e^{-vc_n\l_n(x)}-(1-F(v))\right|=o(1)
\right\}
$
in \eqv(2.2.14), and setting
$
\wt\O^{\tau}_{n}=\O^{\tau}_{n,1}\cap\O^{\tau}_{n,2}\cap\O^{\tau}_{n,3}
$.
\endproof

Condition (A3) may not always be easy to handle. Here is an alternative:

\noindent{\bf Condition (A3').} There exists a sequence of functions $\varepsilon_n\geq 0$  satisfying
$
\displaystyle\lim_{\d\rightarrow 0}\limsup_{n\rightarrow \infty}\varepsilon_n(\d)=0
$
such that, for some $0<\d_0\leq 1$, for all $0<\d\leq\d_0$ and all $t>0$,
$$
E_{\mu_n}\left(
\sum_{j=1}^{\lfloor a_n t\rfloor}\sum_{x\in\VV_n}p_n^{\o}(J_n^{\o}(j-1),x)\frac{\1_{\{(c_n\l^{\o}_n(x))^{-1}\leq \d\}}}{c_n\l^{\o}_n(x)}
\right)
\leq t\varepsilon_n(\d)\,.
\Eq(1.A3')
$$

\lemma{\TH(2.2.A3')}{\it A sufficient condition for (A3) is (A3').}

\proof We will show that if Condition (A1) and Condition (A3') then so is Condition (A3).
As in the proof of Theorem \thv(1.3.theo1) we write $\bar\nu(u)\equiv\nu(u,\infty)$
and let $\bar\nu_n^{J,t}(u)$ be defined through \eqv(2.2.5). Then \eqv(1.A3) of Condition (A3) becomes
$
E_{\mu_n}\bigl(\int_{0}^{\d}du\bar\nu_n^{J,t}(u)\bigr)\leq t\varepsilon_n(\d)
$.
Clearly,
$$
\textstyle
\int_{0}^{\d}du\bar\nu_n^{J,t}(u)
=
\sum_{j=1}^{\lfloor a_n t\rfloor}\sum_{x\in\VV_n}p_n^{\o}(J_n^{\o}(j-1),x)
\frac{1-e^{-\d c_n\l^{\o}_n(x)}}{c_n\l^{\o}_n(x)}\,.
\Eq(2.2.A3'.1)
$$
Now on the one hand, since $\frac{1-e^{-y}}{y}\leq e^{\rho}e^{-y}$, $0\leq y\leq \rho$,
$$
\textstyle
\frac{1-e^{-\d c_n\l^{\o}_n(x)}}{c_n\l^{\o}_n(x)}\1_{\{c_n\l^{\o}_n(x)\leq \d\}}
\leq \d e^{\rho}e^{-\d c_n\l^{\o}_n(x)}\1_{\{\d c_n\l^{\o}_n(x)\leq \rho\}}
\leq \d e^{\rho}e^{-\d c_n\l^{\o}_n(x)}\,,
\Eq(2.2.A3'.2)
$$
for all $\rho>0$, while on the other hand
$
\frac{1-e^{-\d c_n\l^{\o}_n(x)}}{c_n\l^{\o}_n(x)}\1_{\{\d c_n\l^{\o}_n(x)\geq \rho\}}
\leq\frac{\1_{\{(c_n\l^{\o}_n(x))^{-1}\leq \d/\rho\}}}{c_n\l^{\o}_n(x)}
$.
Inserting these two bounds in \eqv(2.2.A3'.1) yields
$$
\textstyle
\int_{0}^{\d}du\bar\nu_n^{J,t}(u)
\leq
\d e^{\rho}\bar\nu_n^{J,t}(\d)
+
\sum_{j=1}^{\lfloor a_n t\rfloor}\sum_{x\in\VV_n}p_n^{\o}(J_n^{\o}(j-1),x)\frac{\1_{\{(c_n\l^{\o}_n(x))^{-1}\leq \d/\rho\}}}{c_n\l^{\o}_n(x)}\,.
\Eq(2.2.A3'.4)
$$
Recall that by Condition (A1),
$
E_{\mu_n}\d\bar\nu_n^{J,t}(\d)\leq t\d\bar\nu(\d)+o(1)
$
where $\lim_{\d\rightarrow 0}\d\bar\nu(\d)=0$.
Thus, averaging out \eqv(2.2.A3'.4) and using Condition (A1) together with \eqv(1.A3') of Condition (A3')
to bound the resulting right hand side, we get that, for all $\rho>0$,
$
E_{\mu_n}\bigl(\int_{0}^{\d}du\bar\nu_n^{J,t}(u)\bigr)\leq
t\varepsilon_n(\d/\rho)
+\d e^{\rho}(t\d\bar\nu(\d)+o(1))
$.
Finally, taking e.g\. $\rho=\sqrt\d$,
$
\lim_{\d\rightarrow 0}\limsup_{n\rightarrow\infty}\{t\varepsilon_n(\d/\rho)+\d e^{\rho}(t\d\bar\nu(\d)+o(1))\}=0
$.
Condition (A3) is therefore satisfied.\endproof.

\bigskip
\line{\bf 2.3. Convergence of the time-time correlation function.
\hfill}

We will now exploit the convergence of $\wh S_n(\cdot)$ established above
to prove convergence of the time-time correlation function,
using the continuous-mapping theorem.

\proofof{assertion (ii) of Theorem \thv(1.3.theo2)}
This  pattern of proof is classical
(see \cite{W} section 9.7.2) and relies on the continuity
property of a certain function of the inverse mapping on $D([0,\infty))$,
the so-called overshoot, which we now define.
Let $\eta\in D([0,\infty))$.
For $t>0$ let $\LL_t$ be the time of the first passage to a level beyond $t$;
{\it i\.e\.},
$$
\LL_t(\eta)\equiv\eta^{-1}(t)\equiv\{\inf u\geq 0 \mid \eta(u)>t\}
\Eq(2.2.16)
$$
(with $\LL_t(\eta)=\infty$ if $\eta(u)\leq t$ for all $u$).
Let $D_t(\eta)=\eta(\LL_t(\eta))$ be the first visit to the set
$\left\{\eta(u)\,,u>0\right\}$
after time $t$.  The associated overshoot is the function  $\theta_t(\eta)$ defined through
$$
\theta_t(\eta)=D_t(\eta)-t\,.
\Eq(2.2.17)
$$
With this definition the time-time correlation function \eqv(1.3.3) may be rewritten as
$$
\CC_{n}(t,s)
=
\PP_{\mu_n}\left(\left\{\wh S_n(u)\,,u>0\right\}\cap (t, t+s)=\emptyset\right)
=
\PP_{\mu_n}\left(\theta_t\bigl(\wh S_n\bigr)\geq s\right)\,.
\Eq(2.2.18)
$$
%
%
Similarly, \eqv(1.3.theo1.3) can be rewritten as
$$
\CC_{\infty}(t,s)
=\PP\left(\left\{S(u)\,,u>0\right\}\cap (t, t+s)=\emptyset\right)
=
\PP\left(\theta_t(S)\geq s\right)\,.
\Eq(2.2.19)
$$
As announced,
the motivation behind this approach
is that the overshoot function is an almost
surely continuous function on $ D([0,\infty))$
with respect to L\'evy motions having almost surely diverging paths
(see \cite{W}, Theorem 13.6.5 p\.447).
Hence, if  \eqv(1.3.theo2.1) holds true $\P$-almost surely,
the continuous mapping theorem (applied for each fixed $\o$ that belongs to the set
of full measure for which $S_n\Rightarrow S$ obtains) readily yields that
$\P$-almost surely, uniformly
%
in $0\leq t<t+s$,
$
\lim_{n\rightarrow\infty}\PP_{\mu_n}\bigl(\theta_t(\wh S_n)\geq s\bigr)
=\PP\bigl(\theta_t(\wh S)\geq s\bigr)
$.
Assume now that
$
\wh S_n\Rightarrow \wh S
$
in $\P$-probability. Note that for each continuous  bounded function $g$ on $[0,\infty)$
the function $g\circ \theta_t$ is a continuous  bounded function on $D([0,\infty))$.
Thus, by \eqv(2.2.15), for each $\o\in\wt\O^{\tau}_{n}$ and large enough $n$,
$$
\bigl|\EE\bigl(g\circ \theta_t(\wh S_n)\bigr)-\EE\bigl(g\circ \theta_t(\wh S)\bigr)\bigr|=o(1)\,.
\Eq(2.2.20)
$$
From this and the definition of weak convergence it follows that
$
\lim_{n\rightarrow\infty}\PP_{\mu_n}\bigl(\theta_t(\wh S_n)\geq s\bigr)
=\PP\bigl(\theta_t(\wh S)\geq s\bigr)
$
in $\P$-probability. Since the sequence of subsets $\wt\O^{\tau}_{n}$ does not depend on $t$ and $s$,
convergence holds uniformly in $0\leq t<t+s$,
in $\P$-probability.

It remains to express $\PP\bigl(\theta_t(\wh S)\geq s\bigr)$ in terms of $\CC_{\infty}(t,s)$ and $F$.
If $\s=0$ then $\wh S=S$, and by \eqv(2.2.19),
$
\PP\bigl(\theta_t(\wh S)\geq s\bigr)=\CC_{\infty}(t,s)
$,
which proves \eqv(1.3.theo1.4). Otherwise, from the assumption that
$\s$ and $S(\cdot)$ in \eqv(1.3.theo2.1)  are independent r\.v\.'s
on  $(\O^{X}, \FF^{X}, \PP)$ for each fixed $\o\in\O^{\t}$, we get, conditioning on $\s$, that
$$
\PP\bigl(\theta_t(\wh S)\geq s\bigr)
=1-F(t+s)+\int_{0}^t\PP\bigl(\theta_{t-v}(S)\geq s\bigr)dF(v)
=1-F(t+s)+\int_{0}^t\CC_{\infty}(t-v,s)dF(v)\,.
\Eq(2.2.22)
$$
Since \eqv(2.2.22) holds true for each $\o\in\O^{\t}$ uniformly in $0\leq t<t+s$, \eqv(1.3.theo1.4') obtains
uniformly in $0\leq t<t+s$, and inherits the convergence mode of
$\PP_{\mu_n}\bigl(\theta_t(\wh S_n)\geq s\bigr)$, that is to say, the convergence mode of $\wh S_n$.
The proof of assertion (ii) of Theorem \thv(1.3.theo2) is now complete.
\endproof

\bigskip
\line{\bf 2.4. The special case of the asymmetric trap model on the complete graph.\hfill}
\line{\bf\hskip .85truecm (Convergence to renewal processes.)
\hfill}

In this section we focus on the asymmetric trap model defined in \eqv(3.1.2)-\eqv(3.1.3)
when $G_n(\VV_n, \EE_n)$ is the complete graph on $n$ vertices, a loop being attached to each vertex,
and when the landscape if formed of arbitrarily distributed i.i.d\. positive traps.
As already observed in Section 1, Theorem \thv(1.3.theo1) and Theorem \thv(1.3.theo2)
only cover situations where the auxiliary time scale $a_n$
of the re-scaled clock process \eqv(1.3.2) diverges with $n$,
leaving out the case of constant $a_n$.
It is obvious that in that latter case the partial-sum structure
of the clock process must be preserved in the limit, whenever a limit exists.
In Theorem \thv(2.4.theo1) and Theorem \thv(2.4.theo2) below we specialize
the results of Theorem \thv(1.3.theo1) and Theorem \thv(1.3.theo2)
to the asymmetric trap model on the complete graph,
and complement them with sufficient conditions for convergence of the re-scaled clock processes
to a partial-sum process, more precisely, to a renewal process.

For constant $a_n$ the sample paths of $S_n$ are increasing
functions on $[0,\infty)$ that have discontinuities at all integer time points.
The natural topological space in which to interpret
weak convergence of $S_n$ is, here, the space $\R^{\infty}$
of infinite sequences equipped with the usual Euclidean topology (see {\it e\.g\.} \cite{Bi} section 3).
We will use the arrow $\Rrightarrow$ to denote weak convergence in that space.
As in Theorem \thv(1.3.theo1), weak convergence in Skorohod topology on $D([0,\infty))$ will be denoted by $\Rightarrow$.
Set $r_n=c_n^{1/(1-a)}$ and define
$$
\nu_n(u,\infty)=a_n
\frac{
\sum_{x\in\VV_n}\t^a(x)
\exp\{-u(r_n/\t(x))^{(1-a)}\}
}
{
\sum_{x\in\VV_n}\t^a(x)
}\,,\quad u\geq 0\,.
\Eq(2.4.theo1.1)
$$

\theo{\TH(2.4.theo1)} {\it
Consider the asymmetric trap model on the complete graph
on time scale $c_n$.
The following holds for any choice of the initial distribution $\mu_n$.

\item{(i)} If
there exists
a sequence $a_n$ satisfying $a_n\uparrow\infty$ as $n\uparrow\infty$,
a $\s$-finite measure $\nu$ on $(0,\infty)$ satisfying $\int_{(0,\infty)}(1\wedge u)\nu(du)<\infty$,
and a function $\varepsilon\geq 0$ satisfying
$
\displaystyle\lim_{\d\rightarrow 0}\varepsilon(\d)=0
$,
such that, either $\P$-almost surely or in $\P$-probability, for all $u>0$,
$$
\lim_{n\rightarrow \infty}\nu_n(u,\infty)=\nu(u,\infty)\,,
\Eq(2.4.theo1.2)
$$
and, for all $0<\d\leq\d_0$, for some $0<\d_0\leq 1$,
$$
\limsup_{n\rightarrow \infty}\int_0^{\d}\nu_n(u,\infty)du\leq \varepsilon(\d)\,,
\Eq(2.4.theo1.2bis)
$$
%
then,
w.r.t\. the same convergence mode,
$$
S_n(\cdot)\Rightarrow S(\cdot)=\sum_{t_k\leq \cdot}\xi_k\,,
\Eq(2.4.theo1.3)
$$
where $\{(t_k,\xi_k)\}$ are the marks of a Poisson process on $[0,\infty)\times (0,\infty)$
with mean measure $dt\times d\nu$.

\item{(ii)} If, taking $a_n=1$, there exists a probability distribution $\nu$ on $(0,\infty)$
such that, either $\P$-almost surely or in $\P$-probability,  \eqv(2.4.theo1.2) is verified for all $u\geq 0$,
 then, w.r.t\. the same convergence mode,
$$
S_n(\cdot)\Rrightarrow R(\cdot)=\sum_{k\leq \cdot}\xi_k\,,
\Eq(2.4.theo1.4)
$$
where $\{\xi_k, k\geq 1\}$ are independent r\.v\.'s with identical
distribution $\nu$.
}

In the sequel we will adopt the terminology used in \cite{Fe}
and call the sequence
$\{R(k)\,,k\in\N\}$
a {\it renewal process of inter-arrival distribution $\nu$}
(equivalently, of {\it inter-arrival times $\xi_k$}).
As in Theorem \thv(1.3.theo2) the extra Condition (A0)
on the convergence of the initial increment $\s_n$
enables us to deduce convergence of the full
clock process $\wh S_n(\cdot)$ from that of $S_n(\cdot)$.

\theo{\TH(2.4.theo2)} {\it \item{(i')}  If, in addition to the assumptions of
assertion (i) of Theorem \thv(2.4.theo1), Condition (A0) is satisfied
w.r.t\. the same convergence mode as in \eqv(2.4.theo1.2), then, in this convergence mode,
denoting by $\s$ the random variable of (possibly random) distribution function $F$,
the following holds:
For $S(\cdot)$ defined in \eqv(2.4.theo1.3),
$$
\wh S_n(\cdot)\Rightarrow \wh S(\cdot)=\s+S(\cdot)\,,
\Eq(2.4.theo2.1)
$$
where $\s$ and $S(\cdot)$ are independent.
Moreover for $\CC_{\infty}(t,s)$ defined in \eqv(1.3.theo1.3),
for all $0\leq t<t+s$,
$$
\lim_{n\rightarrow\infty}\CC_{n}(t,s)=1-F(t+s)+\int_{0}^t\CC_{\infty}(t-v,s)dF(v)\,.
\Eq(2.4.theo2.3)
$$
In particular, if $\s=0$,
$$
\lim_{n\rightarrow\infty}\CC_{n}(t,s)=\CC_{\infty}(t,s)\,.
\Eq(2.4.theo2.4)
$$

\item{(ii')} Substituting the assumptions of assertion (ii) of Theorem \thv(2.4.theo1)
to those of assertion (i) in the
statement
of assertion (i') above,
and leaving the definition of $\s$ unchanged,
the following holds:
For $R(\cdot)$ defined in \eqv(2.4.theo1.4),
$$
\wh S_n(\cdot)\Rrightarrow \wh R(\cdot)=\s+R(\cdot)\,,
\Eq(2.4.theo2.5)
$$
where $\s$ and $R(\cdot)$ are independent.
Moreover, \eqv(2.4.theo2.3)-\eqv(2.4.theo2.4) hold true with $\CC_{\infty}(t,s)$ defined through
$$
\CC_{\infty}(t,s)=\PP\left(\left\{R(k)\,,k\in\N\right\}
\cap (t, t+s)=\emptyset\right)\,,\quad 0\leq t<t+s\,.
\Eq(2.4.theo2.6)
$$

}

Thus,  when $a_n$ diverges, $\wh S_n(\cdot)$ converges to a delayed subordinator,
and it converges to  a delayed renewal process otherwise.

\remark As in assertion (ii) of Theorem \thv(1.3.theo2),
the statement that $\s$ and $S(\cdot)$ are independent in \eqv(2.4.theo2.1)
has the precise meaning that for each fixed $\o\in\O^{\t}$, $\s$ and $S(\cdot)$
are independent random variables on the probability space $(\O^{X}, \FF^{X}, \PP)$.
The same remark applies to the statement that $\s$ and $R(\cdot)$ in \eqv(2.4.theo2.5) are independent.


Specializing the previous theorem to the case where the initial distribution $\mu_n$ is the invariant
measure $\pi_n$ of the jump chain (see \eqv(4.4)) yields the following:

\corollary{\TH(2.4.cor1)} {\it Let $\mu_n=\pi_n$.
Under the assumptions of assertion (i) (respectively, assertion (ii)) of Theorem \thv(2.4.theo1),
w.r.t\. the same convergence mode as in \eqv(2.4.theo1.2) (equivalently,  \eqv(2.4.theo1.3), respectively, \eqv(2.4.theo1.4)),
$$
\lim_{n\rightarrow\infty}\CC_{n}(t,s)=\CC_{\infty}(t,s)\,,\quad 0\leq t<t+s\,,
\Eq(2.4.cor1.1)
$$
where $\CC_{\infty}(t,s)$ is defined in \eqv(1.3.theo1.3) (respectively, \eqv(2.4.theo2.6)).
}

Clearly, all sequences of initial distribution $\mu_n$ such that $\s=0$ in
Theorem \thv(2.4.theo2)
(i\.e\.,  all sequences of $\mu_n$'s such that Condition (A0) is satisfied with $F(v)=1$, $v\geq 0$)
give the same limiting time-time correlation function as the special choice $\mu_n=\pi_n$.
This is tantamount to the proof of the next corollary.

\corollary{\TH(2.4.cor2)} {\it Corollary \thv(2.4.cor1) remains valid for all sequences
of initial distribution $\mu_n$ such that,
w.r.t\. the same
convergence mode as in \eqv(2.4.theo1.2) (equivalently, \eqv(2.4.cor1.1)), for all $v\geq 0$,
$$
\lim_{n\rightarrow\infty}\sum_{x\in\VV_n}\mu_n(x)
\exp\{-v(r_n/\t(x))^{(1-a)}\}
=0\,,
\Eq(2.4.cor2.1)
$$
where $r_n=c_n^{1/(1-a)}$.
}

\remark
One may also interpret weak convergence of $S_n$ in the space $D^u([0,\infty))$
of c\`adl\`ag functions on $[0,\infty)$ equipped with the uniform topology.
Indeed if $\wt D([0,\infty))$ denotes the subspace of $D([0,\infty))$ consisting
of increasing paths having discontinuities at each integer times,
one easily sees that $\wt D([0,\infty))$ is a separable subspace of $D^u([0,\infty))$.

\proofof{Theorem \thv(2.4.theo1)}
The first assertion of Theorem \thv(2.4.theo1) is an elementary
specialization of Theorem \thv(1.3.theo1) to the asymmetric trap model on the complete graph.
Simply note that
$$
\sum_{j=1}^{\lfloor a_n t\rfloor}
\sum_{x\in\VV_n}p_n(J_n(j-1),x)
\exp\{-uc_n\l_n(x)\}
=
\frac{\lfloor a_n t\rfloor}{a_n}\nu_n(u,\infty)\,,
\Eq(1.3.theo1.5)
$$
where the r.h.s\. is chain independent. Thus, if $a_n$ is a diverging sequence,
\eqv(1.A1) and \eqv(1.A2) of Conditions (A1) and (A2)  of Theorem \thv(1.3.theo1)
reduce, respectively, to
$$
\nu_n(u,\infty)\rightarrow\nu(u,\infty)\,,
\Eq(1.AA1)
$$
$$
\frac{1}{a_n }\left[\nu_n(u,\infty)\right]^2\rightarrow0\,,
\Eq(1.AA2)
$$
as $n\rightarrow\infty$, and, clearly, \eqv(1.AA1) implies \eqv(1.AA2). Similarly,
\eqv(1.A1) of Condition (A3) becomes \eqv(2.4.theo1.2bis).

\remark Note that, setting
$
h_n(v)=\sum_{x\in\VV_n}p_n(v,x)e^{-uc_n\l_n(x)}
$,
\eqv(1.3.theo1.5) can be written as
$$
{\lfloor a_n t\rfloor}^{-1}\sum_{j=1}^{\lfloor a_n t\rfloor}h_n(J_n(j-1))
=
\sum_{y\in\VV_n}\pi_n(y)h_n(y)
=
E_{\pi_n}h_n(J_n(j-1))\,.
$$
In other words the sum appearing in Condition (A1) of Theorem \thv(1.3.theo1) is `ergodic'. A similar
observation holds for Condition (A2).

The novel part of  Theorem \thv(2.4.theo1) is  assertion (ii), whose elementary proof we now give.
Assume first that there exists a probability distribution $\nu$ such that, for all $u\geq 0$,
\eqv(2.4.theo1.2) holds in $\P$-probability.
Set $\xi_{n,i}=c_n^{-1}\l_n^{-1}(J_n(i))e_{n,i}$, $i\geq 0$. Putting $a_n=1$ in \eqv(1.3.2),
$
S_n(t)=\overline S_n(\lfloor t\rfloor)=\sum_{i=1}^{\lfloor t\rfloor}\xi_{n,i}
$.
Notice that for each $\o\in\O^{\t}$, $\{\xi_{n,i}, i\geq 1\}$ is an i.i.d\. sequence
on the probability space $(\O^{X}, \FF^{X}, \PP)$ since, by \eqv(4.6), the chain variables
$(J_n(i), i\in\N)$ form an i.i.d\. sequence, and since
$
\PP_{\mu_n}(\xi_{n,i}>u)=\nu_n(u,\infty)
$
does not depend on $i$. This means that $\overline S_n$ has stationary positive increments.
To prove \eqv(2.4.theo1.4) it thus suffices to prove that,
in $\P$-probability, for each integer $k$ (finite and independent of $n$),
$
\overline S_n(k)@>d>>R(k)
$
(see e.g\. \cite{Bi} p\. 30).
%
%
To this end consider  the Laplace transforms
$
\L_n(k,\theta)=\EE_{\mu_n} e^{-\theta  \overline S_n(k)}
$
and
$
\L(k,\theta)=\EE e^{-\theta R(k)}
$,
$\theta>0$.
From the assumption that, for all $u\geq 0$, \eqv(2.4.theo1.2) holds in $\P$-probability,
it follows that there exists a sequence $\wt\O^{\t}_n\subset\O^{\t}$ satisfying
$\lim_{n\rightarrow\infty}\P(\wt\O^{\t}_n)=1$,
and such that, for all large enough $n$,
$$
\sup_{u\geq 0}\left|\PP_{\mu_n}(\xi_{n,i}>u)-\nu(u,\infty)\right|=o(1)\,,\quad 1\leq i\leq n\,,
\Eq(1.3.theo1.6)
$$
for all $\o\in \wt\O^{\t}_n$.
Let now  $\o\in \wt\O^{\t}_n$ be fixed, where $n$ will be taken as large as needed. By independence,
$
\L_n(k,\theta)=\left(\EE e^{-\theta \xi_{n,i}}\right)^k
$.
From the integration by parts formula
$
\EE_{\mu_n}e^{-\theta\xi_{n,i}}
= 1-\theta\int_{0}^{\infty}e^{-\theta u}\PP_{\mu_n}(\xi_{n,i}>u)du
$,
it follows that
$$
\left|\EE_{\mu_n} e^{-\theta\xi_{n,i}}-\EE e^{-\theta\xi_{i}}\right|
\leq \sup_{u\geq 0}\left|\PP_{\mu_n}(\xi_{n,i}>u)-\nu(u,\infty)\right|\,.
\Eq(1.3.theo1.6')
$$
Thus, by \eqv(1.3.theo1.6), for all $n$ large enough, for each $k$,
$
\sup_{\theta>0}\left|\L_n(k,\theta)-\L(k,\theta)\right|=o(1)
$.
Now, by Feller's continuity theorem (see e.g\. \cite{Fe}, XIII.1, Theorem 2a), this implies that,
for all $n$ large enough, for each $k$,
$
\sup_{u>0}\left|\PP_{\mu_n}(\xi_{n,i}>u)-\PP(\xi_{i}>u)\right|=o(1)
$.
Since this holds true for each fixed $\o\in \wt\O^{\t}_n$,
it is tantamount to the statement that, for each $k$,
$
\overline S_n(k)@>d>>R(k)
$
in $\P$-probability.
The proof of assertion (ii) when \eqv(2.4.theo1.2) holds in $\P$-probability is now complete.
The proof in the case of $\P$-almost sure convergence is an elementary modification of it whose details we skip.
The proof of  Theorem \thv(2.4.theo1) is now done.
\endproof

\proofof{Theorem \thv(2.4.theo2)}
We first deal with assertion (i').
Eq\. \eqv(2.4.theo2.1) is proved just as \eqv(1.3.theo2.1) of Theorem \thv(1.3.theo2).
Assuming that for each $\o\in\O^{\t}$, $\s$ and $S(\cdot)$ in \eqv(2.4.theo2.1) are independent
random variables on the probability space $(\O^{X}, \FF^{X}, \PP)$,
\eqv(2.4.theo2.3) is proved in the same way as \eqv(1.3.theo1.4') of Theorem \thv(1.3.theo2),
and the special case $\s=0$ of \eqv(2.4.theo2.4) is nothing but \eqv(1.3.theo1.4).

Let us show that the above independence assumption is verified.
For this let $\o\in\O^{\t}$ be fixed. Note that by \eqv(4.6)
the jump chain $(J_n(i), i\in\N)$ becomes stationary in exactly one step. Namely,
for any initial distribution $\mu_n$, for all $i\geq 1$,
$
P_{\mu_n}(J_n(i)=x)=\pi_n(x)
$,
$x\in\VV_n$.
Thus, for each $n$, $\s_n$ and $\{\overline S_n(k)\,,k\geq 1\}$ in \eqv(1.3.2'')
are independent r\.v\.'s on $(\O^{X}, \FF^{X}, \PP)$.
This in turn implies that, for each $n$, $\s_n$ and $\{S_n(t)\,,t>0\}$
in the r\.h\.s\. of \eqv(1.3.2') are independent r\.v\.'s on $(\O^{X}, \FF^{X}, \PP)$.
Thus $\s$ and $S(\cdot)$ are independent, and since this is true for each $\o\in\O^{\t}$,
the claim follows.

We skip the proof of assertion (i''), which
is a re-run of the proof of assertion (i') (and, upstream from it,
of Theorem \thv(1.3.theo2)) in the simpler setting of discrete time process.
\endproof

\proofof{Corollary \thv(2.4.cor1)} 
Since
$
1-\PP_{\mu_n}(\s_n<v)=\sum_{x\in\VV_n}\mu_n(x)e^{-vc_n\l_n(x)}
$
(see e.g\. the proof of assertion (i) of Theorem \thv(1.3.theo2)) it follows from
\eqv(2.4.theo1.1) and the choice $\mu_n=\pi_n$ that
$$
1-\PP_{\mu_n}(\s_n<v)
=\frac{1}{a_n}\nu_n(u,\infty)\,.
\Eq(2.4.cor1.2)
$$
Suppose first that the assumptions of assertion (i) of Theorem \thv(2.4.theo1) are verified.
In view of \eqv(2.4.theo1.2) and \eqv(2.4.cor1.2),
$
1-\PP_{\mu_n}(\s_n<v)
\rightarrow 0
$
for all $v\geq 0$,
so that Condition (A0) is satisfied with $F(v)=1$, $v\geq 0$,
w.r.t\. the same convergence mode as in \eqv(2.4.theo1.2).
Eq\. \eqv(2.4.cor1.1) then follows from \eqv(2.4.theo2.4).
Suppose next that the assumptions of assertion (ii) of Theorem \thv(2.4.theo1)
are verified. Reasoning as above we readily see that
Condition (A0) is satisfied with $F(v)=\nu(u,\infty)$, $v\geq 0$,
w.r.t\. the same convergence mode as in \eqv(2.4.theo1.2).
Thus, by \eqv(2.4.theo2.5), the first increment $\s$ of the limiting renewal process
$\wh R$ has the same distribution as the inter-arrival times $\xi_k$ of $R$.
Hence, for all $0\leq t<t+s$,
$$
\PP\left(\left\{\s+R(k)\,,k\in\N\right\}\cap (t, t+s)=\emptyset\right)
=\PP\left(\left\{R(k)\,,k\in\N\right\}\cap (t, t+s)=\emptyset\right)
=\CC_{\infty}(t,s)\,,
\Eq(2.4.cor1.3)
$$
where the last equality is \eqv(2.4.theo2.6).
Since
$\s$ and $R(\cdot)$ in  \eqv(2.4.theo2.5) are independent,
we also have,  conditioning  on $\s$ and using \eqv(2.4.theo2.3), that,
for all $0\leq t<t+s$,
$$
\PP\left(\left\{\s+R(k)\,,k\in\N\right\}\cap (t, t+s)=\emptyset\right)
=1-F(t+s)+\int_{0}^t\CC_{\infty}(t-v,s)dF(v)
=\lim_{n\rightarrow\infty}\CC_{n}(t,s)\,.
\Eq(2.4.cor1.4)
$$
Equating the r.h.s\. of \eqv(2.4.cor1.2) to the r.h.s\. of \eqv(2.4.cor1.3) gives \eqv(2.4.cor1.1).
The proof of Corollary \thv(2.4.cor1) is done.
\endproof

\bigskip


\chap{3. Bouchaud's asymmetric trap model on the complete graph.}3

We now begin the investigation of Bouchaud's asymmetric trap model on the complete graph,
which will occupy the rest of the paper. Our aim here is to illustrate the scope and usefulness of
the abstract results of the previous sections
by solving, for the first time, a simple model of mean field
type which is not a time change of a simple random walk.
More realistic models and dynamics will be considered in the
companion papers \cite{G1} and \cite{G2}.

As in Bouchaud's symmetric ($a=0$) trap model on the complete graph encountered in Section 1,
the landscape of the asymmetric ($a>0$) model is made of i.i.d\. heavy tailed r.v\.'s.
This model appeared in \cite{BRM} where it was proposed and studied on various graphs.
The first rigorous results were obtained for the graph $\Z$  in \cite{BC1}.
There, it is shown that the time-time correlation function \eqv(1.1.8) does not exhibit an arcsine aging regime but
is subaging, and has the same ($a$-dependent) aging regime for all $a\in[0,1]$.
The very recent work \cite{BC} suggests  that on the contrary, on the graphs $\Z^d$, $d\geq 3$,
the asymmetry parameter, $a$, has no relevance on the aging phenomenon.
All these results contrast with the case of the complete graph
where the asymmetry parameter will be seen to trigger a dynamical phase transition.
More precisely, we will show that there exists a positive threshold value in $a$ below
which the model exhibits an ($a$-dependent) arcsine aging regime, whereas above it
arcsine aging is interrupted. Moreover, this phenomenon occurs ``on all time scales'',
i.e\. from time scale one up to, and including, the time scale of stationarity.
To make this picture complete we will show how, on the time scale of stationarity,
the model can be driven from an arcsine aging regime to its stationary regime.

In the rest of this section we describe the model
and some of its properties (Subsection 3.1) and state our main results, first, on the convergence
of the time-time correlation function  (Subsection 3.2), and next, on the clock process
(Subsection 3.3). All these results are obtained for a special choice of the initial
distribution $\mu_n$, namely for $\mu_n=\pi_n$.\note{
The role of the initial distribution will be investigated in a elsewhere.}
The proofs of these results will be presented in separate
sections
(Sections 4 to 7).

\bigskip
\line{\bf 3.1. The model.\hfill}

Let us now specify the model. Here $G_n(\VV_n,\EE_n)$ is the complete graph on $\VV_n=\{1,\dots,n\}$ that has a loop at
each vertex. The random landscape $(\t(x), x\in\VV_n)$ is a sequence of i\.i\.d\. random variables
whose distribution belongs to the domain of attraction of a positive stable
law with parameter $\a\in(0,1)$. This means that there exists a function $L$,
slowly varying at infinity, such that
$$
\P(\tau(x)>u)=u^{-\a}L(u)\,,\,\,\, u\geq 0\,.
\Eq(3.2.1)
$$
We denote by $\DD(\a)$ the class of all such probability distributions.
With a slight abuse of notation we write $\t\in\DD(\a)$ whenever \eqv(3.2.1) holds.
With these choices
Bouchaud's asymmetric trap model on the complete graph is the chain $X_{n}$ defined by
\eqv(3.1.2)-\eqv(3.1.3): given a parameter $0\leq a<1$,
the holding time parameters take the form\note{
For convenience we scale down the parameters $\l_n(x)$ in \eqv(3.1.2) by the factor $\sum_{x\in\VV_n}\t^a(x)$.
The same scaling is traditionally used in the symmetric case.
}
$$
\l_n(x)= (\t(x))^{-(1-a)}\,,\quad\forall\, x\in \VV_n\,\,,
\Eq(4.5)
$$
and the jump chain, $J_n$, has transition probabilities
$$
p_n(x,y)=\pi_n(y)\,,\quad\forall(x,y)\in \EE_n\,.
\Eq(4.6)
$$

It is easy to see that the chain $X_n$ has a unique reversible invariant measure, denoted by $\GG_{\a,n}$,
which is the Gibbs measure \eqv(1.1.2), that is
$$
\GG_{\a,n}(x)=\frac{\t(x)}{{\textstyle\sum_{x\in\VV_n}\t(x)}}\,,\quad x\in \VV_n\,.
\Eq(4.3)
$$
Clearly the jump chain $J_n$ also has a unique reversible invariant measure, $\pi_n$, given by
$$
\pi_n(x)=\frac{\t^a(x)}{{\textstyle\sum_{x\in\VV_n}\t^a(x)}}\,,\quad x\in \VV_n\,.
\Eq(4.4)
$$
Therefore $\pi_n$ is nothing but the Gibbs measure with parameter $\frac{\a}{a}$; namely,
$$
\pi_n=^d\GG_{\b,n}\,,\quad\beta\equiv\a/a\in(0,\infty)\,,
\Eq(4.6')
$$
where $=^d$ denotes equality in distribution.

Let us thus take a closer look at the Gibbs measure.
Its behavior changes at the critical value $\a=1$.
When $\a<1$
the order statistics of the Gibbs weights converges in distribution to
Poisson-Dirichlet distribution with parameter $\a$.
To formulate this result we need a little notation.
We use the abbreviation $\PRM(\l)$ for ``Poisson random measure with intensity measure $\l$''.
Let $\mu$ be the measure given by
$$
\mu(x,\infty)=x^{-\a}\,,\quad x>0\,,
\Eq(4.stat.1)
$$
and denote by $\{\g_k\}$ the marks of $\PRM(\mu)$
on $(0,\infty)$. Next denote by
$
\bar\g_{1}\geq \bar\g_{2}\geq\dots
$
the ranked Poisson marks. Then
Poisson-Dirichlet distribution with parameter $\a$
can be represented as
the distribution of the sequence
$$
\bar w_{1}\geq \bar w_{2}\geq\dots\text{where}\bar w_k=\frac{\bar\g_k}{\sum_{l}\bar\g_l}\,.
\Eq(4.stat.2)
$$
If we now label $\t_{n}(x(1))\geq\dots\geq\t_{n}(x(n))$ the landscape variables
arranged in decreasing order of magnitude, then, as $n\rightarrow\infty$,
$$
\bigl(\GG_{\a,n}(x(k))\bigr)_{k\geq 1}\overset d\to\longrightarrow\left(\bar w_k\right)_{k\geq 1}\,,
\Eq(4.stat.3)
$$
(see \cite{PY}, Proposition 10). This readily implies that
almost all the mass of the Gibbs measure is supported by the points $x(k)$ with largest
weights (i\.e\. with deepest traps).
In contrast, when $\a>1$, no single point carries a positive mass asymptotically. In particular,
it is not hard to show that
$$
\lim_{n\rightarrow\infty}\sup_{x\in\VV_n}\GG_{\a,n}(x)=0\text{in $\P$-probability.}
\Eq(4.stat.4)
$$
Here the Gibbs measure ``resembles a uniform measure''.
This dichotomy in the behavior of the Gibbs measure reflects the low
and high temperature regimes of the Gibbs measure of the REM,
the parameter $\a$ playing the role of the inverse of a temperature.

It is now easy to understand, on a heuristic level, that the chain $X_n$
undergoes a dynamical phase transition at the value $a=\a$.
Indeed by \eqv(4.6'), $\pi_n$ undergoes a static phase transition at
the value $a=\a$ which, in view of \eqv(4.6), will be reflected on the dynamics
as follows: when $a>\a$
the jump chain should resemble a symmetric random walk, and may explore the entire landscape;
in contrast, when $a<\a$ the jump chain will quickly go and visit a
trap of extreme depth from which it will not be able to escape, unless time is measured
on the scale of stationarity.

\bigskip
\line{\bf 3.2. Aging of $\CC_n(t,s)$.\hfill}

In this subsection we state our main results on the asymptotic behavior of
the time-time correlation function  $\CC_n(t,s)$ of
\eqv(1.1.8).
These results cover all choices of $a$ and $\a$ with $0<\a<1$, $0\leq a<1$, and $a\neq \a$,
and any choice of the time scale $c_n$ up to and including the time scale of stationarity.
We focus here on the case where the initial distribution $\mu_n$ is the invariant
measure $\pi_n$ of the jump chain.

To understand how the choice of $c_n$ affects $\CC_n(t,s)$ observe that scaling time
amounts to scaling the landscape. Namely, choosing the time scale $c_n$ the form
$
c_n=r_n^{1-a}
$,
we see from \eqv(4.5) that
$
(c_n\l_n(x))^{-1}= \bigl(\frac{\t(x)}{r_n}\bigr)^{(1-a)}
$.
This relation prompts us to call $r_n$ a {\it space scale}.
\note{
There should be no confusion between this and the volume or size, $n$, of the graph.
}
We will distinguish three types of space scales: the
{\it constant} scales (which simply are constant sequences),
the {\it intermediate}, and the {\it extreme} scales.

\definition{\TH(4.def1)}
We say that a positive and diverging sequence $r_n$ is:
\item{(i)} an {\it intermediate space scale} if there exists an increasing and diverging sequence $b_n>0$ such that
$$
\frac{b_n}{n}=o(1)  \text{and} b_n\P(\t(x)\geq r_n)\sim 1\,,
\Eq(4.8)
$$
\item{(ii)} an {\it extreme space scale} if there exists an increasing and diverging sequence $0<b_n\leq n$ such that
$$
\frac{b_n}{n}\sim 1 \text{and} b_n\P(\t(x)\geq r_n)\sim 1\,.
\Eq(4.9)
$$

As the next lemma shows, these scales are well separated:

\lemma{\TH(4.lemma4)} {\it Let $r_n^{cst}$, $r_n^{int}$ and $r_n^{ext}$ denote, respectively,
a constant, an intermediate and an extreme space scale. Then
$
r_n^{cst}\ll r_n^{int}\ll r_n^{ext}
$.
}

The proof of Lemma \thv(4.lemma4) is postponed to the end of Subsection 7.1.
We now present three Theorems. The first of them establishes that if $a<\a$ then
$\CC_{n}(t,\rho t)$
exhibits aging of arcsine type on all time scales.

\newpage

\theo{\TH(4.theo3)} [Arcsine aging regime.] {\it Assume that $a<\a$. Let $c_n=r_n^{1-a}$ and take $\mu_n=\pi_n$.

\item{(i)}
If $r_n$ is a constant space scale then, $\P$-almost surely, for all $\rho>0$,
$$
\lim_{t\rightarrow\infty}\lim_{n\rightarrow\infty}\CC_{n}(t,\rho t)=\asl_{\sfrac{\a-a}{1-a}}(1/1+\rho)\,.
\Eq(4.theo3.1)
$$

\item{(ii)}
If $r_n$ is an intermediate space scale then, in $\P$-probability, for all $t\geq 0$ and all $\rho>0$,
$$
\lim_{n\rightarrow\infty}\CC_{n}(t,\rho t)=\asl_{\sfrac{\a-a}{1-a}}(1/1+\rho)\,.
\Eq(4.theo3.2)
$$
(When $a=0$
this statement
holds $\P$-almost surely whenever $\frac{b_n}{n}\log n=o(1)$.)

\item{(iii)}
If $r_n$ is an extreme space scale then, for all $\rho>0$,  in $\P$-probability,
$$
\lim_{t\rightarrow 0+}\lim_{n\rightarrow\infty}\CC_{n}(t,\rho t)=\asl_{\sfrac{\a-a}{1-a}}(1/1+\rho)\,.
\Eq(4.theo3.3)
$$
}

\remark One can show that the above statements remain valid when $\a=1$
with the arcsine density replaced by the delta mass at 1. Thus the time correlation
function vanishes in the limit. The precise asymptotics of this decay will
be studied elsewhere.

The next theorem shows that when $a>\a$, none of the time scale and limiting procedures
of Theorem \thv(4.theo3) yields aging:

\theo{\TH(4.theo2)} [Stranded in deep traps.] {\it Assume that $a>\a$. Let $c_n=r_n^{1-a}$ and take $\mu_n=\pi_n$.

\item{(ii)} If $r_n$ is a constant or intermediate space scale then, for all $0\leq t<t+s$,
$$
\lim_{n\rightarrow\infty}\CC_{n}(t,s)=1
\text{in $\P$-probability.}
\Eq(4.theo2.1)
$$

\item{(ii)}
If $r_n$ is an extreme space scale then, for all $\rho>0$,
$$
\lim_{t\rightarrow 0}\lim_{n\rightarrow\infty}\CC_{n}(t,\rho t)=1
\text{in $\P$-probability.}
\Eq(4.theo2.2)
$$
}

At a heuristic level Theorem \thv(4.theo2) is easy to understand.
For $a>\a$ the initial distribution $\mu_n$ behaves like a ``low temperature''
Gibbs measure, namely $\mu_n=^d\GG_{\b,n}$, $\beta=\a/a<1$.
This means that almost all its mass is carried by
by traps whose size is of the order of extreme space scales. Now the mean waiting time
in such deep traps diverges as as $n$ diverges whenever time is measured on a scale
which is small compared to the extreme scales: the chain gets stranded.

The last theorem below is valid for all $0\leq a<1$. It states that, as expected, on extreme time scales,
taking the infinite volume limit first,
the process reaches stationarity as $t\rightarrow\infty$.
As before let $\{\g_k\}$ denote the marks of $\PRM(\mu)$ on $(0,\infty)$, and define
$$
{\CC}^{sta}_{\infty}(s)=
\sum_{k}\frac{\g_k}{\sum_{l}\g_l}e^{-s\g_k^{-(1-a)}}
\,,\quad s\geq 0\,.
\Eq(4.theo4.0)
$$

\theo{\TH(4.theo4)} [Crossover to stationarity.] {\it   Let $c_n=r_n^{1-a}$ where
$r_n$ is an extreme space scale. The following holds for all $0\leq a<1$, $a\neq \a$:

\item{(i)} If  $\mu_n=\GG_{\a,n}$ then, for all $s\leq t<t+s$,
$$
\lim_{n\rightarrow\infty}{\CC}_n(t,s)\overset{d}\to={\CC}^{sta}_{\infty}(s)\,,
\Eq(4.theo4.1)
$$
where $\overset{d}\to=$ denotes equality in distribution.

\item{(ii)} If  $\mu_n=\pi_n$, for all $s>0$,
$$
\lim_{t\rightarrow\infty}\lim_{n\rightarrow\infty}{\CC}_n(t,s)\overset{d}\to
={\CC}^{sta}_{\infty}(s)\,.
\Eq(4.theo4.2)
$$
}

Comparing \eqv(4.theo3.3) and  \eqv(4.theo4.2) we see that
when
time goes from $0$ to $\infty$, for $a<\a$, the chain moves out of
an arcsine aging regime and crosses over to its stationary regime.
Aging is then interrupted.

\bigskip
\line{\bf 3.3. Convergence of the clock process.\hfill}

In this subsection we state  the convergence properties of the pure clock process $S_n$ of \eqv(1.3.2) from which
the asymptotic behavior of the time-time correlation function will later be deduced
(see the proofs of Theorem \thv(4.theo3), \thv(4.theo2) and \thv(4.theo4) in Section 7).
These properties will themselves be deduced from Theorem \thv(2.4.theo1) whose notations we now use.
We group them according to the choice of  space scale (constant, intermediate, and extremal),
and show that for each space scale the nature of the limiting clock process changes at the critical
value $a=\a$.

\proposition{\TH(4.prop1)} [Constant scales.] {\it  Take $\mu_n=\pi_n$. If If $r_n$ is a constant space scale then the following holds:
set
$
\nu^{cst,+}=\d_{\infty}
$,
and, for $a<\a$ and $\t\in\DD(\a)$, let $\nu^{cst,-}$ be the measure
on $(0,\infty)$ defined through
$$
\nu^{cst,-}(u,\infty)=\E\frac{\t^a}{\E\t^a} e^{-u/\t^{(1-a)}}\,,\quad u>0\,.
\Eq(4.prop1.2)
$$
\item{(i)} If $a<\a$ then
$
S_n(\cdot)\Rrightarrow R^{cst,-}(\cdot)
$
$\P$-almost surely, where $R^{cst,-}$ is the renewal process of inter-arrival distribution $\nu^{cst,-}$.
\item{(ii)} If $a>\a$ then,
$
S_n(\cdot)\Rrightarrow R^{cst,+}(\cdot)
$
in $\P$-probability, where $R^{cst,+}$ is the degenerate renewal process of inter-arrival distribution $\nu^{cst,+}$.
}

The lemma below shows that $\nu^{cst,-}(u,\infty)$ is regularly varying at infinity with index $\sfrac{\a-a}{1-a}$.

\lemma{\TH(4.lemma1)} {\it If $a<\a$
then $\int_{0}^{\infty}\nu^{cst,-}(u,\infty)du=\infty$ and
$$
\nu^{cst,-}(u,\infty)=u^{-\frac{\a-a}{1-a}}\ell(u)\G\bigl(\sfrac{\a-a}{1-a}\bigr)/\E\t^a
\,,\quad u\rightarrow\infty\,,
\Eq(4.lem1.1)
$$
where $\ell(u)$ is slowly varying at infinity (here $\G$ denotes the gamma function.)
}

\proposition{\TH(4.prop2)} [Intermediate scales.] {\it Take $\mu_n=\pi_n$.
If $r_n$ is an intermediate space scale then the following holds:
set
$
\nu^{int,+}=\d_{\infty}
$,
and, for $a<\a$, let $\nu^{int,-}$ be the measure
on $(0,\infty)$ defined through
$$
\nu^{int,-}(u,\infty)=u^{-\frac{\a-a}{1-a}}\sfrac{\a}{1-a}\G\bigl(\sfrac{\a-a}{1-a}\bigr)/\E\t^a\,,\quad u>0\,.
\Eq(4.prop2.2)
$$

\item{(i)} If $a<\a$ then,
$
S_n(\cdot)\Rightarrow S^{int,-}(\cdot)
$
in $\P$-probability, where $S^{int,-}$ is the stable subordinator of L\'evy measure $\nu^{int,-}$.
\item{(ii)} If $a>\a$ then,
$
S_n(\cdot)\Rrightarrow R^{int,+}(\cdot)
$
in $\P$-probability, where $R^{int,+}$
is the degenerate renewal process of inter-arrival distribution $\nu^{int,+}$.
}

To formulate the results on extreme scales recall that for $\mu$ defined in \eqv(4.stat.1),
$\{\g_k\}$ denote the marks of $\PRM(\mu)$ on $(0,\infty)$, and
introduce the re-scaled landscape variables:
$$
\g_n(x)=r_n^{-1}\t(x)\,\quad x\in\VV_n\,.
\Eq(4.prop3.1)
$$

\proposition{\TH(4.prop3)} [Extreme  scales.] {\it
If $r_n$ is an extreme space scale then both the sequence of re-scaled landscapes
$(\g_n(x),\, x\in\VV_n)$,
$n\geq 1$, and the marks of $\PRM(\mu)$ can be represented on a common probability space
$(\O, \FF, \bold P)$ such that, in this representation,
denoting by $\bold S_n$ the process \eqv(1.3.2),
the following holds.
For $a<\a$, resp\. $a>\a$, let $\nu^{ext,-} $, resp\. $\nu^{ext,+} $ be
the random measures on $(0,\infty)$ defined on $(\O, \FF, \bold P)$ through
$$
\eqalign{
\nu^{ext,-} (u,\infty)&=\frac{1}{\E\t^a}\sum_{k}\g_k^ae^{-u\g_k^{-(1-a)}}\,,\quad u>0\,,
\cr
\nu^{ext,+} (u,\infty)&=\sum_{k}\frac{\g_k^a}{\sum_{l}\g_l^a}e^{-u\g_k^{-(1-a)}}\,,\quad u>0\,.
}
\Eq(4.prop3.2)
$$
Take $\mu_n=\pi_n$. Then, $\bold P$-almost surely,
$$
\eqalign{
&\bold S_n(\cdot)\Rightarrow S^{ext,-}(\cdot)\text{if}a<\a\,,
\cr
&\bold S_n(\cdot)\Rrightarrow R^{ext,+}(\cdot)\text{if}a>\a\,,
}
\Eq(4.prop3.4)
$$
where $S^{ext,-}$ is the subordinator of L\'evy measure $\nu^{ext,-}$,
and $R^{ext,+}$ is the renewal process of inter-arrival distribution $\nu^{ext,+}$.
}

\lemma{\TH(4.lemma2)}{\it Let $\nu^{ext,\pm}$ be as in \eqv(4.prop3.1) and define
$
m^{\pm}=\int_{0}^{\infty}\nu^{ext,\pm} (u,\infty)du
$.
We have
$$
\eqalign{
m^-&=
\frac{\sum_{k}\g_k}{\E\t^a}\text{if}a<\a\,,
\cr
m^+&=\frac{\sum_{k}\g_k}{\sum_{k}\g_k^a}\text{if}a>\a\,,
}
\Eq(4.lem2.2)
$$
and each of the sums appearing in the right hand side is $\bold P$-almost surely finite. Moreover
$$
\frac{1}{m^-}\int_{0}^{\infty}\nu^{ext,-}(u,\infty)du=\frac{1}{m^+}\int_{0}^{\infty}\nu^{ext,+}(u,\infty)du
={\CC}^{sta}_{\infty}(s)\,,
\Eq(4.lem2.2')
$$
for ${\CC}^{sta}_{\infty}$ defined as in \eqv(4.theo4.0).
}

Here the subordinator $S^{ext,-}$ is not stable.
However $\nu^{ext,-} (u,\infty)$ is regularly varying at $0^+$ with index $\sfrac{\a-a}{1-a}$:

\lemma{\TH(4.lemma3)}{\it If $a<\a$ then, $\bold P$-almost surely,
$$
\nu^{ext,-} (u,\infty)\sim
u^{-\frac{\a-a}{1-a}}
\sfrac{\a}{1-a}\G\bigl(\sfrac{\a-a}{1-a}\bigr)\,,
\quad u\rightarrow 0^+\,.
\Eq(4.lem3.1)
$$
}

\lemma{\TH(4.lemma5)}{\it If $a>\a$ then, $0\leq \nu^{ext,+}(u,\infty)\leq 1$ and, $\bold P$-almost surely,
$$
\lim_{u\rightarrow 0}\nu^{ext,+}(u,\infty)=1\,,\quad \lim_{u\rightarrow\infty}\nu^{ext,+}(u,\infty)=0\,.
\Eq(4.lem5.1)
$$
}

The rest of the paper is organized as follows.
The proofs of Proposition \thv(4.prop1), (i), Proposition  \thv(4.prop2), (i),
and Proposition  \thv(4.prop3), which rely on very different tools, are presented
in three separate sections (Section 4,  5, and 6 respectively).
The proofs of Theorems \thv(4.theo3), \thv(4.theo2), and \thv(4.theo4) as well as
Proposition \thv(4.prop1), (ii), Proposition  \thv(4.prop2), (ii),
which all use the results of Section 6, are gathered in Section 7.

\bigskip


\chap{4. Constant scales.}4

In this short section we prove Lemma \thv(4.lemma1) and the first assertion of Proposition \thv(4.prop1).
The proof of the latter relies on assertion (ii) of Theorem \thv(2.4.theo1) on
the convergence of the clock process $S_n$ to a process of partial sums.
The proof of the second assertion of Proposition \thv(4.prop1),
which relies on results from Section 6, is postponed to Section 7.

\proofof{Proposition \thv(4.prop1), (i)}
Take $r_n=a_n=1$
in \eqv(2.4.theo1.1) and $\nu=\nu^{cst,-}$ in \eqv(2.4.theo1.2).
Note that if $\t\in\DD(\a)$ then $\E\t^a<\infty$ for all $a<\a$,
so that $\E\t^ae^{-u/\t^{(1-a)}}\leq \E\t^a<\infty$ for all $u\geq 0$.
Thus, for  all $u\geq 0$, the strong law of large numbers applies to both the numerator and denominator of
\eqv(2.4.theo1.1), yielding
$
\lim_{n\rightarrow\infty}\nu_n(u,\infty)=\nu^{cst,-}(u,\infty)
$
$\P$-almost surely.  Together with the monotonicity of $\nu_n(u,\infty)$ and the continuity of the limiting function
$\nu^{cst,-}(u,\infty)$ in $u$, this entails that there exists a subset $\O^{\t}_1\subset\O^{\t}$ of
the sample space $\O^{\t}$ of the $\t$'s with the property that $\P(\O^{\t}_1)=1$, and such that, on $\O^{\t}_1$,
$$
\lim_{n\rightarrow\infty}\nu_n(u,\infty)=\nu^{cst,-}(u,\infty)
\,,\quad\forall\,u\geq 0\,.
\Eq(4.prop1.3)
$$
The conditions of assertion (ii) of Theorem \thv(2.4.theo1) are thus satisfied $\P$-almost surely.
Assertion (i) of Proposition \thv(4.prop1) is proven.\endproof

\proofof{Lemma \thv(4.lemma1)} We assume throughout that $a<\a$. For $u\geq 0$ and $y\geq 0$ set
$
\vp_u(y)=y^ae^{-u/y^{(1-a)}}
$.
Integrating by parts,
$
\E\vp_u(\t)
=\int_0^\infty\vp'_u(x)\P(\t>x)dx
$.
Performing the change of variable $x=u^{1/(1-a)}y$ and noting that $\vp'_u\left(u^{1/(1-a)}y\right)=u^{-1}\vp'_1(y)$, we get
$
\E\vp_u(\t)
=
u^{\frac{a}{1-a}}
\int_0^\infty\vp'_1(y)\bigl(\t>u^{\frac{1}{1-a}}y\bigr)dy
$.
Since, for $\t\in\DD(\a)$,
$
u^{\frac{\a}{1-a}}\P\bigl(\t>u^{\frac{1}{1-a}}\bigr)=\ell(u)
$
for some function
$\ell(u)$ that varies slowly at infinity, we further get
$$
u^{\frac{\a-a}{1-a}}\E\vp_u(\t)
=\ell(u)
\int_0^\infty\vp'_1(y)
\bigl[\P\bigl(\t>u^{\frac{1}{1-a}}y\bigr)/\P\bigl(\t>u^{\frac{1}{1-a}}\bigr)\bigr]
dy\,.
\Eq(4.lem1.2)
$$
Next, for $\t\in\DD(\a)$,
$
{\P\bigl(\t>u^{\frac{1}{1-a}}y\bigr)}/{\P\bigl(\t>u^{\frac{1}{1-a}}\bigr)}\rightarrow y^{-\a}
$
as $u\rightarrow\infty$.
Because of the monotonicity the approach is uniform, and so
$$
\int_0^\infty\vp'_1(y)
\bigl[\P\bigl(\t>u^{\frac{1}{1-a}}y\bigr)/\P\bigl(\t>u^{\frac{1}{1-a}}\bigr)\bigr]
dy
\rightarrow
\int_0^\infty \vp'_1(y)y^{-\a}dy
=\G\bigl(\sfrac{\a-a}{1-a}\bigr)\,,
\Eq(4.lem1.3)
$$
as
$u\rightarrow\infty$.
Combining \eqv(4.lem1.2) and \eqv(4.lem1.3)  yields \eqv(4.lem1.1) of Lemma \thv(4.lemma1).
From this and the assumption that $a<\a$ the claim that
$\int_{0}^{\infty}\nu^{cst,-}(u,\infty)du=\infty$
readily follows. The lemma is proven.\endproof

\bigskip


\chap{5. Intermediate scales.}5

In this section we prove the first assertion of Proposition \thv(4.prop2)
using assertion (i) of Theorem \thv(2.4.theo1).
(The proof of the second assertion of Proposition \thv(4.prop2),
which relies on results from Section 6, will be given in Section 7.)
The key ingredient of the proof is Proposition \eqv(6.prop1) below
that establishes control on the quantity $\nu_n(u,\infty)$ from \eqv(2.4.theo1.1).

\proposition{\TH(6.prop1)}{\it Let $r_n$ be an intermediate space scale and
choose $a_n\sim r_n^{-a}b_n$ in \eqv(2.4.theo1.1).
Assume that $a<\a$ and let $\nu^{int,-}$ be defined in \eqv(4.prop2.2). There exists a
sequence of subsets $\O^{\tau}_{1,n}\subset\O^{\tau}$
with $\P(\O^{\tau}_{1,n})\geq 1-o(1)$ such that for all $n$ large enough, on $\O^{\tau}_{1,n}$, the following holds
for all $u>0$:
$$
\left|\nu_{n}(u,\infty)-m_n(u)\right|<(b_n/n)^{1/3}\s(u)\,,
\Eq(6.prop1.1)
$$
where $m_n$ is a sequence of positive decreasing functions that satisfy,
$$
\lim_{n\rightarrow\infty}m_n(u)=\nu^{int,-}(u,\infty)\,,
\Eq(6.lem1.2)
$$
and where
$
\s^2(u)= c_0+c_1u^{-1+c_2}\
$
for some constants $0\leq  c_0,c_1<\infty$ and $0<c_2\leq 1$ that depend only on $\a$ and $a$.
}


The proof of Proposition \thv(6.prop1) relies on a weak law of large numbers for a triangular array which we now introduce.
For fixed $u\geq 0$ set
$$
\vp_u(y)=y^ae^{-u/y^{(1-a)}}\,,\quad y\geq 0\,,
\Eq(6.lem1.3)
$$
and, denoting by $\g_n(x)=r_n^{-1}\t(x)$, $x\in\VV_n$, the re-scaled landscape variables,
consider the array of row-wise independent random variables, $\{Z_{n,u}(x), x\in\VV_n, n\geq 1\}$, defined by
$$
Z_{n,u}(x)=(b_n/n)\vp_u(\g_n(x))\,.
\Eq(6.lem1.4)
$$
With this notation $\nu_{n}(u,\infty)=\sum_{x\in\VV_n}Z_{n,u}(x)$.

\lemma{\TH(6.lemma1)}{\it
Given $\rho>0$ set $m_n(u)=\sum_{x\in\VV_n}\E Z_{n,u}(x)\1_{\{Z_{n,u}(x)<\rho\}}$. Then,
under the assumptions of Proposition \thv(6.prop1), the following holds for each fixed $u>0$:
moreover $m_n(u)$ satisfies \eqv(6.lem1.2); moreover
for all large enough $n$, there exists constants $0<c_0,c_1,c_3,c_4<\infty$ and $0< c_2\leq 1$, such that
$$
\P\left(\left|\sum_{x\in\VV_n}Z_{n,u}(x)-m_n(u)\right|\geq z\right)
\leq \frac{1}{z^2}\left(\frac{b_n}{n}\right)\left(c_0+c_4 u^{-1+c_2}\right)+c_3\left(\frac{b_n}{n}\right)^{c_1}\,.
\Eq(6.lem1.1)
$$
%
}

\remark When $a=0$  one may prove, using a classical exponential Tchebychev inequality, that
$
\P\bigl(\bigl|\nu_n(u,\infty)-\E[\nu_n(u,\infty)]\bigr|
\geq 2\sqrt{a_nL/n}\sqrt{\E[\nu_n(2u,\infty)]}\bigr)
\leq e^{-L}
$
for all $L\geq 0$ such that $a_nL/n=o(1)$,
and that
$
\lim_{n\rightarrow\infty}\E[\nu_n(u,\infty)]=\nu^{int,-}(u,\infty)
$.
Using this one may show that the first assertion of Proposition \thv(4.prop2) holds
$\P$-a\.s\. if $(b_n/n)\log n=o(1)$ and in $\P$-probability otherwise. See \cite{G2}
for the details of the proof in the analogous case of the REM dynamics.


Let us collect here the information on
the slow variation properties of the function $\vp_u$ and its inverse
which will be needed in the proof of Lemma \thv(6.lemma1).
We use the notations of Appendix A.3 on regular variations.
Assume that $a>0$.
Clearly $\vp_u$ is strictly increasing and $\vp_u\in R_a$. Thus
$\vp_u^{-1}$ is well defined, strictly increasing, and, by Lemma \thv(A.3.lemma3),
 $\vp_u^{-1}\in R_{1/a}$.
Let the functions $\phi_u$  be defined through
$$
\vp_u^{-1}(y)=y^{1/a}\phi_u(y)\,.
\Eq(A.3.lemma4.2)
$$
Then $\phi_u\in R_0$. The following elementary lemma, stated without proof, gives its explicit form for small $u$.

\lemma{\TH(A.3.lemma4)}{\it Let $0<a<1$. Then, for all $u>0$,
(i) $\phi_u(y)\geq 1$ for all $y\geq 0$,
and
(ii) $\phi_u(y)\leq e^{1/a}$ for all $y\geq u^{1/(1-a)}$. Moreover,
(iii) if $u\geq v$ then $\phi_u(y)\leq \phi_v(y)$ for all $y\geq 0$.
}

\proofof{Lemma \thv(6.lemma1)}
The cases $a=0$ and $a>0$ will be treated separately. Assume first that $a>0$.
As in the proof of the weak law of large numbers for triangular arrays
of independent random variables, our first step consists in writing that,
given $\rho>0$,
$$
\P\left(\left|\sum_{x\in\VV_n}Z_{n,u}(x)-m_n(u)\right|\geq z\right)
\leq
\sum_{x\in\VV_n}\P(Z_{n,u}(x)>\rho)
+\frac{1}{z^2}\sum_{x\in\VV_n}Var\left(Z_{n,u}(x)\1_{\{Z_{n,u}(x)<\rho\}}\right)
\Eq(6.lem1.1')
$$
(see e\.g\. \cite{F}, Section VII.7).
We will now show that for each truncation level satisfying  $\rho\geq 1$, the right hand side of \eqv (6.lem1.1')
goes to zero as $n\rightarrow\infty$.

In what follows, $c_i$, $i\geq 1$, designate finite positive constants that may depend on
the parameters $a$ and $\a$, but not on $u$, and whose value may change from line to line.
We sometimes write $\vp$, $\phi$, and $Z_n(x)$ instead $\vp_u$, $\phi_u$, and $Z_{n,u}(x)$
if no confusion may arise.
We begin by establishing that there exists $c_1>0$ such that, for all large enough $n$,
$$
\sum_{x\in\VV_n}\P(Z_{n,u}(x)>\rho)\leq
2\rho^{-(1+c_1)}(b_n/n)^{c_1}\,,\quad \forall u>0\,.
\Eq(6.lem1.5)
$$
By definition of $Z_n(x)$, for $\t(x)\in\DD(\a)$,
$$
\eqalign{
\P(Z_{n,u}(x)>\rho)
&=\P\left(\t(x)>r_n \vp^{-1}\bigl(\sfrac{n}{b_n}\rho\bigr)\right)
\cr
&=\left(r_n \vp_u^{-1}\bigl(\sfrac{n}{b_n}\bigr)\right)^{-\a}L\left(r_n \vp_u^{-1}\bigl(\sfrac{n}{b_n}\rho\bigr)\right)
\cr
&=\left(r_n \bigl(\sfrac{n}{b_n}\rho\bigr)^\frac{1}{a}\phi\bigl(\sfrac{n}{b_n}\rho\bigr)\right)^{-\a}
L\left(r_n \bigl(\sfrac{n}{b_n}\rho\bigr)^\frac{1}{a}\phi\bigl(\sfrac{n}{b_n}\rho\bigr)\right)\,,
}
\Eq(6.lem1.7)
$$
where $\phi(y)$ is defined in \eqv(A.3.lemma4.2). Using furthermore that
$r_n^{\a}\P(\t(x)>r_n)=L(r_n)$, we obtain
$$
\eqalign{
\sum_{x\in\VV_n}\P(Z_n(x)>\rho)&=\rho^{-\frac{\a}{a}}\left(\frac{b_n}{n}\right)^{\frac{\a}{a}-1}\left[b_n\P(\t(x)>r_n)\right]
\frac
{
L\left(r_n \bigl(\sfrac{n}{b_n}\rho\bigr)^\frac{1}{a}\phi\bigl(\sfrac{n}{b_n}\rho\bigr)\right)
}{
\phi\bigl(\sfrac{n}{b_n}\rho\bigr)^{\a}L(r_n)
}\,.
}
\Eq(6.lem1.8)
$$
By definition of intermediate space scales the sequence $b_n$ satisfies
$b_n\P(\t(x)\geq r_n)\sim 1$ and $\frac{b_n}{n}=o(1)$.
To control the quotient in the r\.h\.s\. of \eqv(6.lem1.7) we use the bound
$\phi_u\bigl(\sfrac{n}{b_n}\rho\bigr)\geq 1$
of Lemma \thv(A.3.lemma4) (valid for all $u\geq 0)$ together with Lemma \thv(A.3.lemma2) to deduce that
there exist positive sequences
$\e_n$ and $\d_n$ that verify $\e_n\downarrow 0$, $\d_n\downarrow 0$ as $n\uparrow\infty$
and such that, for all $n$ large enough,
$$
\frac
{
L\left(r_n \bigl(\sfrac{n}{b_n}\rho\bigr)^\frac{1}{a}\phi\bigl(\sfrac{n}{b_n}\rho\bigr)\right)
}{
\phi\bigl(\sfrac{n}{b_n}\rho\bigr)^{\a}L(r_n)
}
\leq
(1+\d_n)\left(\frac{b_n}{n}\right)^\frac{\e_n}{a}\,.
\Eq(6.lem1.9)
$$
Inserting \eqv(6.lem1.9) in \eqv(6.lem1.8) we get that
$
\sum_{x\in\VV_n}\P(Z_n(x)>\rho)
\leq
2\rho^{-\frac{\a-\e_n}{a}}(b_n/n)^{\frac{\a-\e_n}{a}-1}
$
where, since $\a/a>1$, all $n$ large enough, $(\a-\e_n)/a-1>c_1$ for some $c_1>0$.
Notice that this bound is independent of $u$ and holds true for all $u>0$.
This establishes \eqv(6.lem1.5).

\noindent
We now turn to the variance term in \eqv(6.lem1.1').
We will establish that for all $u>0$ the following holds: for all large enough $n$,
there exist constants $0\leq  c_0,c_1, c_3, c_4<\infty$ and $0<c_2\leq 1$,
that depend only on $\a$ and $a$, and such that,
$$
\sum_{x\in\VV_n}Var\left(Z_{n,u}(x)\1_{\{Z_{n,u}(x)<\rho\}}\right)\leq
b_n/n\left(c_0+c_4 u^{-1+c_2}\right)+c_3(b_n/n)^{c_1}\,.
\Eq(6.lem1.6)
$$
Clearly the left hand side of \eqv(6.lem1.6) is bounded above by
$
\sum_{x\in\VV_n}\E Z^2_n(x)\1_{\{Z_n(x)<\rho\}}
$. By the integration by parts formula for truncated random variables,
$$
\eqalign{
\sum_{x\in\VV_n}\E Z^2_n(x)\1_{\{Z_{n,u}(x)<\rho\}}
\leq &
2n\left(\frac{b_n}{n}\right)^2\int_0^{\rho\frac{n}{b_n}}y\P(\vp_u(\g_n(x))>y)dy
\cr
=&
2\frac{b_n}{n}[b_n\P(\t(x)>r_n)]\int_0^{\vp_u^{-1}(\rho\frac{n}{b_n})}\vp_u(z)\vp_u'(z)h_n(z)dz\,,
\cr
}
\Eq(6.lem1.11)
$$
where, in the last line, we performed the change of variable $y=\vp_u(z)$, and defined
$$
h_n(z)=\frac{\P(\t(x)>r_nz)}{\P(\t(x)>r_n)}\,.
\Eq(6.lem1.12)
$$
To further express the integral in the last line of \eqv(6.lem1.11) we split it into
$I_n(u)=I_n'(u)+I_n''(u)$,
$$
I_n'(u)=2\int_0^{1}\vp_u(z)\vp_u'(z)h_n(z)dz
\,,\quad
I_n''(u)=2\int_1^{\vp_u^{-1}(\frac{n}{b_n}\rho)}\vp_u(z)\vp_u'(z)h_n(z)dz\,.
\Eq(6.lem1.13)
$$
To deal with $I_n'(u)$ we use that
$
h_n(z)\rightarrow z^{-\a}
$,
$n\rightarrow\infty$,
where the convergence is uniform in $z$ as $0\leq z\leq 1$, since for each $n$,
$h_n(z)$ is a monotone function, and since the limit, $z^{\a}$, is continuous.
Thus for all $\e>0$ there exists $n(\e)$ such that for all $n\geq n(\e)$
$$
\left|I'_n(u)-2\int_0^{1}\frac{\vp_u(z)\vp'_u(z)}{z^\a}dz\right|
\leq \e2\int_0^{1}\vp_u(z)\vp'_u(z)dz
\leq \frac{\e}{(2a+1)}\,,\quad \forall u>0\,.
\Eq(6.lem1.15)
$$
%
%
Integrating by parts,
$
2\int_0^{1}\frac{\vp_u(z)\vp'_u(z)}{z^\a}dz=
\vp_u^2(1)+\a\int_0^{1}\frac{\vp_u^2(z)}{z^{1+\a-\e_n}}dz
$.
Performing the change of variable $z=y^{-1/(1-a)}$, the last integral may be written as,
$
i(u)=\frac{\a}{1-a}\int_{1/2}^{\infty}z^{\b-1}e^{-2uz}dz
$
where we set $\b:=\frac{\a-2a}{1-a}$. Now, if $\b>0$,
$
i(u)
\leq\frac{\a}{1-a}\int_{0}^{\infty}z^{\b-1}e^{-uz}dz
=(2u)^{-\b}\frac{\a}{1-a}\G(\b)
$
where $-\b+1>0$, whereas if $\b\leq0$,
$
i(u)
\leq\frac{\a}{1-a}\int_{1/2}^{\infty}z^{-|\b|-1}dz
=\frac{\a}{\b(1-a)}2^{|\b|}
$.
Combining these observations with \eqv(6.lem1.15) we conclude that for all $u>0$ the following holds:
for all large enough $n$ there exist constants $0\leq  c_0,c_4<\infty$ and $0<c_2\leq 1$, that depend only on $\a$ and $a$,
and such that
$$
I'_n(u)\leq c_0+c_4 u^{-1+c_2}\,.
\Eq(6.lem1.15')
$$

To bound $I_n''(u)$ we note that
$
h_n(z)=x^{-\a}(L(r_nz)/L(r_n))
$
and use that by Lemma \thv(A.3.lemma2), for each $x>1$ and large enough $n$,
$$
(1-\d_n)z^{-\a-\e_n}\leq h_n(z)\leq (1+\d_n)z^{-\a+\e_n}\,,
\Eq(6.lem1.17)
$$
for some positive sequences $\e_n$ and $\d_n$ satisfying $\e_n\downarrow 0$, $\d_n\downarrow 0$ as $n\uparrow\infty$.
Thus
$$
I''_n(u)\leq 2 (1+\d_n)\int_1^{\vp_u^{-1}(\frac{n}{b_n}\rho)}\frac{\vp_u(z)\vp_u'(z)}{z^{\a-\e_n}}dz\,,
\Eq(6.lem1.18)
$$
Integrating by parts
$
I''_n(u)\leq
 (1+\d_n)\frac{\vp_u^2(z)}{z^{\a-\e_n}}\biggr|_1^{\vp_u^{-1}(\frac{n}{b_n}\rho)}
+
 (1+\d_n)\int_1^{\vp_u^{-1}(\frac{n}{b_n}\rho)}(\a-\e_n)\frac{\vp_u^2(z)}{z^{1+\a-\e_n}}dz
$.
Using Lemma \thv(A.3.lemma4), (i), we easily see that for all $u>0$,
$
\frac{\vp_u^2(z)}{x^{\a-\e_n}}\bigr|_1^{\vp_u^{-1}(\frac{n}{b_n}\rho)}\leq(\frac{n}{b_n}\rho)^{2-\frac{\a-\e_n}{a}}
$.
Next, 
$$
\eqalign{
\int_1^{\vp_u^{-1}(\frac{n}{b_n}\rho)}\frac{\vp_u^2(x)}{x^{1+\a-\e_n}}dx
=&
\frac{1}{1-a}\int_{(\vp_u^{-1}(\frac{n}{b_n}\rho))^{-(1-a)}}^1 x^{\frac{\a-2a-\e_n}{1-a}-1}e^{-2ux}dx
\cr
\leq &
\sfrac{\a-\e_n}{|\a-2a-\e_n|}\left[1+\bigl(\vp_u^{-1}\bigl(\sfrac{n}{b_n}\rho\bigr)\bigr)^{2a-\a+\e_n}\right]
\leq
\sfrac{\a-\e_n}{|\a-2a-\e_n|}
\bigl[1+
\bigl(e\sfrac{n}{b_n}\rho\bigr)^{2-\frac{\a-\e_n}{a}}
\bigr]
}
\Eq(6.lem1.21)
$$
which is valid for all $u>0$. Indeed, if $u\leq v:=(\frac{n}{b_n}\rho)^{1-a}$, then the last ineqality follows from
Lemma \thv(A.3.lemma4), (ii); if on the contrary $u> v$, then, by Lemma \thv(A.3.lemma4), (iii),
$\vp_u^{-1}(\frac{n}{b_n}\rho)\leq\vp_v^{-1}(\frac{n}{b_n}\rho)$,
whereas by Lemma \thv(A.3.lemma4), (ii), for all $y\geq \frac{n}{b_n}\rho$,
$
\vp_v^{-1}(y)\leq(e\frac{n}{b_n}\rho)^{1/a}
$.
Collecting our bounds we get, reasoning as in the proof of \eqv(6.lem1.5),
that for all large enough $n$,
$
\frac{b_n}{n}I''_n(u)
\leq c_5\frac{b_n}{n}+c_6\left(\frac{b_n}{n}\right)^{\frac{\a-\e_n}{a}-1}
\leq c_5\frac{b_n}{n}+c_6\left(\frac{b_n}{n}\right)^{c_1}
$
for some finite constants, $c_5,c_6,c_1>0$, and all $u>0$. Inserting this in \eqv(6.lem1.11).

Inserting our bounds on  $I'_n(u)$ and $\frac{b_n}{n}I''_n(u)$ in
 \eqv(6.lem1.11)  yields \eqv(6.lem1.6).
Finally, inserting  \eqv(6.lem1.5) and \eqv(6.lem1.6) in \eqv(6.lem1.1')
gives \eqv(6.lem1.1).

It now remains to prove \eqv(6.lem1.2). Write
$$
m_n(u)=\sum_{x\in\VV_n}\E Z_{n,u}(x)\1_{\{Z_n(x)<\rho\}}=n\E Z_{n,u}(x)-n\E Z_{n,u}(x)\1_{\{Z_{n,u}(x)\geq\rho\}}
\Eq(6.lem1.23)
$$
Integration by parts yields
$$
n\E Z_{n,u}(x)
=[b_n\P(\t(x)>r_nz)]\int_0^{\infty}\vp_u'(z)h_n(z)dz
:=[b_n\P(\t(x)>r_nz)]J_n(u)\,,
\Eq(6.lem1.24)
$$
where $h_n$ is given by \eqv(6.lem1.12). As in \eqv(6.lem1.13),  write
$J_n(u)=J_n'(u)+J_n''(u)$, where
$$
J_n'(u)=\int_0^{1}\vp_u'(x)h_n(x)dx
\,,\quad
J_n''(u)=\int_1^{\infty}\vp_u'(x)h_n(x)dx\,,
\Eq(6.lem1.25)
$$
We will treat $J_n$ much in the same way as we treated $I_n$.
On the one hand, proceeding as we did to establish \eqv(6.lem1.15), we obtain that
$
\lim_{n\rightarrow\infty}J'_n(u)= \int_0^{1}\frac{\vp_u'(x)}{x^\a}dx
$
for all $u>0$. On the other hand, using \eqv(6.lem1.17),
$$
(1-\d_n)\int_1^{\infty}\frac{\vp_u'(x)}{x^{\a+\e_n}}dx
\leq
J''_n(u)
\leq
(1+\d_n)\int_1^{\infty}\frac{\vp_u'(x)}{x^{\a-\e_n}}dx\,,
\Eq(6.lem1.27)
$$
where  $0<\e_n, \d_n\downarrow 0$
as $n\uparrow\infty$. Since
$
\int_0^{\infty}\frac{\vp_u'(x)}{x^{\a}}dx
=u^{-\frac{\a-a}{1-a}}\sfrac{\a}{1-a}\G\bigl(\sfrac{\a-a}{1-a}\bigr)
=\nu^{int,-}(u,\infty)
$,
which is finite for $u>0$,
dominated convergence applies and yields,
$
\lim_{n\rightarrow\infty}J''_n(u)=\int_1^{\infty}\frac{\vp_u'(x)}{x^\a}dx
$.
Putting together our results we get that,
$$
\lim_{n\rightarrow\infty}n\E Z_{n,u}(x)
=\int_0^{\infty}\frac{\vp_u'(x)}{x^\a}dx
=\nu^{int,-}(u,\infty)\,,\quad u>0\,.
\Eq(6.lem1.29)
$$

We now want to show that
$
n\E Z_{n,u}(x)\1_{\{Z_{n,u}(x)\geq\rho\}}\rightarrow 0
$
as $n\rightarrow\infty$.
Integration by parts and the change of variable  $y=\vp_u(z)$ yields
$$
\eqalign{
n\E Z_{n,u}(x)\1_{\{Z_{n,u}(x)\geq\rho\}}
&=
b_n\int_{\frac{n}{b_n}\rho}^{\infty}\P(\vp_u(\g_n(x))>y)dy
+
b_n\left(\frac{n}{b_n}\right)\rho\P(\vp_u(\g_n(x))>\rho)
\cr
&=
[b_n\P(\t(x)>r_n)]J''_n\bigl(\vp_u^{-1}(\rho\sfrac{n}{b_n})\bigr)
+\rho n\P\left(Z_{n,u}(x)\geq\rho\right)
\,,
}
\Eq(6.lem1.30)
$$
where $J''_n$ is defined in \eqv(6.lem1.25).
To deal with $J''_n(\vp_u^{-1}(\rho\sfrac{n}{b_n}))$ we use the upper bound \eqv(6.lem1.27)
(which is valid for all $n$ large enough) and, proceeding as in the paragraph below \eqv(6.lem1.18)
(but replacing $\vp_u^2$ by $\vp_u$), we readily obtain that
$
J''_n(\vp_u^{-1}(\rho\sfrac{n}{b_n}))
\leq (1+o(1))(\frac{n}{b_n}\rho)^{1-\frac{\a-\e_n}{a}}
$.
Since we already established (see \eqv(6.lem1.5)) that
$
n\P\left(Z_{n,u}(x)\geq\rho\right)
\leq 2\rho^{-\frac{\a}{a}}(\frac{n}{b_n})^{1-\frac{\a-\e_n}{a}}
\rightarrow 0
$,
the claim that
$n\E Z_{n,u}(x)\1_{\{Z_{n,u}(x)\geq\rho\}}\rightarrow 0$
is established. Inserting this and \eqv(6.lem1.29)
in \eqv(6.lem1.23)
proves \eqv(6.lem1.2).

The proof of Lemma \thv(6.lemma1) for $a>0$ is now complete.
The
case $a=0$ is treated in the same way (see also the remark below Lemma \thv(6.lemma1))
with the difference that
the function $\vp(y)$ and its inverse now become
$\vp(y)=e^{-u/y}$ and  $\vp^{-1}(y)=-\frac{u}{\log y}$, $y\geq 0$.
We omit the details of this elementary adaptation.
\endproof

\proofof{of Proposition \thv(6.prop1)} With the notations of Lemma \thv(6.lemma1)
we may rewrite \eqv(2.4.theo1.1) as
$$
\nu_n(u,\infty)=\frac{a_nr_n^{a}}{b_n}
\frac{
\sum_{x\in\VV_n}Z_n(x)
}
{
\frac{1}{n}\sum_{x\in\VV_n}\t^a(x)
}\,.
\Eq(6.prop1.32)
$$
By the strong law of large numbers, $\frac{1}{n}\sum_{x\in\VV_n}\t^a(x)\rightarrow \E\t^a(x)=1$ $\P$-almost surely,
and by assumption on $a_n$, ${a_nr_n^{a}}/{b_n}\rightarrow 1$.
It thus follows from Lemma \thv(6.lemma1) that, setting
$
\s^2(u)=c_0+c_4 u^{-1+c_2}
$
and choosing
$
z=(b_n/n)^{1/3}\s(u)
$
in \eqv(6.lem1.1), for each fixed $u>0$,
$$
\P\left(\left|\nu_n(u,\infty)-m_n(u)\right|\geq (b_n/n)^{1/3}\s(u)\right)
\leq (b_n/n)^{1/3}+c_3(b_n/n)^{c_1}\,.
\Eq(6.lem1.33)
$$
We now want to make use of Lemma \thv(app.A.4) with
$
X_n(u)=\nu_n(u,\infty)
$,
$
f_n(u)=m_n(u)
$,
$
g_n(u)=\s(u)
$,
$
\eta_n= (b_n/n)^{1/3}$,
and
$
\rho_n= (b_n/n)^{1/3}+c_3(b_n/n)^{c_1}
$.
Indeed $\s(u)$ is a positive decreasing function, so is $m_n$ for each $n$, and the properties \eqv(app.A.4.2) are readily checked. Thus,
$$
\lim_{n\rightarrow 0}\P\left(
\sup_{u>0}\left\{\left|\nu_n(u,\infty)-m_n(u)\right|\geq(b_n/n)^{1/3}\s(u)\right\}\right)
=0\,.
\Eq(6.prop1.34)
$$
Proposition \thv(6.prop1) is proven.
\endproof


We are now ready to prove Proposition \thv(4.prop2).

\proofof{Proposition \thv(4.prop2), (i)}
Let $r_n$ be an intermediate space scale. Assume  that $a<\a$.
Choose $a_n\sim r_n^{-a}b_n$ and $\nu=\nu^{int,-}$ in assertion (i) of Theorem \thv(2.4.theo1).
By Proposition  \thv(6.prop1),
Condition \eqv(2.4.theo1.2) is satisfied in $\P$-probability. To see that Condition \eqv(2.4.theo1.2bis)
also is satisfied
we again make use of Lemma \thv(app.A.4), choosing this time (with the notation of Proposition \thv(6.prop1))
$
X_n(\d)=\int_0^{\d}\nu_n(u,\infty)du
$,
$
f_n(\d)=\int_0^{\d}m_n(u)du
$,
$
g_n(\d)=\int_0^{\d}\s^2(u)du=c_0\d+(c_4/c_2) \d^{c_2}
$
where $0<c_2\leq 1$, and
$
\eta_n= (b_n/n)^{1/3}
$.
Clearly, $f_n(\d)$ and $g_n(\d)$ are positive increasing functions for each $n$,
and the leftmost relation in \eqv(app.A.4.2) is satisfied, albeit with reversed inequality,
for all $l\geq 1/\d_0$ and small enough $\d_0\leq 1$.
Moreover, it follows from  Proposition \thv(6.prop1) that, setting
$$
\textstyle
A_n(\d)=
\left\{\left|\int_0^{\d}\nu_n(u,\infty)du-\int_0^{\d}m_n(u)du\right|
\geq (b_n/n)^{1/3}(c_0\d+(c_4/c_2) \d^{c_2})\right\}\,,
\Eq(6.prop1.4')
$$
there exists a sequence $0<\rho_n\downarrow 0$ such that,
for all $n$ large enough, for all $\d\leq \d_0$,
$
\P\left(A_n(\d)\right)\leq\rho_n
$.
Therefore Lemma \thv(app.A.4) applies, yielding
$
\lim_{n\rightarrow 0}\P\left(\sup_{0<\d<\d_0}A_n(\d)\right)=0
$.
Now by \eqv(6.lem1.2), for all $\e>0$ and all large enough $n$,
$
\bigl|\int_0^{\d}m_n(u)du-\int_0^{\d}\nu^{int,-}(u,\infty)du\bigr|\leq \e\d
$,
while by \eqv(4.prop2.2)
$
\int_0^{\d}\nu^{int,-}(u,\infty)du=c_5\d^{c_6}
$
for some constant $0<c_5,c_6<\infty$  (that depend only $\a$ and $a$).
Hence we have established that there exists
$\O^{\tau}_{2,n}\subset\O^{\tau}$ with $\P(\O^{\tau}_{2,n})\geq 1-o(1)$
such that for all $n$ large enough, on $\O^{\tau}_{2,n}$, 
for all $0<\d\leq\d_0$,
$
\int_0^{\d}\nu_n(u,\infty)du\leq c_7\d^{c_8}
$,
where $0<c_7,c_8<\infty$ are constants (that depend only on $\a$ and $a$).
Condition \eqv(2.4.theo1.2bis) is thus satisfied in $\P$-probability. So,
all conditions of assertion (i) of Theorem \thv(2.4.theo1) are satisfied in $\P$-probability.
The proof of assertion (i) of Proposition \thv(4.prop2) is done.\endproof

\bigskip


\chap{6. Extreme scales.}6

Let us motivate the strategy we will implement in this section.
Consider the re-scaled sequence $\g_n(x)=r_n^{-1}\t(x)$, $x\in\VV_n$.
For each $n$ form the point process
$
\Upsilon_n=\sum_{x\in\VV_n}\1_{\g_n(x)}
$,
and let
$\Upsilon=\sum_{k=1}^{\infty}\1_{\g_k}$ be $\PRM(\mu)$ with $\mu$  given by \eqv(4.stat.1).
It is well known that when
$(\t(x), x\in\VV_n)$ are i.i.d\. r.v\.'s  equi-distributed with $\t\in\DD(\a)$,
$\Upsilon_n$ converges weakly to $\Upsilon$ if and only if $r_n$ is an extreme space scale.
\note{
By e.g\. [Re] Proposition 3.21 p.154, $\Upsilon_n$ converges weakly
to $\Upsilon$ if and only if the sequence $r_n$ satisfies
$
\lim_{n\rightarrow\infty} n\P(\tau(x)>r_n u)= u^{-\a}
$,
$u>0$; in view of Definition \thv(4.def1),
this is the same as saying that $r_n$ is an extreme space scale.
}
Thus, on extreme scales, the convergence of (appropriate almost sure continuous) functionals of
$\Upsilon_n$ simply follows from the weak convergence  of $\Upsilon_n$, using the Continuous Mapping Theorem.
However, this only yields convergence in distribution, which is not enough for our needs.

The usual way out of this difficulty is to think of weak convergence from Skorohod's
representation Theorem and replace the sequence $(\g_n(x), x\in\VV_n)$
by a new sequence with  identical distribution, but almost sure convergence properties.
This strategy was first implemented in the context of an aging system by Fontes et {\it al.} \cite{FIN},
and fruitfully
used in many subsequent papers and various models (see the review paper \cite{BC2}). Technically, all these works rely
on one specific choice of extreme scale.
In the present paper we have no such restrictions.

In Subsection 7.1, we give an explicit representation
of the re-scaled landscape (Lemma \thv(7.lemma1)) which is valid for all extreme scales,
and establish its convergence properties
(Proposition \thv(7.prop1)). In Subsection 7.2 we consider the model obtained
by substituting the representation for the original landscape and prove Proposition \thv(4.prop3).
(The second assertion of Propositions \thv(4.prop1) and \thv(4.prop2) will be proved there as well).
The final Subsection 7.3 contains the proofs  of  Lemma \thv(4.lemma2) and Lemma \thv(4.lemma3).

\bigskip
\line{\bf 6.1. A representation of the re-scaled landscape.\hfill}

The representation we now introduce
is due to Lepage {\it et al\.} \cite{LWZ} and relies on an elementary property of
order statistics.
Let $\bar\t_{n}(1)\geq\dots\geq\bar\t_{n}(n)$ and $\bar\g_{n}(1)\geq\dots\geq\bar\g_{n}(n)$ denote, respectively,
the landscape and re-scaled landscape variables
 $(\t(x), x\in\VV_n)$ and $(\g_n(x), x\in\VV_n)$
arranged in decreasing order of magnitude. For $u\geq 0$ set
$G(u)=\P(\tau(x)>u)$ and
$$
G^{-1}(u):=\inf\{y\geq 0 : G(y)\leq u\}\,.
\Eq(7.1.1)
$$
Let $(E_i, i\geq 1)$ be a sequence of i\.i\.d\. mean one exponential random variables
defined on a common probability space $(\O^{E}, \FF^{E}, \bold P)$.
We will now see that both the ordered landscape variables and the limiting point process $\Upsilon$
can be expressed in terms of this sequence. Set, for $k\geq 1$,
$$
\eqalign{
\G_k&=\sum_{i=1}^k E_i\,,\cr
\g_k&=\G_k^{-1/\a}\,,
}
\Eq(7.1.2)
$$
and, for $1\leq k\leq n$, $n\geq 1$,
$$
\g_{nk}=r_n^{-1}G^{-1}(\G_k/\G_{n+1})\,.
\Eq(7.1.3)
$$

\lemma{\TH(7.lemma1)} {\it  For each $n\geq 1$,
$
(\bar\g_{n}(1),\dots,\bar\g_{n}(n))\overset d\to=(\g_{n1},\dots,\g_{nn})\,.
$
}

\proof Note that $G$ is non-increasing and right-continuous so that $G^{-1}$
is non-increasing and right-continuous. It is well known that if the random variable $U$
is a uniformly distributed on $[0,1]$ we may write $\t(0)\overset d\to=G^{-1}(U)$
(see e\.g\. \cite{Re}, page 4).
In turn it is well known (see \cite{Fe}, Section III.3) that if $(U(k), 1\leq k\leq n)$
are independent random variables uniformly distributed on $[0,1]$ then, denoting by
$\bar U_{n}(1)\leq\dots\leq \bar U_{n}(n)$ their ordered statistics,
$
(\bar U_{n}(1),\dots,\bar U_{n}(n))\overset d\to=(\G_1/\G_{n+1},\dots,\G_n/\G_{n+1})
$.
Combining these two facts yields the claim of the lemma.
\endproof

Next,
let $\Upsilon$ be the point process
in $M_P(\R_+)$
which has counting function
$$
\Upsilon([a,b])=\sum_{i=1}^{\infty}\1_{\{\g_k\in [a,b]\}}\,.
\Eq(7.1.5)
$$
\lemma{\TH(7.lemma4)} {\it $\Upsilon$ is a Poisson random measure on $(0,\infty)$
with mean measure $\mu$ given by \eqv(4.stat.1).
}

\proof
The point process 
$
\G=\sum_{i=1}^{\infty}\1_{\{\G_k\}}
$
defines a homogeneous Poisson random measure on $[0,\infty)$ and thus,
by the mapping theorem (\cite{Re}, Proposition 3.7), setting $T(x)=x^{-1/\a}$ for $x>0$,
$\Upsilon=\sum_{i=1}^{\infty}\1_{\{T(\G_k)\}}$ is Poisson random measure on $(0,\infty)$
with mean measure $\mu(x,\infty)=T^{-1}(x)$.
\endproof

Then, on the fixed probability space $(\O^{E}, \FF^{E}, \bold P)$,
all random variables of interest
will have an almost sure limit.

\proposition{\TH(7.prop1)} {\it Let $r_n$ be an extreme space scale. Let $f:(0,\infty)\rightarrow[0,\infty)$ be
a continuous function that obeys
$$
\int_{(0,\infty)}\min(f(u), 1)d\mu(u)<\infty\,.
\Eq(7.1.6)
$$
Then, $\bold P$-almost surely,
$$
\lim_{n\rightarrow\infty}\sum_{k=1}^{n}f(\g_{nk})=\sum_{k=1}^{\infty}f(\g_{k})<\infty\,.
\Eq(7.1.7)
$$
}

The proof is inspired from the proof of Proposition 3.1 of \cite{FIN}. It relies on the following two lemmata. Set
$$
g_n(x)=r_n^{-1}G^{-1}(x/n)\,.
\Eq(7.1.8)
$$
\lemma{\TH(7.lemma5)} {\it For any fixed $x<\infty$, $g_n(x)\rightarrow x^{-1/\a}$ as $n\rightarrow\infty$.
}

\lemma{\TH(7.lemma6)} {\it For any $\d>0$ there exit constants $0<C',C''<\infty$ such that, for $n\geq C''$,
$$
g_n(x)\leq  C'x^{-(1-\d)/\a}\,,\quad
C''\leq x\leq n\,.
\Eq(7.1.9)
$$
}

In the sequel we use the notation and results of Appendix A.3 on regular variations.

\proofof{Lemma \thv(7.lemma5)}  Observe first that by assumption
$G\in R_{-\a}$. Thus, by Lemma \thv(A.3.lemma3),  $G^{-1}\in R_{-1/\a}(0+)$.
Observe next that if $r_n$ is an extreme space scale, taking $b_n=n$ in \eqv(4.9)
yields $nG(r_n)\sim 1$, and invoking again Lemma \thv(A.3.lemma3),
$r_n^{-1}G^{-1}(1/n)\sim 1$. Using these two observations we may write
$$
g_n(x)=\frac{r_n^{-1}G^{-1}(x/n)}{r_n^{-1}G^{-1}(1/n)}[r_n^{-1}G^{-1}(1/n)]
=x^{-1/\a}[1+o(1)]\frac{\ell(x/n)}{\ell(1/n)}\,,
\Eq(7.1.10)
$$
for some function $\ell$ which is slowly varying at $0+$.
But this implies that for any $x<\infty$, $g_n(x)\rightarrow x^{-1/\a}$ as $n\rightarrow \infty$.
\endproof

\proofof{Lemma \thv(7.lemma6)} Let $\l\in(0,1)$ be a constant whose value will be chosen later, and
assume that $1/\l\leq  x\leq n$. By
\eqv(7.1.10),
$$
g_n(\l x)=\l^{-1/\a} x^{-1/\a}[1+o(1)]\frac{\ell(\l x/ n)}{\ell(1/n)}\,.
\Eq(7.1.11)
$$
By the Representation Theorem (Theorem \thv(A.3.theo1) of Appendix A.3) adapted to the case of functions
that are slowly varying at zero, the quotient in the right hand side of \eqv(7.1.11) may be written as
$$
\frac{\kappa(n/(\l x))}{\kappa(n)}
\exp\left\{\int_{n/(\l x)}^{n}\frac{\varepsilon(y)}{y}dy\right\}\quad(\l x\geq 1)\,,
\Eq(7.1.12)
$$
where $\kappa(y)\rightarrow \kappa\in(0,\infty)$, and
$\varepsilon(y)\rightarrow 0$ as $y\rightarrow\infty$.
Now
$$
\left|\int_{n/(\l x)}^{n}\frac{\varepsilon(y)}{y}dy\right|
\leq \d'\left|\int_{n/(\l x)}^{n}\frac{1}{y}dy\right|
\leq \d'|\log(1/(\l x))|\,,
\Eq(7.1.13)
$$
where, since $x\leq n$, $\d'=\d'(\l)=\sup\{|\varepsilon(y)|, y\geq 1/\l\}$.
Thus, since $\l x\geq 1$, the exponential in \eqv(7.1.12) is bounded above by
$$
(\l x)^{\d'}\,.
\Eq(7.1.14)
$$
Given $\d>0$ we may now choose $\l\in(0,1)$ in such a way that $\d'(\l)<\d/\a$ and that
$\kappa(y)\in[\kappa/2,\kappa]$ for $y\geq 1/\l$. The lemma now follows from \eqv(7.1.11)-\eqv(7.1.14)
with $C'=4\l^{-(1+\d')}$ and $C''=1/\l$.
\endproof

\proofof{Proposition \thv(7.prop1)} By the Strong law of large numbers there exists a subset $\wt\O^{E}\subset\O^{E}$
of full measure such that, for all  $n$ large enough and all  $\o\in\wt\O^{E}$,
$
\G_n=n(1+\l_n)
$
where
$
\l_n=o(1)
$.
From now on we assume that $\o\in\wt\O^{E}$. Thus
$$
\sum_{i=1}^{n}f(\g_{ni})=\sum_{i=1}^{n}f\left(r_n^{-1}G^{-1}(\G_i/[n(1+\l_n)])\right)\,.
\Eq(7.1.15)
$$
Let us first consider the case $f(x)=x$, $x>0$.
Recall the notation $\g_i=\G^{-1/\a}_i$. For $y>0$ set
$
I(y)=\{i\geq 1 : \g_i\geq y\}
$,
$
I^c(y)=\{i\geq 1 : \g_i< y\}
$
and, for $\d>0$ and large enough $n$ write:
$$
\sum_{i=1}^{n}\g_{ni}=\sum_{i\in I(\d)}\g_{ni}+\sum_{i\in I(n^{-1/\a})\setminus I(\d)}\g_{ni}+\sum_{i\in I^c(n^{-1/\a})}\g_{ni}\,.
\Eq(7.1.16)
$$
From Lemma \thv(7.lemma5) and \eqv(7.1.11) it follows that,
$$
\sum_{i\in I(\d)}\g_{ni}\rightarrow \sum_{i\in I(\d)}\g_i\,,\quad n\rightarrow\infty\,.
\Eq(7.1.17)
$$

Next, by Lemma \thv(7.lemma6), for all $0<\d<1$ and some constant $0<C<\infty$, we have
$$
\sum_{i\in I(n^{-1/\a})\setminus I(\d)}\g_{ni}
\leq  \sum_{i\in I(n^{-1/\a})\setminus I(\d)} C\G_i^{-(1-\d)/\a}
=\sum_{i\in I(n^{-1/\a})\setminus I(\d)} C\g_i^{(1-\d)}\,.
\Eq(7.1.18)
$$
The last sum in \eqv(7.1.18) is bounded above by
$$
W_\d= \sum_{i:\g_i\leq \d} C\g_i^{(1-\d)}\,.
\Eq(7.1.19)
$$
Just as in \cite{FIN} page 601,
we now claim that, with $\d>0$ chosen such that $\d+\a<1$,
$
W:=\lim_{\d\rightarrow 0}W_\d=0
$
$\bold P$-almost surely. To prove this note that $W$ is well defined by monotonicity, and is non-negative.
By standard Poisson calculation,
$$
\bold E(W_\d)=\a\int_{0}^{\d}w^{1-\d}w^{-(1+\a)}dw\leq \frac{\a}{1-(\d+\a)}\d^{1-(\d+\a)}\,,
\Eq(7.1.20)
$$
so that $\bold E(W_\d)\rightarrow 0$ as $\d\rightarrow 0$. By dominated convergence, $\bold E(W)=0$, and
the claim follows.

Finally, for $i\in I^c(n^{-1/\a})$, we have $\G_i/n\leq 1$. Since $G$ is right-continuous non-increasing and since,
being the tail of a probability distribution, $G(x)\rightarrow 0$ as $x\rightarrow 1$, we have, for large enough $n$,
$$
G^{-1}(\G_i/[n(1+\l_n)])\leq G^{-1}(1/(1+\l_n))\leq 1\,.
\Eq(7.1.21)
$$
Thus
$$
\sum_{i\in I^c(n^{-1/\a})}\g_{ni}
\leq \sum_{i\in I^c(n^{-1/\a})} r_n^{-1}
\leq n r_n^{-1}\,.
\Eq(7.1.22)
$$
Now remember from the proof of Lemma \thv(7.lemma5) that $G^{-1}\in R_{-1/\a}(0+)$
and $r_n^{-1}G^{-1}(1/n)\sim 1$. Therefore
$
r_n^{-1}n^{1/\a}\ell(1/n)\sim 1
$
where the function $\ell$ is slowly varying at $0+$, and together with \eqv(7.1.22) this yields,
$$
\sum_{i\in I^c(n^{-1/\a})}\g_{ni}\rightarrow 0\,,\quad n\rightarrow\infty\,.
\Eq(7.1.24)
$$

Combining the previous estimates we obtain that, on a subset of $\O^{E}$ of full measure,
$$
\lim_{n\rightarrow\infty}\sum_{i=1}^{n}\g_{ni}
=\lim_{\d\rightarrow 0}\sum_{i:\g_i\geq \d}\g_i
=\sum_{i=1}^{\infty}\g_{i}\,.
\Eq(7.1.25)
$$

Proving \eqv(7.1.7) is now simple. We only indicate the main modifications.
By assumption \eqv(7.1.6)
the sum in the right hand side of \eqv(7.1.7) is almost surely finite. It moreover follows
from \eqv(7.1.6) that there exists $\b>\a$ such that $f(x)\leq x^{\b}$ for all $x$
sufficiently small. Thus, to deal with the second term in the right hand side of \eqv(7.1.16),
write
$$
\sum_{i\in I(n^{-1/\a})\setminus I(\d)}f(\g_{ni})
\leq\sum_{i\in I(n^{-1/\a})\setminus I(\d)} C\g_i^{\beta(1-\d)}
\Eq(7.1.26)
$$
instead of in \eqv(7.1.18), choose $\d>0$ small enough so that $\beta(1-\d)-\a>0$, and proceed
as in \eqv(7.1.19)-\eqv(7.1.20). Similarly, to bound the third term in the right hand side of
\eqv(7.1.16), write
$$
\sum_{i\in I^c(n^{-1/\a})}f(\g_{ni})
\leq n r_n^{-\b}
\Eq(7.1.27)
$$
instead of \eqv(7.1.22) and proceed as in \eqv(7.1.21)-\eqv(7.1.24).
Turning to the first term in the right hand side of
\eqv(7.1.16), we obviously have, proceeding as in \eqv(7.1.17),
$$
\sum_{i\in I(\d)}f(\g_{ni})\rightarrow \sum_{i\in I(\d)}f(\g_i)\,,\quad n\rightarrow\infty\,.
\Eq(7.1.28)
$$
We may then conclude just as in \eqv(7.1.24). The proposition is thus proven.
\endproof


We conclude this section with the proof of Lemma \thv(4.lemma4)

\proofof{Lemma \thv(4.lemma4)} We have to establish that
$
\underline r_n\ll \bar r_n\ll r_n
$
where $\underline r_n$, $\bar r_n$ and $r_n$ denote, respectively,
a constant, an intermediate and an extreme space scale.
Let us first prove that $\bar r_n/r_n=o(1)$.
Using that $G^{-1}\in R_{-1/\a}(0+)$
(see the proof of Lemma \thv(7.lemma5)), it follows from Definition \thv(4.def1) that
$$
\bar r_n/r_n
=
\frac{G^{-1}\bigl(\bar b_n^{-1}(1+o(1))\bigr)}{G^{-1}\bigl(b_n^{-1}(1+o(1))\bigr)}
=(1+o(1))\left(\frac{b_n}{\bar b_n}\right)^{1/\a}
\frac{\ell\bigl(\bar b_n^{-1}(1+o(1))\bigr)}{\ell\bigl(b_n^{-1}(1+o(1))\bigr)}\,,
\Eq(4.lem4.2)
$$
for some function $\ell$ which is slowly varying at $0+$.
From our assumption on $\bar r_n$ and $r_n$, and Definition \thv(4.def1) it is plain that,
with obvious notations,
$$
1\ll \bar b_n\ll b_n\leq n\,.
\Eq(4.lem4.1)
$$
Hence, by \eqv(4.lem4.2) and Lemma \thv(A.3.lemma2),
$
0\leq \bar r_n/r_n\leq (1+\d_n)\left(b_n/\bar b_n\right)^{1/\a-\e_n}
$,
where $\e_n\downarrow 0$, $\d_n\downarrow 0$ as $n\uparrow\infty$.
From this and \eqv(4.lem4.1) the claim follows.
It remains to prove  that $\underline r_n/\bar r_n=o(1)$.
Since by definition $\underline r_n$ is a constant,
it suffices to show that
$
G^{-1}\bigl(\bar b_n^{-1}(1+o(1))\bigr)\uparrow\infty
$
as $n\uparrow\infty$. Now this is plain since $G^{-1}\in R_{-1/\a}(0+)$,
and since, by assumption, $\bar b_n\uparrow\infty$ as $n\uparrow\infty$.
Lemma \thv(4.lemma4) is proven.\endproof


\bigskip
\line{\bf 6.2. Proof of Proposition \thv(4.prop3).\hfill}

In this subsection we consider the model obtained by substituting the new landscape
$(\g_{nk}, 1\leq k\leq n)$ for the original (re-scaled) landscape $(\g_n(x), x\in\VV_n)$.
We assume throughout that $r_n$ is an extreme space scale.
As for short and intermediate space scales, the proof of Proposition \thv(4.prop3)
relies on Theorem \thv(2.4.theo1).
To distinguish the quantity  $\nu_n(u,\infty)$, expressed in \eqv(2.4.theo1.1) in the original landscape variable, from its
expression in the new landscape variables, we call the latter $\bold v_n(u,\infty)$. Therefore
$$
\bold v_n(u,\infty)=a_n
\frac{
\sum_{k=1}^n(r_n\g_{nk})^ae^{-u/\g_{nk}^{(1-a)}}
}
{
\sum_{k=1}^n(r_n\g_{nk})^a
}\,,\quad u\geq 0\,.
\Eq(4.prop3.5)
$$
We first treat the numerator in \eqv(4.prop3.5). For $u\geq 0$ set
$$
\vp_u(y)=y^ae^{-u/y^{(1-a)}}\,,\quad y\geq 0\,.
\Eq(4.prop3.6)
$$
We want to apply Proposition \thv(7.prop1) to the sum
$
\sum_{k=1}^n\vp_u(\g_{nk})
$.
For this let $x^*$ be defined through $\vp_u(x^*)=1$. Noting that $0<x^*\leq 1$ for $0\leq a< 1$ and $u\geq 0$,
a simple calculation yields
$
\int_{(0,\infty)}\min(\vp_u(y), 1)d\mu(y)
=\frac{\a}{1-a}\int_{1/x^*}^{\infty}y^{-\frac{1-\a}{1-a}}e^{-uy}dy
+(x^*)^{-\a}
$,
which is always finite if $u>0$, regardless of the respective size of $a$ and $\a$. Thus, for all $u>0$, $\bold P$-almost surely,
$$
\lim_{n\rightarrow\infty}\sum_{k=1}^n\vp_u(\g_{nk})=\sum_{k=1}^{\infty}\vp_u(\g_{k})
<\infty\,.
\Eq(4.prop3.7)
$$
In contrast, the behavior of the denominator in \eqv(4.prop3.5) will depend on whether $a$ is larger or smaller than $\a$.

\noindent{\bf The case $a>\a$}. Here we have
$
\int_{(0,\infty)}\min(x^a, 1)d\mu(x)<\infty
$,
so that $\bold P$-almost surely,
$$
\lim_{n\rightarrow\infty}\sum_{k=1}^n\g_{nk}^a=\sum_{k=1}^{\infty}\g_{k}^a
<\infty\,.
\Eq(4.prop3.8)
$$
In that case, choosing $a_n=1$ in \eqv(4.prop3.5),
we get, collecting \eqv(4.prop3.7) and \eqv(4.prop3.8), that for all $u>0$, $\bold P$-almost surely,
$$
\lim_{n\rightarrow\infty}\bold v_n(u,\infty)
=
\lim_{n\rightarrow\infty}
\frac{
\sum_{k=1}^n\vp_u(\g_{nk})
}
{
\sum_{k=1}^n\g_{nk}^a
}
=
\frac{\sum_{k=1}^{\infty}\vp_u(\g_{nk})}{\sum_{k=1}^{\infty}\g_{nk}^a}
=\nu^{ext,+} (u,\infty)\,.
\Eq(4.prop3.9)
$$
It is plain that $\nu^{ext,+}$ is a  probability measure with continuous density:
indeed it is an infinite mixture of exponential densities,
the coefficients of the mixture being the weights ${\g_{k}^a}/{\sum_{l}\g_{l}^a}$ of
Poisson-Dirichlet random probability measure with parameter $\a/a$.
From the monotonicity of $\bold v_n(u,\infty)$ and the continuity of the limiting function $\nu^{ext,+}(u,\infty)$
we conclude that there exists a subset $\O^{E}_1\subset\O^{E}$ of
the sample space $\O^{E}$ of the $\g$'s with the property that $\bold P(\O^{E}_1)=1$, and such that, on $\O^{E}_1$,
$$
\lim_{n\rightarrow\infty}\bold v_n(u,\infty)=\nu^{ext,+} (u,\infty)
\,,\quad\forall\,u\geq 0\,.
\Eq(4.prop3.10)
$$
Condition \thv(2.4.theo1.2) of assertion (i) of Theorem \thv(2.4.theo1) is thus satisfied $\bold P$-almost surely.
To see that Condition \thv(2.4.theo1.2bis) also is satisfied on a set of full measure
we use that on $\O^{E}_1$, by \eqv(4.prop3.10), for all $0<\d\leq\d_0$ and some $0<\d_0\leq 1$,
$
\lim_{n\rightarrow\infty}\int_0^{\d}\bold v_n(u,\infty)du=\int_0^{\d}\nu^{ext,+}(u,\infty)du
$.
Again the monotonicity of $\int_0^{\d}\bold v_n(u,\infty)du$ and the continuity of the limiting function
allow us to conclude that there exists of a subset $\O^{E}_2\subset\O^{E}$ with the property that
$\bold P(\O^{E}_2)=1$, and such that, on $\O^{E}_2$,
$
\lim_{n\rightarrow\infty}\int_0^{\d}\bold v_n(u,\infty)du=\int_0^{\d}\nu^{ext,+}(u,\infty)du
$
for all $0<\d\leq\d_0$.
We may thus pass to the limit $\d\rightarrow 0$ and write
$
\lim_{\d\rightarrow 0}\lim_{n\rightarrow\infty}\int_0^{\d}\bold v_n(u,\infty)du=\lim_{\d\rightarrow 0}\int_0^{\d}\nu^{ext,+}(u,\infty)du
$.
Now by \eqv(4.prop3.2),
$$
\int_0^{\d}\nu^{ext,+}(u,\infty)du
=
\sum_{k}\frac{\g_k}{\sum_{l}\g_l^a}\Bigl(1-e^{-\d\g_k^{-(1-a)}}\Bigr)
\leq
\d e^{\sqrt\d}\nu^{ext,+}(\d,\infty)
+
\sum_{k}\frac{\g_k}{\sum_{l}\g_l^a}\1_{\{\g_k\leq \d^{{1}/{2(1-a)}}\}}\,,
\Eq(4.prop3.10')
$$
where we proceeded as in \eqv(2.2.A3'.1)-\eqv(2.2.A3'.4) to derive the upper bound.
Now from this bound, Lemma \thv(4.lemma3), and \eqv(7.1.19)-\eqv(7.1.20), it follows that
$\lim_{\d\rightarrow 0}\int_0^{\d}\nu^{ext,+}(u,\infty)du=0$ $\bold P$-almost surely.
All the assumptions of assertion (i) of Theorem \thv(2.4.theo1) are thus satisfied $\bold P$-almost surely.
The proof of Proposition \thv(4.prop3) in the case  $a>\a$ is complete.

\noindent{\bf The case $a<\a$}. Since $\E(r_n\g_{nk})^a<\infty$ then, clearly, the sum
$\sum_{k=1}^n\g_{nk}^a$  is wrongly normalized. Here we rewrite
\eqv(4.prop3.5) in the form
$$
\bold v_n(u,\infty)=\frac{a_nr_n^{a}}{n}
\frac{
\sum_{k=1}^n\vp_u(\g_{nk})
}
{
\frac{1}{n}\sum_{k=1}^n(r_n\g_{nk})^a
}\,,\quad u\geq 0\,.
\Eq(4.prop3.11)
$$
As in the case of short and intermediate space scales we want to control the denominator via a strong law of large numbers.
One easily checks that since the law of the variables $r_n\g_{nk}$ is independent of $n$
(namely, for each $n$ and $k$, $\bold P(r_n\g_{nk}>u)=\P(\t>u)$, $\t\in\DD(\a)$)
the classical proof by Etemadi \cite{E} goes through, yielding
%
%
$$
\lim_{n\rightarrow\infty}\frac{1}{n}\sum_{k=1}^n(r_n\g_{nk})^a=\E\t^a<\infty
\text{$\bold P$-almost surely.}
\Eq(4.prop3.12)
$$
Thus, choosing $a_nn\sim nr_n^{-a}$ in \eqv(4.prop3.12) (or, equivalently, $a_n\sim r_n^{-a}b_n$),
we get that for all $u>0$, $\bold P$-almost surely,
$$
\lim_{n\rightarrow\infty}\bold v_n(u,\infty)
=
\lim_{n\rightarrow\infty}
\frac{
\sum_{k=1}^n\vp_u(\g_{nk})
}
{
\frac{1}{n}\sum_{k=1}^n(r_n\g_{nk})^a
}
=
\sum_{k=1}^{\infty}\frac{\vp_u(\g_{nk})}{\E\t^a}
=\nu^{ext,-}(u,\infty)<\infty\,.
\Eq(4.prop3.13)
$$
Now using Lemma \thv(4.lemma3) one easily checks that
$\int_{(0,\infty)}(1\wedge u)\nu^{ext,-}(du)<\infty$,
and  that $\nu^{ext,-}(u,\infty)$ is continuous, on a subset of full measure.
Again we conclude that
there exists a subset $\O^{E}_2\subset\O^{E}$ of
the sample space $\O^{E}$ of the $\g$'s with the property that $\bold P(\O^{E}_2)=1$, and such that, on $\O^{E}_2$,
$$
\lim_{n\rightarrow\infty}\bold v_n(u,\infty)=\nu^{ext,-}(u,\infty)
\,,\quad\forall\,u\geq 0\,.
\Eq(4.prop3.14)
$$
The conditions of assertion (ii) of Theorem \thv(2.4.theo1) are thus satisfied $\bold P$-almost surely.
Proposition \thv(4.prop3) is proven in the case  $a<\a$.
Of course, taking the intersection $\O^{E}_1\cap\O^{E}_2$, the two convergence results of
\eqv(4.prop3.4) can be stated simultaneous on a common full measure set. \endproof

\bigskip
\line{\bf 6.3. Proof of  Lemma  \thv(4.lemma2) and Lemma  \thv(4.lemma3).\hfill}

The proof of  Lemma  \thv(4.lemma2) is elementary.
We skip it and focus on the more involved proof of Lemma \thv(4.lemma3).
%
Recall from \eqv(4.prop3.6) that, for $u\geq 0$,
$\vp_u(y)=y^ae^{-u/y^{(1-a)}}$, $y\geq 0$, and write $\vp_1\equiv\vp$.
Set $u^{-\frac{\a-a}{1-a}}=m$. By \eqv(4.prop3.2) we may write
$$
u^{\frac{\a-a}{1-a}}\nu^{ext, -}(u,\infty)=\frac{1}{m}\sum_{k=1}^{\infty}\vp(m^{1/\a}\g_k)\,.
\Eq(4.lem3.2)
$$
Assertion (i) of the lemma will thus be proven if we can prove that
$$
\lim_{m\rightarrow\infty}\frac{1}{m}\sum_{k=1}^{\infty}\vp(m^{1/\a}\g_k)=
\sfrac{\a}{1-a}\G\bigl(\sfrac{\a-a}{1-a}\bigr)\text{$\bold P$-almost surely.}
\Eq(4.lem3.3)
$$
Note that  for this it is enough to take the limit along the integers since,
$\vp(m^{1/\a}\g_k)$ being a strictly increasing function of $m$,
$$
\frac{\lfloor m \rfloor}{m}\frac{1}{\lfloor m \rfloor}\sum_{k=1}^{\infty}\vp(\lfloor m \rfloor^{1/\a}\g_k)
\leq
\frac{1}{m}\sum_{k=1}^{\infty}\vp(m^{1/\a}\g_k)
\leq
\frac{\lceil m \rceil}{m}\frac{1}{\lceil m \rceil}\sum_{k=1}^{\infty}\vp(\lceil m \rceil^{1/\a}\g_k)\,.
\Eq(4.lem3.4)
$$
The proof now proceeds as follows.
Given a threshold function $M\equiv M(m)$ (to be chosen later)
let $\PRM(\mu^+_{M})$ and $\PRM(\mu^-_{M})$
be the Poisson point processes with points $\{\g^\pm_k\}$ whose intensity measures are defined through
$$
\mu^-_{M}(A)=\mu(A\cap(0,M/m^{1/\a})) \text{and} \mu^+_{M}(A)=\mu(A\cap[M/m^{1/\a},\infty))
\Eq(4.lem3.5)
$$
for any Borel set $A\subseteq(0,\infty)$.
(In other words $\PRM(\mu^+_{M})$ and $\PRM(\mu^-_{M})$ are $\PRM(\mu)$
restricted to the sets $(0,M/m^{1/\a})$ and $[M/m^{1/\a},\infty)$ respectively).
Using these two processes we break the middle sum in \eqv(4.lem3.4) into
$
\frac{1}{m}\sum_{k=1}^{\infty}\vp(m^{1/\a}\g^-_k)
+
\frac{1}{m}\sum_{k=1}^{\infty}\vp(m^{1/\a}\g^+_k)
$.
We will show that if $M$ is of the form
$
M=\varepsilon\bigl(\frac{m}{\log m}\bigr)^{\frac{1}{\a}}
$,
for some small enough $0<\varepsilon<1$, then, $\bold P$-almost surely,
$$
\lim_{m\rightarrow\infty}\frac{1}{m}\sum_{k=1}^{\infty}\vp(m^{1/\a}\g^-_k)=
\sfrac{\a}{1-a}\G\bigl(\sfrac{\a-a}{1-a}\bigr)\,,
\Eq(4.lem3.6)
$$
$$
\text{and}\lim_{m\rightarrow\infty}\frac{1}{m}\sum_{k=1}^{\infty}\vp(m^{1/\a}\g^+_k)=
0\,.\,\,\,\,\,\,\,\,\,\,\,\,\,\,\,\,\,\,\,\,\,\,\,\,\,\,\,\,\,\,\,\,\,\,\,\,\,\,\,\,\,\,\,\,\,\,\,\,\,\,\,
\Eq(4.lem3.7)
$$

We first prove \eqv(4.lem3.6). The  boundedness of  the Poisson points $\g^-_k$  enables us to
use a classical large deviation upper bound. Set
$$
A_m=\left\{\left|
\frac{1}{m}\sum_{k=1}^{\infty}\vp(m^{1/\a}\g^-_k)-
\bold E \frac{1}{m}\sum_{k=1}^{\infty}\vp(m^{1/\a}\g^-_k)
\right|\geq\d_m\right\}\,,
\Eq(4.lem3.8)
$$
where
$
{\d_m}=2\bigl(\frac{\log m}{m}\bigr)^{1-\frac{a}{\a}}
$.
By Tchebychev exponential inequality,
for all $\l>0$,
$$
\bold P\left(A_m\right)\leq 2 \exp\left\{-\l\d_m-\bold E ({\l}/{m})\sum_{k=1}^{\infty}\vp(m^{1/\a}\g^-_k)
+\log\bold E \exp\left\{({\l}/{m})\sum_{k=1}^{\infty}\vp(m^{1/\a}\g^-_k)\right\}
\right\}\,.
\Eq(4.lem3.10)
$$
Simple Poisson point process calculations yield
$
\bold E \frac{1}{m}\sum_{k=1}^{\infty}\vp(m^{1/\a}\g^-_k)=\s^{(1)}_M
$,
where
$$
\s^{(1)}_M
=\sfrac{\a}{1-a}\int_{1/M^{1-a}}^{\infty}y^{\sfrac{\a-a}{1-a}-1}e^{-y}dy
=(1-o(1))\sfrac{\a}{1-a}\G\bigl(\sfrac{\a-a}{1-a}\bigr)\,,
\Eq(4.lem3.9)
$$
and
$$
\log\bold E \exp\left\{({\l}/{m})\sum_{k=1}^{\infty}\vp(m^{1/\a}\g^-_k)\right\}
=-\int_{0}^{\infty}(1-e^{\frac{\l}{m} \vp(m^{1/\a}x)})d\mu^-_{M}(x)
\,.
\Eq(4.lem3.11)
$$
Furthermore, for all $l>1$,
$
\int_{0}^{\infty}\vp^k(m^{1/\a}x))d\mu^-_{M}(x):=m\s^{(k)}_M
$, where
$$
\s^{(l)}_M=\sfrac{\a}{1-a}\int_{1/M^{1-a}}^{\infty}y^{\sfrac{\a-la}{1-a}-1}e^{-y}dy
\,.
\Eq(4.lem3.12)
$$
In the worst situation $\a<la$ for all $l>1$ (indeed if $\a\geq la$,
then $\s^{(l)}_M\leq\sfrac{\a}{1-a}\G\bigl(\sfrac{\a-la}{1-a}\bigr)<\infty$).
Let us thus assume that $\a<la$ for all $l>1$. In this case,
$
\s^{(l)}_M\leq\bar\s^{(l)}_M:=\sfrac{\a}{(1-a)(2a-\a)}M^{al-\a}
$,
and so,
$$
-\int_{0}^{\infty}(1-e^{\frac{\l}{m} f(m^{1/\a}x)})d\mu(x)
\leq
\s^{(1)}_M\l+
\bar\s^{(2)}_M\frac{\l^2}{4m}e^{\frac{\l M^a}{2m}}\,.
\Eq(4.lem3.14)
$$
Inserting this bound in \eqv(4.lem3.11), plugging the result in \eqv(4.lem3.10), and choosing $\l=\d_m{2m}/{\bar\s^{(2)}_M}$,
we obtain
$$
\bold P\left(A_m\right)\leq
2\exp\left\{-{\d^2_mm}/{\bar\s^{(2)}_M}\left(2-e^{{2\d_mM^a}/{\bar\s^{(2)}_M}}\right)\right\}\,.
\Eq(4.lem3.15)
$$
If we now take
$
{\d^2_m}=4\bigl(\frac{\log m}{m}\bigr)^{2(1-\frac{a}{\a})}
$
and
$
M=\varepsilon\bigl(\frac{m}{\log m}\bigr)^{\frac{1}{\a}}
$,
$0<\varepsilon<1$,
then
$$
\eqalign{
&{\d^2_mm}/{\bar\s^{(2)}_M}
=
4\sfrac{(1-a)(2a-\a)}{\a}\left({1}/{\varepsilon}\right)^{2a-\a}\log m\,,
\cr
&{2\d_mM^a}/{\bar\s^{(2)}_M}
=4\sfrac{(1-a)(2a-\a)}{\a}\varepsilon^{\a-a}\,,
}
\Eq(4.lem3.16)
$$
(recall that by  assumption $2a>\a$ and $a<\a$).
Choosing $\varepsilon$ sufficiently small so as to guarantee that
$$
{\d^2_mm}/{\bar\s^{(2)}_M} \geq 6\text{and}{2\d_mM^a}/{\bar\s^{(2)}_M}\leq \log(4/3)\,,
\Eq(4.lem3.17)
$$
the bound \eqv(4.lem3.15) becomes
$
\bold P\left(A_m\right)\leq\frac{2}{m^2}\
$.
Thus $\sum_m \bold P\left(A_m\right)\leq\infty$ which,
invoking the first Borel-Cantelli Lemma, proves \eqv(4.lem3.6).

From now on we take
$
M=\varepsilon\bigl(\frac{m}{\log m}\bigr)^{{1}/{\a}}
$
and assume that $\varepsilon$ satisfies \eqv(4.lem3.17).
It remains to prove \eqv(4.lem3.7). Using that $\vp(x)\leq x^a$, $x\geq 0$, we have
$$
\frac{1}{m}\sum_{k=1}^{\infty}\vp(m^{1/\a}\g^+_k)
=\frac{1}{m}\sum_{k=1}^{\infty}\vp(m^{1/\a}\g_k)\1_{\{\g_k>M/m^{1/\a}\}}
\leq  m^{-(1-{a}/{\a})}\sum_{k=1}^{\infty}\g_k^a\1_{\{\g_k>{\varepsilon}/{(\log m)^{\frac{1}{\a}}}\}}\,.
\Eq(4.lem3.18)
$$
We further decompose the last sum in the r.h.s\. above into $\SS^-(m)+\SS^+(m)$, where
$$
\eqalign{
\SS^-(m)=&m^{-(1-{a}/{\a})}\sum_{k=1}^{\infty}\g_k^a
\1_{\{{\varepsilon}/{(\log m)^{{1}/{\a}}}<\g_k\leq 1\}}\,,
\cr
\SS^+(m)=&m^{-(1-{a}/{\a})}\sum_{k=1}^{\infty}\g_k^a\1_{\{\g_k>1\}}\,.
}
\Eq(4.lem3.19)
$$
To deal with $\SS^-(m)$ we write
$$
\eqalign{
\SS^-(m)
\leq
\frac{m^{a/\a}}{m}
\sum_{k=1}^{\infty}\1_{\{{\varepsilon}/{(\log m)^{{1}/{\a}}}<\g_k\leq 1\}}
=
\frac{m^{a/\a}{\mu((\log m)^{{1}/{\a}}, 1]}}{m}
\frac{\sum_{k=1}^{\infty}\1_{\{{\varepsilon}/{(\log m)^{{1}/{\a}}}<\g_k\leq 1\}}}
{\mu((\log m)^{{1}/{\a}}, 1]}
}\,.
\Eq(4.lem3.20)
$$
Since
$
\mu((\log m)^{{1}/{\a}}, 1]=\log m/\varepsilon^{\a}-1\uparrow\infty
$
as
$m\uparrow\infty$, it follows from the strong law of large numbers
for non-homogeneous Poisson processes (see \cite{Ki} p\. 51) that
$$
\lim_{m\rightarrow\infty}
\frac{\sum_{k=1}^{\infty}\1_{\{{\varepsilon}/{(\log m)^{{1}/{\a}}}<\g_k\leq 1\}}}
{\mu((\log m)^{{1}/{\a}}, 1]}
=1
\text{$\bold P$-almost surely.}
\Eq(4.lem3.21)
$$
and since
$
m^{a/\a-1}\mu((\log m)^{{1}/{\a}}, 1]= o(1)
$,
as follows from the assumption that $a<\a$,
we get that
$
\lim_{m\rightarrow\infty}\SS^-(m)
$
$\bold P$-almost surely.
To treat $\SS^+(m)$ note that
$
\int_{(0,\infty)}\min(u^a\1_{u>1}, 1)d\mu(u)<\infty
$.
Thus, by Campbell's Theorem,
$
\sum_{k=1}^{\infty}\g_k^a\1_{\{\g_k>1\}}<\infty
$
$\bold P$-almost surely. From this and the fact that
$
m^{{a}/{\a}-1}=o(1)
$,
we get that
$
\lim_{m\rightarrow\infty}\SS^+(m)
$
$\bold P$-almost surely. Collecting our results yields that
$
\lim_{m\rightarrow\infty}\frac{1}{m}\sum_{k=1}^{\infty}\vp(m^{1/\a}\g^+_k)=0
$
$\bold P$-almost surely, and establishes \eqv(4.lem3.7).

The proof of Lemma \thv(4.lemma3) is complete.\endproof

\bigskip


\chap{7. Proof of Theorems \thv(4.theo3), \thv(4.theo2), \thv(4.theo4) and Propositions \thv(4.prop1) and \thv(4.prop2), (ii).}7

In this section we prove the three theorems of Subsection 3.2 as well as
Proposition \thv(4.prop1), (ii), and Proposition \thv(4.prop2), (ii), of
Subsection 3.3.
This is where the results of Appendix A.2
on renewal theory in discrete and continuous time come into play.

\proofof{Theorem \thv(4.theo3)} Let $a<\a$ and set $\mu_n=\pi_n$.

\noindent{\it (i) Constant space scale.} Assume that $r_n$ is a constant space scale.
By Proposition \thv(4.prop1), (i),
and Corollary \thv(2.4.cor1), $\P$-almost surely,
$$
\lim_{n\rightarrow\infty}\CC_{n}(t,s)=\CC^{cst,-}_{\infty}(t,s)\quad\forall\, 0\leq t<t+s\,,
\EQ(8.theo3.1)
$$
where
$
\CC^{cst,-}_{\infty}(t,s)=\PP\left(\left\{R^{cst,-}(k)\,,k\in\N\right\}\cap (t, t+s)=\emptyset\right)
$,
and where $R^{cst,-}$ is the renewal process of inter-arrival distribution $\nu^{cst,-}$
defined in \eqv(4.prop1.2). By Lemma \thv(4.lemma1), $\nu^{cst,-}$ is regularly varying at
infinity with index $-\sfrac{\a-a}{1-a}$.
Since $a<\a$ and $0<\a<1$, $0<\sfrac{\a-a}{1-a}<1$.
Thus, by Dynkin-Lamperti Theorem in discrete time
\cite{Dyn, Lam} (see also Appendix A.2),
$
\lim_{t\rightarrow\infty}\CC^{cst,-}(t,\rho t)=\asl_{\sfrac{\a-a}{1-a}}(1/1+\rho)
$
for all $\rho>0$. Taking $s=\rho t$ in \eqv(8.theo3.1) and passing to the limit $t\rightarrow\infty$
yields the claim of Theorem \thv(4.theo3), (i).

\noindent{\it (ii) Intermediate space scale.} Assume that $r_n$ is an intermediate space scale. It follows
from  Corollary \thv(2.4.cor1) and Proposition \thv(4.prop2), (i), that, in $\P$-probability,
$$
\lim_{n\rightarrow\infty}\CC_{n}(t,s)=\CC^{int,-}_{\infty}(t,s)\quad\forall\, 0\leq t<t+s\,,
\EQ(8.theo3.2)
$$
where
$
\CC^{int,-}_{\infty}(t,s)=\PP\left(\left\{S^{int,-}(u)\,,u>0\right\}\cap (t, t+s)=\emptyset\right)
$.
Here $S^{int,-}$ is a stable subordinator of index $\sfrac{\a-a}{1-a}$. Thus, by
Dynkin-Lamperti Theorem in continuous time (see \eqv(A.2.theo2.3) of Theorem \eqv(A.2.theo2) in Appendix A.2),
for all $t\geq 0$ and all $\rho>0$,
$
\CC^{int,-}_{\infty}(t,\rho t)=\asl_{\sfrac{\a-a}{1-a}}(1/1+\rho)
$.
Taking $s=\rho t$ in \eqv(8.theo3.2) then yields \eqv(4.theo3.2).
The statement below  \eqv(4.theo3.2)
follows from the remark below Lemma \thv(6.lemma1). Theorem \thv(4.theo3), (ii), is proven.

\noindent{\it (iii) Extreme space scale.}
In this paragraph we use the representation of the landscape introduced in Subsection 6.1.
Assume that $r_n$ is an extreme space scale and consider the model obtained by substituting
the representation \eqv(7.1.3) for the original (ranked and re-scaled) landscape.
We will use bold letters to distinguish objects defined in this representation
from the original ones. Namely, we denote by $\bold S_n$ the pure
clock process \eqv(1.3.2), by $\wh{\bold S}_n$ the full
clock process \eqv(1.3.2'),
and by  $\bold C_{n}(t,s)$ the corresponding time correlation function \eqv(1.3.3).
Clearly, by Lemma \thv(7.lemma1),
$$
\CC_n(t,s)\overset{d}\to=\bold C_{n}(t,s)\text{for all $n\geq 1$ and all $0\leq t<t+s$.}
\EQ(8.theo3.5)
$$
Now, by Corollary \thv(2.4.cor1) and Proposition \thv(4.prop3), $\bold P$-almost surely,
$$
\lim_{n\rightarrow\infty}\bold C_{n}(t,s)=\CC^{ext,-}_{\infty}(t,s)\quad\forall\, 0\leq t<t+s\,,
\EQ(8.theo3.4)
$$
where
$
\CC^{ext,-}_{\infty}(t,s)=\PP\left(\left\{S^{ext,-}(u)\,,u>0\right\}\cap (t, t+s)=\emptyset\right)
$,
and  $S^{ext,-}$ is the (random) subordinator of (random) L\'evy measure $\nu^{ext,-}$ defined in \eqv(4.prop3.2).
Moreover, by Lemma \thv(4.lemma3),  $\nu^{ext,-}$ is  $\bold P$-almost surely regularly varying at
infinity with index $-\sfrac{\a-a}{1-a}$. Thus, by Dynkin-Lamperti Theorem in continuous time
(Theorem \eqv(A.2.theo2) of Appendix A.2),
applied for fixed $\o$ in the set of full measure for which Lemma \thv(4.lemma3) holds, we get that,
$\bold P$-almost surely,
$$
\lim_{t\rightarrow 0+}\CC^{ext,-}(t,\rho t)=\asl_{\sfrac{\a-a}{1-a}}(1/1+\rho)\quad\forall\,\rho>0\,.
\EQ(8.theo3.6)
$$
Finally, by \eqv(8.theo3.5) with  $s=\rho t$, using successively \eqv(8.theo3.4)
and \eqv(8.theo3.6) to pass to the limit $n\rightarrow\infty$ and $t\rightarrow 0+$,
we obtain that, for all $\rho>0$,
$
\lim_{t\rightarrow 0+}\lim_{n\rightarrow\infty}\CC_n(t,\rho t)\overset{d}\to=\asl_{\sfrac{\a-a}{1-a}}(1/1+\rho)
$.
Since convergence in distribution to a constant implies convergence in probability,
the claim of Theorem \thv(4.theo3), (iii) follows.
The proof of Theorem \thv(4.theo3) is complete.
\endproof

All the proofs stated in the remainder of this section are based on the approach
used in the proof of Theorem \thv(4.theo3), (iii), above:
we will first seek almost sure results for the model obtained by substituting the representation \eqv(7.1.3)
for the original landscape, and next transfer them to the original model using Lemma \thv(7.lemma1).

\proofof{Theorem \thv(4.theo4)} Assume that $r_n$ is an extreme space scale.
As in the proof of Theorem \thv(4.theo3), (iii), consider the model with landscape \eqv(7.1.3)
and denote by $\bold S_n$, $\wh{\bold S}_n$, and $\bold C_{n}(t,s)$
the corresponding pure clock process \eqv(1.3.2), full clock process \eqv(1.3.2'),
and time correlation function \eqv(1.3.3).
To keep the notation simple we do not introduce new symbols for the chains $X_n$ and $J_n$.
In particular their invariant measures, denoted as before by $\GG_{\a,n}$ and $\pi_n$, are the random
measures on $(\O^{E}, \FF^{E}, \bold P)$ on $\VV_n$ defined through:
$$
\eqalign{
\GG_{\a,n}(k)&=\frac{\g_{nk}}{\sum_{k=1}^n\g_{nk}}\,,\quad k\in\VV_n\,,\cr
\pi_n(k)&=\frac{\g_{nk}^a}{\sum_{k=1}^n\g_{nk}^a}\,,\quad k\in\VV_n\,.
}
\Eq(4.theo4'.0)
$$
The proof of Theorem \thv(4.theo4) makes use of the following statement and its proof.

\theo{\TH(4.theo4')} {\it   Let $r_n$ be an extreme space scale.
The following holds for all $0\leq a<1$:

\item{(i)} If  $\mu_n=\GG_{\a,n}$ then, $\bold P$-almost surely, for all $\leq t<t+s$,
$$
\lim_{n\rightarrow\infty}\bold C_{n}(t,s)={\CC}^{sta}_{\infty}(s)\,.
\Eq(4.theo4'.00)
$$

\item{(ii)} If $\mu_n=\pi_n$, $\bold P$-almost surely, for all $s>0$,
$$
\lim_{t\rightarrow\infty}\lim_{n\rightarrow\infty}\bold C_{n}(t,s)
={\CC}^{sta}_{\infty}(s)\,.
\Eq(4.theo4'.000)
$$
}

The proof of Theorem \thv(4.theo4') relies on results from renewal theory
in the ``classical'' setting referred to, in appendix A.2, as the ``finite mean life time'' case
(see Theorem \thv(A.2.theo2), (ii), and Theorem \thv(A.2.theo2'), (ii), of appendix A.2.1,
and Theorem \thv(A.2.theo2") of appendix A.2.2).

\proofof{Theorem \thv(4.theo4')} We first prove assertion (ii).
By Corollary \thv(2.4.cor1) and Proposition \thv(4.prop3), $\bold P$-almost surely,
$$
\lim_{n\rightarrow\infty}\bold C_{n}(t,s)=\CC^{ext,\pm}_{\infty}(t,s)\quad\forall\, 0\leq t<t+s\,,
\EQ(8.theo4'.1)
$$
where
$$
\eqalign{
\CC^{ext,-}_{\infty}(t,s)&=\PP\left(\left\{S^{ext,-}(u)\,,u>0\right\}\cap (t, t+s)=\emptyset\right)\text{if $a<\a$,}\cr
\CC^{ext,+}_{\infty}(t,s)&=\PP\left(\left\{R^{ext,+}(u)\,,u>0\right\}\cap (t, t+s)=\emptyset\right)\text{if $a>\a$,}\cr
}
\EQ(8.theo4'.2)
$$
and where $S^{ext,-}$ is the subordinator of L\'evy measure $\nu^{ext,-}$,
and $R^{ext,+}$ is the renewal process of inter-arrival distribution $\nu^{ext,+}$,
$\nu^{ext,\pm}$ being defined in \eqv(4.prop3.2).
In view of Lemma \thv(4.lemma2) we are now in the classical setting
``finite mean life time'' renewal theory so that the claim of assertion (ii) follows
from \eqv(8.theo4'.2), using  Dynkin-Lamperti Theorem for ``finite mean life time''.
More precisely, if $a<\a$ then, by Lemma \thv(4.lemma2) and Theorem \thv(A.2.theo2), (ii),
we have that, $\bold P$-almost surely, for each fixed $s>0$,
$$
\lim_{t\rightarrow\infty}\CC^{ext,-}_{\infty}(t,s)
=\frac{1}{m^-}\int_s^{\infty}\nu^{ext,-}(x,\infty)dx
={\CC}^{sta}_{\infty}(s)
\,,
\EQ(8.theo4'.3)
$$
where ${\CC}^{sta}_{\infty}$ is defined as in \eqv(4.theo4.0).
Similarly,
if $a>\a$ then, by Lemma \thv(4.lemma2) and Theorem \thv(A.2.theo2'), (ii), we have that,
$\bold P$-almost surely, for each fixed $s>0$,
$$
\lim_{t\rightarrow\infty}\CC^{ext,+}_{\infty}(t,s)
=\frac{1}{m^+}\int_s^{\infty}\nu^{ext,+}(x,\infty)dx
={\CC}^{sta}_{\infty}(s)
\,.
\EQ(8.theo4'.4)
$$
Assertion (ii) of the theorem is thus proven.

To prove assertion (i) we first need to check Condition (A0) (see \eqv(1.A0)) when  $\mu_n=\GG_{\a,n}$.
By \eqv(4.theo4'.0),
$$
1-F_n(v):=\sum_{x\in\VV_n}\GG_{\a,n}(x)e^{-vc_n\l_n(x)}
=\sum_{k}\frac{\g_{nk}}{\sum_{l}\g_{nl}}e^{-s\g_{nk}^{-(1-a)}}\,.
$$
A straightforward application of Proposition \thv(7.prop1) then yields that, for all $0\leq a<1$,
$\bold P$-almost surely,
$
\lim_{n\rightarrow\infty}(1-F_n(v))=(1-F^{sta}(v)):={\CC}^{sta}_{\infty}(s)
$.
Hence, by Theorem \thv(2.4.theo2), $\bold P$-almost surely,
denoting by $\s^{sta}$ the random variable with distribution function $F^{sta}$,
$$
\eqalign{
&\wh{\bold S}_n(\cdot)\Rightarrow \wh S^{ext,-}(\cdot)=\s^{sta}+S^{ext,-}(\cdot)\text{if $a<\a$,}\cr
&\wh{\bold S}_n(\cdot)\Rrightarrow \wh R^{ext,+}(\cdot)=\s^{sta}+R^{ext,+}(\cdot)\text{if $a>\a$,}
}
\EQ(8.theo4'.5)
$$
where in the first line (respec\. the second line) $\s^{sta}$
is independent of $S^{ext,-}$ (respec\. $R^{ext,+}$).
In view of \eqv(8.theo4'.5) and the just proven assertion (ii) of Theorem \thv(4.theo4'),
we are now in the realm of stationary processes, and the conclusion will follow
from Theorem \thv(A.2.theo2") of Appendix A.2.2.

More precisely, let $\o\in\wt\O^E$ be fixed, where $\wt\O^E$ denotes the set of full measure for which \eqv(8.theo4'.5) obtains.
If $a<\a$, by \eqv(8.theo4'.3) and the definition of $F^{sta}$,
$
F^{sta}(s)=\lim_{t\rightarrow\infty}\CC^{ext,-}_{\infty}(t,s)
$.
By Theorem \thv(A.2.theo2"), (ii), and the first line of \eqv(8.theo4'.5) it then follows that
$
\wh S^{ext,-}\overset{d}\to=S^{ext,-}
$,
so that
$
\PP\bigl(\theta_t(\wh S^{ext,-})\geq s\bigr)=\PP\bigl(\theta_t(S^{ext,-})\geq s\bigr)
$.
Recalling from the proof of Theorem \thv(2.4.theo2) (see \eqv(2.2.22))
that \eqv(2.4.theo2.3) can equivalently be written as
$
\lim_{n\rightarrow\infty}{\CC}_n(t,s)=\PP\bigl(\theta_t(\wh S)\geq s\bigr)
$,
we finally get that, for all $\leq t<t+s$,
$
\lim_{n\rightarrow\infty}\bold C_{n}(t,s)=\PP\bigl(\theta_t(S^{ext,-})\geq s\bigr)={\CC}^{sta}_{\infty}(s)
$.
Since this holds true for all $\o\in\wt\O^E$, the claim of assertion (i) follows.
If $a>\a$ it similarly follows from \eqv(8.theo4'.4), Theorem \thv(A.2.theo2"), (i),
and the second line of \eqv(8.theo4'.5) that, $\bold P$-almost surely, for all $\leq t<t+s$,
$
\lim_{n\rightarrow\infty}\bold C_{n}(t,s)=\PP\bigl(\theta_t(R^{ext,+})\geq s\bigr)={\CC}^{sta}_{\infty}(s)
$.
This concludes the prove of assertion (i). The proof of Theorem \thv(4.theo4') is complete.\endproof

We may now conclude the proof of Theorem \thv(4.theo4). The first assertion directly
follows from \eqv(8.theo3.5) and Theorem \thv(4.theo4'), (i).  The second assertion follows from
from \eqv(8.theo3.5) with  $s=\rho t$, using successively \eqv(8.theo4'.1)
and \eqv(8.theo4'.3) to pass to the limit $n\rightarrow\infty$ and $t\rightarrow 0+$
(just as in the proof of Theorem \thv(4.theo3), (iii)).
The proof of Theorem \thv(4.theo4) is done.\endproof

It remains to prove Theorem \thv(4.theo2), and
assertion (ii) of Proposition \thv(4.prop1) and of Proposition \thv(4.prop2).

\proofof{Proposition \thv(4.prop1) and \thv(4.prop2), (ii)}
In the sequel we will use the symbol $\bar r_n$ to denote a constant or intermediate space scale
and keep the symbol $r_n$ for extreme scales.
Let $a>\a$ and assume that $c_n=\bar r_n^{(1-a)}$ where $\bar r_n$ is a constant
or intermediate space scale. Proceeding as in \eqv(4.prop3.5)
to express  $\nu_n(u,\infty)$ in the landscape representation
\eqv(7.1.3),
and denoting  by $\bold v_n(u,\infty)$ the resulting quantity,
we get, choosing $a_n=1$ and setting $\rho_n=\bar r_n/r_n$,
$$
\bold v_n(u,\infty)=
\frac{
\sum_{k=1}^n\g_{nk}^ae^{-u\rho_n/\g_{nk}^{(1-a)}}
}
{
\sum_{k=1}^n\g_{nk}^a
}\,,\quad u\geq 0\,.
\Eq(4.prop1-2.1)
$$
By Lemma \thv(4.lemma4), for all $\varepsilon>0$ and all $n$ large enough,
$0\leq\rho_n\leq \varepsilon$. Thus, for all $\varepsilon>0$ and all $n$ large enough,
$$
1
\geq
\bold v_n(u,\infty)
\geq
\frac{\sum_{k=1}^n\g_{nk}^ae^{-\varepsilon u/\g_{nk}^{(1-a)}}}{\sum_{k=1}^n\g_{nk}^a}
\,,\quad u\geq 0\,.
\Eq(4.prop1-2.4)
$$
Note that the lower bound of \eqv(4.prop1-2.4) is nothing but \eqv(4.prop3.5)
evaluated at $\varepsilon u$. Using \eqv(4.prop3.10) to pass to the limit
$n\rightarrow\infty$ in \eqv(4.prop1-2.1) yields that, $\bold P$-almost surely,
for all $\varepsilon>0$,
$$
1
\geq
\lim_{n\rightarrow\infty}\bold v_n(u,\infty)
\geq
\nu^{ext,+} (\varepsilon u,\infty)
\,,\quad u\geq 0\,,
\Eq(4.prop1-2.5)
$$
where $\nu^{ext,+}$ is defined in \eqv(4.prop3.2). Since $\nu^{ext,+} (0,\infty)=1$, passing to the limit
$\varepsilon\rightarrow 0$ in \eqv(4.prop1-2.5) finally yields that, $\bold P$-almost surely,
$
\lim_{n\rightarrow\infty}\bold v_n(u,\infty)=1
$,
$u\geq 0$. By Lemma \thv(7.lemma1), for each $n\geq 1$,
$
\bold v_n(u,\infty)\overset{d}\to=\nu_n(u,\infty)
$.
Therefore, for all $u\geq 0$,
$$
\lim_{n\rightarrow\infty}\nu_n(u,\infty)=1\text{in $\P$-probability.}
$$
Using  the monotonicity of $\nu_n$ 
it readily follows from a subsequence argument, that, for all $\e>0$,
$
\lim_{n\rightarrow\infty}\P\left(
\sup_{u>0}\left|\nu_{n}(u,\infty)-1\right|>\e
\right)=0
$.
Therefore, by Theorem \thv(2.4.theo1), (ii), $S_n(\cdot)\Rrightarrow R^{*,-}(\cdot)$
in $\P$-probability, where $R^{*,-}$ is the degenerate renewal process of inter-arrival distribution
where $\nu^{*,-}=\d_{\infty}$.
Assertion (ii) of Proposition \thv(4.prop1) and of Proposition \thv(4.prop2) are thus proven.
\endproof


\proofof{Theorem \thv(4.theo2)} 
Assume first that $r_n$ is an extreme space scale.
The starting point of the proof of assertion (ii) is \eqv(8.theo4'.1) for $a>\a$.
By the second line of \eqv(8.theo4'.2), conditioning on the first jump of $R^{ext,+}$ yields
$
\CC^{ext,+}_{\infty}(t,s)=1-F(t+s)+\int_{0}^t\CC^{ext,+}_{\infty}(t-v,s)dF(v)
$,
where $F(u)=1-\nu^{ext,+}(u)$ (see e.g\. \eqv(2.2.22)).
On the one hand this implies that
$
1\geq \CC^{ext,+}_{\infty}(t,\rho t)\geq 1-F(t+\rho t)=\nu^{ext,+}(t+\rho t)
$.
On the other hand it easily follows from the definition of $\nu^{ext,+}(u)$ (see \eqv(4.prop3.2))
that $\lim_{t\rightarrow 0}\nu^{ext,+}(t)=1$ $\bold P$-almost surely.
Therefore
$\lim_{t\rightarrow 0}\CC^{ext,+}_{\infty}(t,\rho t)=1$ $\bold P$-almost surely.
Combining this statement with \eqv(8.theo3.5) yields the claim of assertion (ii).

Assume now that $r_n$ is a constant or intermediate space scale and let $a>\a$.
It follows from  Corollary \thv(2.4.cor1) and
either Proposition \thv(4.prop1), (ii), for constant scales,
or Proposition \thv(4.prop2), (ii), for intermediate scales that,
in $\P$-probability,
$$
\lim_{n\rightarrow\infty}\CC_{n}(t,s)=\CC^{*,-}_{\infty}(t,s)\quad\forall\, 0\leq t<t+s\,,
\EQ(8.theo3.2')
$$
where
$
\CC^{*,-}_{\infty}(t,s)=\PP\left(\left\{S^{*,-}(u)\,,u>0\right\}\cap (t, t+s)=\emptyset\right)
$,
and where $R^{*,-}$ is the degenerate renewal process of inter-arrival distribution
where $\nu^{*,-}=\d_{\infty}$.
Conditioning on the first jump and arguing as above readily yields that
$\CC^{*,-}_{\infty}(t,s)=1$ for all $0\leq t<t+s$.
Inserting this result in \eqv(8.theo3.2') proves assertion (i).
The proof of Theorem \thv(4.theo2) is done. \endproof

\bigskip


\def\drift{\roman{\bf d}}

\chap{A. Appendix}{9}

\bigskip
\line{\bf A.1. Subordinators and renewal processes. \hfill}

We summarize here succinctly the
needed information about subordinators and renewal processes.
Classical references are
the book by Bertoin [Be]
and  It\^o's lecture notes  \cite{I} (for subordinators)
 and Feller \cite{Fe} and Bingham et {\it al\.} \cite{BGT} (for renewal processes).

\noindent{\bf Definition.} Subordinators form the sub-class of L\'evy processes
(processes with stationary independent increments) that take values in $[0,\infty)$.
Let $\{S(t), t\geq 0\}$ be a subordinator. Its Laplace transform takes the characteristic form
$$
E(\exp\{-\theta S(t)\})=\exp\{-t\Phi(\theta)\}\,,\quad \theta\geq 0\,,
\Eq(A.2.1)
$$
where $\Phi : [0,\infty)\rightarrow [0,\infty)$ (called the Laplace exponent) is given by
$$
\Phi(\theta)=\drift+\int_{(0,\infty)}\left(1-e^{-\theta x}\right)\nu(dx)\,,
\Eq(A.2.2)
$$
and where $\drift\in\R$ is a constant drift term and $\nu$ (called the L\'evy measure) is a
$\s$-finite measure on $(0,\infty)$ with the property that
$\int_{(0,\infty)}(1\wedge x)\nu(dx)<\infty$.

Stable subordinators with index $\a$ are the important sub-class of subordinators whose
L\'evy measure has the form $\nu(x,\infty)=c x^{-\a}$ for some $0<\a<1$ and $c>0$.

\noindent{\bf The It\^o representation.} The following result due to It\^o (c.f\. \cite{I} page 1.11.2)
establishes the relation between subordinators and an associated Poisson random measure.

\theo{\TH(A.2.theo1)}{\it $S(t)$ can be represented as
$$
S(t)=\drift t+\int_{0<s\leq t}\int_{0<u<\infty}uN(ds\,,du)\,,
\Eq(A.2.theo1.1)
$$
where $N$ is a Poisson random measure on $[0,\infty)\times(0,\infty)$ with intensity measure $dt\times d\nu$.
}

\noindent Let $\{(t_k,\xi_k)\}$ denote the points of $N$: they represent the pairs of jump times and jump size of $S(t)$.
Using them we may write \eqv(A.2.theo1.1) in the nice alternative form
$$
S(t)=\drift t+
\sum_{t_k\leq t}\xi_k\,.
\Eq(A.2.theo1.2)
$$

As the next definition shows, renewal processes can be thought of as a subordinator sampled at equidistant points.

\noindent{\bf Definition.} A renewal process $\{R(n)\,,n\in\N\}$ is a partial sum process
with identical and independent increments taking values in $[0,\infty)$. $R(n)$ is represented as
$$
R(n)=\sum_{k\leq n}\xi_k\,,
$$
where $\{\xi_k, k\geq 1\}$ are independent r\.v\.'s with identical
distribution $\nu$.
The $\xi_k$'s, which stand for the life-time of items, are called {\it inter-arrival times}; their law, $\nu$,
is called the {\it inter-arrival distribution}.

\noindent{\bf Delayed processes.}
a {\it delayed} renewal process corresponding to a renewal process $R$ is the process $\wh R$ defined by $\wh R=\s +R$
where $\s$ is a nonnegative random variable independent of $S$.
Similarly the {\it delayed} subordinator corresponding to a subordinator $S$ is the process $\wh S$ defined by $\wh S=\s +S$
where $\s$ is a nonnegative random variable independent of $S$.
We will say that a renewal process or subordinator {\it pure} when we want to emphasize that $\s=0$.

\bigskip
\line{\bf A.2. Renewal theory. \hfill}

We now summarize what we need to know about renewal theory for subordinators and renewal processes.

\bigskip
\line{\bf A.2.1. The Dynkin-Lamperti Theorem. \hfill}


We will refer to the theorem below as to Dynkin-Lamperti Theorem in continuous time.

Set
$$
\CC_{\infty}(t,s)=\PP\left(\left\{S(u)\,,u>0\right\}\cap (t, t+s)=\emptyset\right)\,,\quad 0\leq t<t+s\,.
\Eq(A.2.theo2.0)
$$
where $S$ is a subordinator of L\'evy measure $\nu$.
As already observed in \eqv(2.2.19), \eqv(A.2.theo2.0)  can be written in the more classical form
$
\CC_{\infty}(t,s)
=
\PP\left(\theta_t(S)\geq s\right)
$,
where $\theta_t(\eta)$ is the overshoot function defined in \eqv(2.2.17).

\theo{\TH(A.2.theo2)}{\it
\item{(i)}[Arcsine law.] If the tail of the L\'evy measure $\nu(x,\infty)$
is regularly varying at infinity with index $\a\in [0,1]$, then
$$
\lim_{t\rightarrow\infty}\CC_{\infty}(t,\rho t)=
\cases
\asl_{\a}(1/1+\rho),
&\hbox{if  $\,\,0<\a<1$} , \,\,\,\cr
1 , &\hbox{if $\,\,\a=0$} , \,\,\,\cr
0 , &\hbox{if $\,\,\a=1$.}\,\,\,\cr
\endcases
\Eq(A.2.theo2.1)
$$
If $\nu(x,\infty)$ is regularly varying at $0^+$ with index $\a\in [0,1]$, then
\eqv(A.2.theo2.1) holds with $t\rightarrow\infty$ replaced by $t\rightarrow 0+$.
\item{} If $\nu(x,\infty)=\kappa x^{-\a}$ for some constant $\kappa>0$ and $\a\in(0,1)$ (that is if
$S$ is a stable subordinator with index $\a\in(0,1)$) then
$$
C_{\infty}(t,\rho t)=\asl_{\a}(1/1+\rho)\text{for all} t>0\,.
\Eq(A.2.theo2.3)
$$
Moreover, in order for $\CC_{\infty}(t,\rho t)$ to converge to the integrated arcsine density \eqv(A.2.theo2.1)
it is necessary and sufficient that $\nu(x,\infty)$ be regularly varying with index $\a\in (0,1)$.

\item{(ii)}[Finite mean life time renewal.] If
$
\int_0^{\infty}\nu(x,\infty)dx=m<\infty
$
and
$S(\cdot)$ is not a compound Poisson process (i\.e\. $S(\cdot)$ is non-arithmetic)
then, for each fixed $s>0$,
$$
\lim_{t\rightarrow\infty}\CC_{\infty}(t,s)=\frac{1}{m}\int_s^{\infty}\nu(x,\infty)dx\,.
\Eq(A.2.theo2.4)
$$
}

\proof
The first half of the theorem  -- namely the arcsine law -- is a
restatement of Theorem 6,
assertion (iii), page 81 of \cite{Ber}. The second half -- the finite mean life time case --
is contained in Theorem 1 of \cite{BvHS}.\endproof

We now state the ``classical'' discrete time Dynkin-Lamperti Theorem.
Set
$$
\CC_{\infty}(t,s)=\PP\left(\left\{R(k)\,,k\in\N\right\}\cap (t, t+s)=\emptyset\right)\,,\quad 0\leq t<t+s\,,
\Eq(A.2.theo2'.0)
$$
where $R$ is a renewal process of inter-arrival distribution $\nu$.
Let $\theta_t(\cdot)$ denote the overshoot function \eqv(2.2.17) in discrete time.
In this setting $\theta_t(R)$ is usually called the residual waiting time. Clearly,
$
\CC_{\infty}(t,s)=\PP\left(\theta_t(R)\geq s\right)
$.
One has (see \cite{Dyn} or \cite{BGT}, section 8.6):

\theo{\TH(A.2.theo2')} [Dynkin, 55(61) and Lamperti, 58]. {\it
\item{(i)}[Arcsine law.] A necessary an sufficient condition
for $\theta_t(R)/t$ to have a non-degenerate limit law is
that $\nu$ is regularly varying at infinity with index $0<\a<1$. In that case,
$$
\lim_{t\rightarrow\infty}\CC_{\infty}(t,\rho t)=\asl_{\a}(1/1+\rho)\,.
\Eq(A.2.theo2'.1)
$$
\item{(ii)}[Finite mean life time renewal.] If
$
\int_0^{\infty}\nu(x,\infty)dx=m<\infty
$
and if $\nu$ is non-latticed then, for each fixed $s>0$,
$$
\lim_{t\rightarrow\infty}\CC_{\infty}(t,s)=\frac{1}{m}\int_s^{\infty}\nu(x,\infty)dx\,.
\Eq(A.2.theo2'.4)
$$
}

\bigskip
\line{\bf A.2.2. Stationarity of delayed processes with ``finite mean life time''. \hfill}

It is well known (see e.g\. \cite{Fe}), and not difficult to prove, that when $\int_0^{\infty}\nu(x,\infty)dx=m<\infty$,
the delayed renewal process $\wh R=\s +R$, whose initial jump is sampled from
the limit law of the residual waiting time $\theta_t(R)$,
is stationary.
A similar statement holds for the delayed subordinator $\wh S:=\s +S$ (see \cite{vHS}). These results are summarized
in the theorem below.

\theo{\TH(A.2.theo2")} {\it Let $F$ denote the distribution function of $\s$.
\item{(i)}[Delayed renewal process] Under the assumptions and with the notations of Theorem \thv(A.2.theo2'), (ii),
if
$
F(s)=\lim_{t\rightarrow\infty}\CC_{\infty}(t,s)
$,
then, denoting by $\wh R$ the delayed renewal process $\wh R=\s +R$,
$$
\wh R\overset{d}\to= R\,.
$$
\item{(ii)}[Delayed subordinator] Under the assumptions and with the notations of Theorem \thv(A.2.theo2), (ii),
if
$
F(s)=\lim_{t\rightarrow\infty}\CC_{\infty}(t,s)
$,
then, denoting by $\wh S$ the delayed subordinator $\wh S=\s +S$,
$$
\wh S\overset{d}\to= S\,.
$$
}


\line{\bf A.3. Regular variations. \hfill}

We assume as known the elementary properties of regularly and slowly varying functions as described
in Section 1 of \cite{BGT} and, in particular, the {\it Uniform Convergence Theorem}
(\cite{BGT}, Theorem 1.2.1)
for slowly varying functions (\cite{BGT}, Theorem 1.3.1).
In the sequel we denote by $R_0$ the class of functions that are slowly varying at infinity, by $R_{\rho}$ the class
of functions that are regularly varying at infinity with index $\rho$, by $R_{\rho}(0+)$
the class of functions that are regularly varying at $0+$, and we set $R=\cup_{\rho\in\R}R_{\rho}$
(\cite{BGT}, Section 1.4.2).
The results below are stated in the setting of slow variations at infinity. They can easily be adapted to
that of slow variations at the origin by using that
a function $f(x)$ is  slowly (regularly) varying at zero if and only if $f(x^{-1})$ is  slowly (regularly)
varying at infinity.
The next two lemmas contain bounds on slowly varying functions  that
will often be needed in Section 5 and 6.

\lemma{\TH(A.3.lemma1)}[\cite{Fe},VIII.8, Lemma 2.] {\it If $\ell\in R_0$ then $x^{-\e}\leq \ell(x)\leq x^{\e}$  for any fixed $\e>0$
and all $x$ sufficiently large.}

We will also frequently use the following bounds of Potter's type.

\lemma{\TH(A.3.lemma2)}{\it Let $\ell\in R_0$ and let $u_n$ and $v_n$ be positive non decreasing sequences such that
$v_n\uparrow\infty$, $u_n\uparrow\infty$ as $n\uparrow\infty$. For any given $x>0$ there exist positive sequences
$\e_n$ and $\d_n$ that verify $\e_n\downarrow 0$, $\d_n\downarrow 0$ as $n\uparrow\infty$
and such that, for all $n$ large enough,
$$
(1-\d_n)\min\left\{\Bigl(\sfrac{v_n}{u_n}x\Bigr)^{\e_n},\Bigl(\sfrac{v_n}{u_n}x\Bigr)^{-\e_n}\right\}\leq
\frac{\ell(v_n x)}{\ell(u_n)}
\leq (1+\d_n)\max\left\{\Bigl(\sfrac{v_n}{u_n}x\Bigr)^{\e_n},\Bigl(\sfrac{v_n}{u_n}x\Bigr)^{-\e_n}\right\}\,.
\Eq(A.3.2)
$$
}


Both these lemmata are  immediate consequences of the {\it Representation Theorem} for slowly
varying functions which we now state.

\theo{\TH(A.3.theo1)}[\cite{BGT}, I.3.1, Theorem 1.3.1.] {\it The function $\ell$ is slowly varying at
infinity if and only if it may be written in the form
$$
\ell(x)=\kappa(x)\exp\left\{\int_{a}^{x}\frac{\varepsilon(y)}{y}dy\right\}\quad(x\geq a)\,,
\Eq(A.3.3)
$$
for some $a>0$, where $\kappa(\cdot)$ is measurable  and $\kappa(x)\rightarrow \kappa\in(0,\infty)$,
$\varepsilon(x)\rightarrow 0$ as $x\rightarrow\infty$.
}

Finally we state an important result about inverse of regularly varying functions. Let $f$ be a function
defined and locally bounded on $[0,\infty)$, and that tends to zero as $x\rightarrow\infty$. Its
{\it generalized inverse}
$$
f^{-1}(x):=\inf\{y\geq 0 : f(y)\leq x\}\,,
\Eq(A.3.1)
$$
is defined on $[f(0),\infty)$. The following result is an (easy) adaptation to the case of functions
$f$ in $R_{\rho}$ with $\rho<0$ of a theorem of \cite{BGT} stated for $\rho>0$.

\lemma{\TH(A.3.lemma3)}[\cite{BGT}, I.5.7, Theorem 1.5.12.]  {\it If $f\in R_{\rho}$ with $\rho<0$,
there exists $g\in R_{1/\rho}(0+)$ with
$$
f(g(x))\sim g(f(x))\sim x\,,\quad x\rightarrow 0\,.
\Eq(A.3.lemma3.1)
$$
Here $g$ (an `asymptotic inverse' of $f$) is determined to within asymptotic equivalence,
and one version of $g$ is $f^{-1}$.
}


\bigskip
\line{\bf A.4. A technical lemma. \hfill}

We conclude this appendix with a technical lemma that will be needed
to prove that Conditions (A1)-(A3) are verified uniformly in $\d$, when convergence
holds in $\P$-probability only.

\lemma{\TH(app.A.4)}{\it Let $\{X_n, n\geq  1\}$, be a sequence of random
variables defined on $(\O^{\t}, \FF^{\t}, \P)$, taking values in the space of positive
decreasing functions on $(0,\infty]$. Assume that there exist decreasing sequences
$\rho_n$ and $\eta_n$ satisfying $0<\rho_n,\eta_n\downarrow 0$ as $n\uparrow\infty$,
and positive decreasing functions
$f_n$ and $g_n$ on $(0,\infty]$ such that, for all $u>0$,
$$
\P\Bigl(
\left|X_n(u)-f_n(u)\right|\geq \eta_ng_n(u)
\Bigr)\leq\rho_n\,.
\Eq(app.A.4.1)
$$
Assume in addition that, for all large enough $n$, there exist constants $0<\kappa, \kappa'<\infty$
and an integer $l_0$ such that, for all $l\geq l_0$,
$$
g_n\Bigl(\frac1{l}\Bigr)\geq \kappa g_n\Bigl(\frac1{l+1}\Bigr)
\text{and}
g_n(l+1)\geq \kappa' g_n(l)\,.
\Eq(app.A.4.2)
$$
Then
$$
\lim_{n\rightarrow 0}\P\left(
\sup_{u>0}\bigl\{\left|X_n(u)-f_n(u)\right|\geq \eta_ng_n(u)\bigr\}\right)
=0
\,.
\Eq(app.A.4.3)
$$
The conclusions of the lemma are unchanged if $f_n$ and $g_n$ positive decreasing functions, that satisfy the relations
\eqv(app.A.4.2) with reversed inequalities.
}

\proof Given a constant $0<c<\infty$
set
$
A_{n,c}(u)=\left\{
c^{-1}\left|X_n(u)-f_n(u)\right|\geq \eta_ng_n(u)
\right\}
$ and write
$
A_{n,1}(u)\equiv A_{n}(u)
$.
Let us first establish that
$$
\textstyle
\lim_{n\rightarrow 0}\P\left(
\sup_{0<u\leq 1/l_0}
A_{n}(u)
\right)
=0
\,.
\Eq(app.A.4.4)
$$
Under the assumptions of the lemma,
$
\bigcup_{\frac{1}{l+1}<u\leq \frac{1}{l}}A_n(u)
\subseteq
A_n\left(\sfrac{1}{l}\right)\cup A_{n,\kappa}\left(\sfrac{1}{l+1}\right)
$,
$l\geq l_0$.
Hence
$$
\textstyle
\P\left(\bigcup_{0<u\leq 1}A_n(u)\right)
=\P\left(\bigcup_{l=l_0}^{\infty}\bigcup_{\frac{1}{l+1}<u\leq \frac{1}{l}}A_n(u)\right)
\leq
\P\left(\bigcup_{l=l_0}^{\infty}A_n\left(\sfrac{1}{l}\right)\right)
+\P\left(\bigcup_{l=l_0}^{\infty} A_{n,\kappa}\left(\sfrac{1}{l+1}\right)\right)
\,.
\Eq(app.A.4.6)
$$
Consider the first term in the right hand side of \eqv(app.A.4.6).
Set
$
Y_n(l)=g^{-1}_n\left(\sfrac{1}{l}\right)\left|X_n\left(\sfrac{1}{l}\right)-f_n\left(\sfrac{1}{l}\right)\right|
%
$,
$l\geq l_0$, and
$
Z_n=\sup_{l\geq l_0}Y_n(l)
$.
Then
$
\P\left(\bigcup_{l=l_0}^{\infty}A_n\left(\frac{1}{l}\right)\right)
=\P\left(Z_n\geq \eta_n\right)
$,
whereas the assumption \eqv(app.A.4.1) becomes
$
\P(Y_n(l)\geq \eta_n)\leq\rho_n
$.
Using a classical subsequence argument, one readily deduces from the latter bound that
$
\P\left(Z_n\geq \eta_n\right)=o(1)
$ as
$n\rightarrow\infty$.

Indeed,
for each $k\geq 1$, choose $m_k$ to be the smallest integer such that $\rho_{m_k}\leq 2^{-k}$.
Without loss of generality we may assume that $m_k$ is a strictly increasing sequence.
Then, for all $n\geq m_k$,
$
\P(Y_n(l)\geq \eta_n)\leq 2^{-k}
$.
In particular, for any subsequence $n_k$ satisfying $m_k\leq n_k<m_{k+1}$, $k\geq 1$, we have
$
\P(Y_{n_k}(l)\geq \eta_{n_k})\leq 2^{-k}
$.
By Borel-Cantelli Lemma,
$
\P\left(\bigcup_{k=1}^{\infty}\bigcap_{m=k}^{\infty}\left\{Y_{n_m}(l)< \eta_{n_m}\right\}\right)=1
$.
Hence
$
\P\left(\bigcup_{k=1}^{\infty}\bigcap_{m=k}^{\infty}\bigcap_{l=l_0}^{\infty}\left\{Y_{n_m}(l)< \eta_{n_m}\right\}\right)=1
$.
Equivalently,
$
\P\left(\bigcap_{k=1}^{\infty}\bigcup_{m=k}^{\infty}\left\{Z_{n_{m}}\geq \eta_{n_m}\right\}\right)=0
$,
that is,
$
\lim_{k\rightarrow\infty}\P\left(\bigcup_{m=k}^{\infty}\left\{Z_{n_{m}}\geq \eta_{n_m}\right\}\right)=0
$.
We thus established that for each $n_k$ such that $m_k\leq n_k\leq m_{k+1}$, $\P\bigl(Z_{n_k}\geq \eta_{n_k}\bigr)=o(1)$,
$k\rightarrow\infty$. From this the desired conclusion follows.
}

Arguing in the same way to deal with the last term in the right hand side of \eqv(app.A.4.6)
yields the claim of \eqv(app.A.4.4). One proves in exactly the same way that
$
\lim_{n\rightarrow 0}\P\bigl(\sup_{u> 1/l_0}A_{n}(u)\bigr)=0
$
using this time that
$
\textstyle
\bigcup_{l\leq u <l+1}A_n(u)
\subseteq
A_{n,\kappa'}(l)\cup A_n(l+1)
$,
$l\geq l_0$. The proof of the lemma is done.\endproof


\bigskip


\end

\ref
\key FM
\by L.R.G\. Fontes and P\. Mathieu
\paper $K$-processes, scaling limit and aging for the trap model in the complete graph
\jour Ann. Probab.
\vol 4
\issue 36
\pages 1322--1358
\yr 2008
\endref

C(erný, Jir(í; Gayrard, Véronique Hitting time of large subsets of the hypercube.
Random Structures Algorithms  33  (2008),  no. 2, 252--267. (Reviewer: George Stoica)

\ref
\key C
\by J\. \v Cern\'y
\paper The behaviour of aging functions in one-dimensional Bouchaud's trap model
\jour Comm. Math. Phys.
\vol 261
\issue  1
\pages 195--224
\yr 2006
\endref

\ref BB
\key
\by J\.-P\. Bouchaud and E\. Bertin
\paper Dynamical ultrametricity in the critical trap model.
\jour J. Phys. A
\vol 13
\issue 35
\pages 3039--3051
\yr 2002
\endref

\bigskip

\end